\documentclass[11pt]{amsart}

\usepackage{hyperref}
\usepackage{amsthm} 
\usepackage{amsmath} 
\usepackage{amssymb} 
\usepackage{pinlabel}
\usepackage{cleveref}

\usepackage{thm-restate}

\usepackage{tikz}
\usepackage{enumerate}
\usepackage[numbered]{bookmark} 

\newcommand{\calA}{\mathcal{A}}
\newcommand{\calB}{\mathcal{B}}
\newcommand{\calC}{\mathcal{C}}
\newcommand{\calD}{\mathcal{D}}
\newcommand{\calE}{\mathcal{E}}
\newcommand{\calF}{\mathcal{F}}

\newcommand{\calP}{\mathcal{P}}

\newcommand{\NN}{\mathbb{N}}

\newcommand{\QQ}{\mathbb{Q}}
\newcommand{\RR}{\mathbb{R}}

\newcommand{\ZZ}{\mathbb{Z}}

\newcommand{\bdy}{\partial} 

\newcommand{\vol}{\operatorname{vol}}

\theoremstyle{plain}
\numberwithin{equation}{section} 
\numberwithin{figure}{section}
\numberwithin{table}{section}
\newtheorem{theorem}{Theorem}[section]
\newtheorem{corollary}[theorem]{Corollary}
\newtheorem{lemma}[theorem]{Lemma}
\newtheorem{conjecture}[theorem]{Conjecture}
\newtheorem*{conjecture*}{Conjecture}
\newtheorem{proposition}[theorem]{Proposition}

\theoremstyle{definition}
\newtheorem{definition}[theorem]{Definition}
\newtheorem{remark}[theorem]{Remark}
\newtheorem*{remark*}{Remark}

\newtheorem*{claim*}{Claim}

\newtheorem*{problem*}{Problem}

\newtheorem{question}[theorem]{Question}
\newtheorem*{question*}{Question}
\newtheorem*{answer*}{Answer}
\newtheorem*{case*}{Case}
\newtheorem*{application*}{Application}

\crefformat{equation}{(#2#1#3)}
\crefmultiformat{equation}{(#2#1#3)}{ and~(#2#1#3)}{, (#2#1#3)}{ and~(#2#1#3)}

\usepackage{standalone}
\newcommand\tLambda{\widetilde{\Lambda}}
\newcommand\tlambda{\tilde{\lambda}}
\newcommand\te{\tilde{e}}
\newcommand\tM{\widetilde{M}}
\newcommand\tS{\tilde{S}}
\newcommand\tE{\tilde{E}}
\newcommand\eps{\varepsilon}
\newcommand\Hom{\operatorname{Hom}}
\newcommand\Frob{\operatorname{Frob}}
\newcommand\inter{\operatorname{int}}
\newcommand\graph{\operatorname{Graph}}
\newcommand\drift{\operatorname{drift}}
\newcommand\stash{\operatorname{Stash}}

\def\co{\colon\thinspace}

\theoremstyle{remark}
\newtheorem{step}{Step}
\newtheorem{proofcase}{Case}

\author{Bal\'azs Strenner}
\date{\today}
\title[Fibrations of 3-manifolds and asymptotic translation length]{Fibrations of 3-manifolds and asymptotic translation length in the arc complex}

\begin{document}

\begin{abstract}
  Given a 3-manifold $M$ fibering over the circle, we investigate how the
  asymptotic translation lengths of pseudo-Anosov monodromies in the arc
  complex vary as we vary the fibration. We formalize this problem by defining
  normalized asymptotic translation length functions $\mu_d$ for every integer
  $d \ge 1$ on the rational points of a fibered face of the unit ball of the
  Thurston norm on $H^1(M;\RR)$. We show that even though the functions $\mu_d$
  themselves are typically nowhere continuous, the sets of accumulation points
  of their graphs on $d$-dimensional slices of the fibered face are rather nice
  and in a way reminiscent of Fried's convex and continuous normalized stretch
  factor function. We also show that these sets of accumulation points depend
  only on the shape of the corresponding slice. We obtain a particularly
  concrete description of these sets when the slice is a simplex. We also
  compute $\mu_1$ at infinitely many points for the mapping torus of the
  simplest hyperbolic braid to show that the values of $\mu_1$ are rather
  arbitrary. This suggests that giving a formula for the functions $\mu_d$
  seems very difficult even in the simplest cases.
\end{abstract}

\maketitle

\section{Introduction}
\label{sec:intro}

To every fibration $M \to S^1$ of a 3-manifold $M$ over the circle, there is an
associated element of $H^1(M;\ZZ)$, the pullback of a generator of
$H^1(S^1;\ZZ)\cong \ZZ$. The integral cohomology classes that correspond to
fibrations of $M$ are organized by the faces of the unit ball of the Thurston
norm $\left\|\cdot\right\|$ on $H^1(M;\RR)$ \cite{Thurston86}: a face $\calF$
can be \emph{fibered}, in which case every integral point in the interior of
the cone $\RR_+\calF$ corresponds to a fibration, or not fibered, in which case
no integral point in $\RR_+\calF$ corresponds to a fibration.

An element $\phi \in H^1(M;\ZZ)$ is \emph{primitive} if it cannot be written in
the form $k \phi'$ for some $\phi' \in H^1(M;\ZZ)$ and integer $k\ge 2$. If an
element $\phi \in H^1(M;\ZZ)$ corresponds to a fibration, then it is primitive
if and only if the fibers are connected. For any $\phi \in H^1(M;\QQ)$, denote
by $\bar{\phi}\in H^1(M;\ZZ)$ the unique primitive integral point on the ray
$\RR_+\phi$.

We can now state a classical result of Fried from 1982, which was the main
motivation for this work. If $M$ admits a complete finite-volume hyperbolic
metric, then the monodromies of the fibrations are pseudo-Anosov mapping
classes by Thurston's Hyperbolization Theorem \cite{Otal96}. For a fibered face
$\calF$ of $M$, define the \emph{normalized entropy function}
\begin{displaymath}
  \xi\co \inter(\calF) \cap H^1(M;\QQ) \to \RR_+
\end{displaymath}
by the formula
\begin{equation}\label{eq:normalized-stretch}
  \xi(\phi) = \left\|\bar{\phi}\right\| \cdot \log{\lambda(\bar{\phi})}
\end{equation}
where $\lambda(\bar{\phi})$ denotes the stretch factor of the pseudo-Anosov
monodromy corresponding to $\bar{\phi}$. In Theorem E of \cite{Fried82}, Fried
proves that the function $\xi$ extends to a convex, continuous function to the
interior of $\calF$ and $\xi(\phi) \to \infty$ as $\phi \to \bdy \calF$. The
goal of this paper is to investigate analogous functions on the rational points
of the fibered faces that are defined not in terms of the stretch factor but
another numerical invariant of pseudo-Anosov maps, the asymptotic translation
length in the arc complex.

The \emph{arc complex} $\calA(S)$ of a connected punctured surface $S$ is a
simplicial complex whose vertices are isotopy classes of properly embedded arcs
in $S$ and whose simplices correspond to collections of disjoint arcs. For two
vertices $\alpha$ and $\beta$ of $\calA(S)$, their distance
$d_{\calA}(\alpha,\beta)$ is defined as the minimal number of edges of a path
in the 1-skeleton of $\calA(S)$ that starts at $\alpha$ and ends at $\beta$.
The \emph{asymptotic translation length} of a mapping class $f$ in the arc
complex is defined as
\begin{displaymath}
  \ell_\calA(f) = \liminf_{n\to \infty} \frac{d_\calA(\alpha,f^n(\alpha))}{n}
\end{displaymath}
where $\alpha$ is any arc.
The number $\ell_\calA(f)$ is a natural invariant encoding geometric
information about the 3-manifold $M$: Futer and Schleimer
\cite{FuterSchleimer14} showed that it is proportional to the height and area
of the boundary of the maximal cusp in $M$.

Based on work of Baik, Shin and Wu \cite{BaikShinWu18}, we define the
\emph{$d$-adic normalized asymptotic translation length function}
\begin{displaymath}
  \mu_d\co \inter(\calF)\cap H^1(M;\QQ) \to \RR_+
\end{displaymath}
by the formula
\begin{equation}\label{eq:normalized-atl}
  \mu_d(\phi) = \left\|\bar{\phi}\right\|^{1+\frac1d} \cdot \ell_\calA(\bar{\phi})
\end{equation}
where $\ell_\calA(\bar{\phi})$ is defined as $\ell_\calA(f)$ where $f$ is the
monodromy of the connected fiber corresponding to $\bar{\phi}$. In order for
$\mu_d$ to be defined, the fibers of $M$ have to be punctured, so $M$ has to be a
cusped 3-manifold. In this paper, we will work under the stronger hypothesis
that the fibered face $\calF$ is \emph{fully-punctured}, meaning that the
singular set of every pseudo-Anosov monodromy in $\RR_+\calF$ is contained in
the set of punctures of the fiber. (If this condition holds for one monodromy
in $\RR_+\calF$, then it holds for all.)

A \emph{$d$-dimensional slice} of a fibered face $\calF$ is an intersection
$\calF \cap \Sigma$ where $\Sigma$ is a $d+1$-dimensional linear subspace of
$H^1(M;\RR)$ intersecting the interior of $\calF$. The slice is \emph{rational}
if $\Sigma \cap H^1(M;\QQ)$ is dense in $\Sigma$.

\begin{theorem}\label{thm:closure}
  Let $M$ be a connected cusped 3-manifold that admits a complete finite-volume
  hyperbolic metric. Let $\calF$ be a fully-punctured fibered face of the unit
  ball of the Thurston norm on $H^1(M;\RR)$. Suppose that
  $1 \le d \le \dim(H^1(M;\RR))-1$ and let $\Omega$ be a rational
  $d$-dimensional slice of $\calF$. Consider
  $\graph(\mu_d|_\Omega) \subset \Omega \times \RR$, the graph of the
  normalized asymptotic translation length function $\mu_d$, restricted to
  $\Omega$.

  There is a continuous function $g\co \inter(\Omega) \to \RR_+$ such that
  $g(\phi) \to \infty$ as $\phi \to \bdy \Omega$ and the set of accumulation points
  of $\graph(\mu_d|_\Omega)$ 
  is
  \begin{displaymath}
    \{(\omega,g(\omega)) : \omega \in \inter(\Omega)\}
  \end{displaymath}
  if $d=1$ and
  \begin{displaymath}
    \{(\omega,r) : \omega \in \inter(\Omega), \,0 \le r \le g(\omega)\} \cup
    (\bdy\Omega \times [0,\infty))
  \end{displaymath}
  if $d\ge 2$.
\end{theorem}

In words, the set of accumulation points is the graph of $g$ if $d=1$ and the
closure of the region under the graph of $g$ if $d \ge 2$.

As an immediate corollary, we have

\begin{corollary}\label{cor:nowhere-cont}
  If $M$, $\calF$, $d\ge 2$ and $\Omega$ are as in \Cref{thm:closure}, then
  $\mu_d|_\Omega$ is a nowhere continuous function.
\end{corollary}

In this sense, the functions $\mu_d$ are therefore very different from Fried's
function $\xi$ which is always continuous. Nevertheless, the properties of
continuity and blowing up at the boundary still make an appearance in
\Cref{thm:closure} for the bounding function $g$.

We derive a formula for $g$ in
\Cref{thm:closure-concrete}. However, it is not clear from this formula whether $g$ is always
convex.

\begin{conjecture}\label{conj:convex}
  Is the function $g$ in \Cref{thm:closure} convex?
\end{conjecture}

When $\Omega$ is a simplex, we are able to describe the function $g$
explicitly. We will show in \Cref{lemma:convex-function} that convexity holds in this case.

\begin{theorem}\label{thm:simplex}
  Let $M$, $\calF$, $d$, $\Omega$ and $g$ be as in \Cref{thm:closure}. Suppose
  $\Omega$ is a simplex with vertices $\omega_1, \ldots, \omega_{d+1}$ and
  define the reparametrization
  \begin{displaymath}
    g^*(\alpha_1, \ldots, \alpha_{d+1}) = g\left(\sum_{i=1}^{d+1}\alpha_i\omega_i \right)
  \end{displaymath}
  of the function $g$ by
  \begin{displaymath}
    \{(\alpha_1,\ldots,\alpha_{d+1}) : \alpha_i>0, \sum_{i=1}^{d+1}\alpha_i =
    1\},
  \end{displaymath}
  the interior of the standard simplex. Let $\Sigma$ be the subspace spanned by
  $\Omega$, let $\Lambda = \Sigma \cap H^1(M;\ZZ)$ be the integral lattice in
  $\Sigma$ and let $\vol_\Lambda$ be the translation-invariant volume form on
  $\Sigma$ with respect to which $\Lambda$ has covolume 1. Then
  \begin{displaymath}
    g^*(\alpha_1, \ldots, \alpha_{d+1}) =\sqrt[d]{\frac{1}{O_d \cdot d!
         \cdot \vol_\Lambda(\Sigma/\langle \omega_1, \ldots,
         \omega_{d+1}\rangle_\ZZ) \cdot \prod_{i=1}^{d+1}\alpha_i}}
  \end{displaymath}
  where $O_d$ is a constant depending only on $d$.

  In the case $d=1$, we have $O_1 = 1$, therefore
  \begin{displaymath}
    g^*(\alpha, 1-\alpha) =\frac{1}{
      \vol_\Lambda(\Sigma/\langle \omega_1,
      \omega_{2}\rangle_\ZZ) \cdot \alpha(1-\alpha)}.
  \end{displaymath}
\end{theorem}

The constant $O_d$ has a concrete interpretation: it is the smallest possible
volume for a $d$-dimensional simplex $\sigma$ in $\RR^d$ with the property that each larger
scaled and translated copy of $\sigma$ ($a\sigma + b$ where $a,b\in \RR$ and $a>1$) contains a point of $\ZZ^d$
in its interior. Although determining the value of $O_d$ for $d\ge 2$ seems to
be an elementary lattice geometry question, we do not even know the value of
$O_2$.

\begin{question}
  What is the value of the constant $O_d$ for $d \ge 2$?
\end{question}

To shed some light on the exact values of the functions $\mu_d$ in addition to
the accumulation points of their graphs, we compute $\mu_1$ at infinitely many
points for the mapping torus of the simplest hyperbolic braid. Both the answer
and the proof are rather ad hoc, suggesting that it is very difficult to
elegantly describe $\mu_d$ even in the simplest cases.

\begin{figure}[ht]
  \centering
  \begin{tikzpicture}[scale=1]
    \def\rad{0.02}
    \draw (0,0) circle (1.6);
    \foreach \i in {1,2,3}
    \filldraw (\i-2,0) circle (\rad);
    \foreach \j in {1,2}{
      \pgfmathsetmacro\i{\j-1.5}
      \draw[->] (\i-0.4,0.2) .. controls (\i-0.2,0.6) and (\i+0.2,0.6) ..
      node[above] {$\sigma_{\j}$}
      (\i+0.4,0.2);
      \draw[->] (\i+0.4,-0.2) .. controls (\i+0.2,-0.6) and (\i-0.2,-0.6) ..
      (\i-0.4,-0.2);
      }
  \end{tikzpicture}
  \caption{The half-twists $\sigma_1$ and $\sigma_2$.}
  \label{fig:sigmas-simple}
\end{figure}
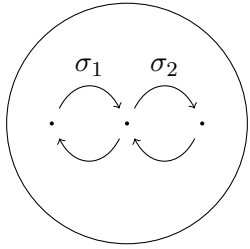

\begin{theorem}\label{thm:example}
  Let $M$ be the mapping torus of the pseudo-Anosov braid
  $f = \sigma_1 \sigma_2^{-1}$ (read in either order) on three strands, see
  \Cref{fig:sigmas-simple}. The fibered face $\calF$ containing $f$ is
  one-dimensional and $f$ corresponds to the midpoint of $\calF$. By choosing a
  linear identification of $\calF$ with $[-1,1]$, we have
  $\mu_1(0) = \frac{8}3$ and
  \begin{displaymath}
    \mu_1(t) =
    \begin{cases}
      \frac{8}3 & \mbox{if } t = \pm \frac12 \\
      4 & \mbox{if } t=\pm \frac13 \\
      \frac{64}{13} & \mbox{if } t=\pm \frac14 \\
      \frac{8}{(1+|t|)^2} & \mbox{if } t=\pm \frac1k \mbox{ when } k\ge 5 \mbox{ is odd} \\
      \frac{8}{1+2|t|-t^2 } & \mbox{if } t=\pm \frac1k \mbox{ when } k\ge 6 \mbox{ is even} \\
    \end{cases}.
  \end{displaymath}
  Moreover,
  \begin{displaymath}
    \lim_{\QQ\ni u\to t} \mu_1(u) = \frac{8}{1-t^2}
  \end{displaymath}
  for all $t \in (-1,1)$. Therefore
  \begin{displaymath}
    \mu_1(t) < \lim_{\QQ\ni u\to t} \mu_1(u)
  \end{displaymath}
  for $t=0$ and all $t = \pm \frac{1}k$ ($k\in \ZZ$, $k\ge 2$) and
  $\mu_1$ is discontinuous at all of these points.
\end{theorem}

In other words, the function $\mu_1$ defined on the 1-dimensional fibered face in
\Cref{thm:example} is discontinuous at every point where we have computed its
value.
We wonder if $\mu_1$ is discontinuous at every rational point of every
1-dimensional slice. More generally:

\begin{question}
  Suppose $M$, $\calF$, $d$, $\Omega$ and $g: \inter(\Omega) \to \RR_+$ are as
  in \Cref{thm:closure}. Does
  \begin{displaymath}
    \mu_d(x) < g(x)
  \end{displaymath}
  hold for every rational point $x$ in the interior of $\Omega$?
\end{question}

It would be interesting to generalize \Cref{thm:closure} in various directions. For example,
one could try to drop the hypothesis that the fibered face $\calF$ is fully-punctured. Instead of the arc complex, one could also consider the curve
complex and define the normalized asymptotic translation length functions
analogously. Our proof has two key ingredients that are specific to the arc
complex:
\begin{itemize}
\item Agol's veering triangulation of 3-manifolds \cite{Agol11} and
\item a theorem of Minsky and Taylor stating that there is a 1-Lipschitz
  retraction from the arc complex $\calA(S)$ to the edges of the veering
  triangulation \cite{MinskyTaylor17}.
\end{itemize}
Generalizing \Cref{thm:closure} to other cases would require replacing these
technical tools with tools suitable in the other cases. Of Sections
\ref{sec:atl-via-graphs}, \ref{sec:cone-estimates}, \ref{sec:frob-lemmas} and
\ref{sec:proof} containing the proof of \Cref{thm:closure}, only Section
\ref{sec:atl-via-graphs} relies crucially on veering triangulations. We use
veering triangulations also for proving \Cref{prop:cone-description}, but as we
remark there, alternative approaches to analogous results already exist. The
remaining parts of \ref{sec:cone-estimates}, \ref{sec:frob-lemmas} and
\ref{sec:proof} should generalize to other cases essentially without
modifications.

\subsection*{Dependence only on shape}
\label{sec:shape}

One interesting property of the functions $\mu_d$ is that up to a constant
factor their bounding function $g$ on any $d$-dimensional slice $\Omega$ only
depends on the shape of $\Omega$. This is in sharp contrast to Fried's
normalized stretch factor function $\xi$, which can take different forms even
on 1-dimensional fibered faces.

\begin{theorem}\label{thm:shape}
  For $i=1,2$, let $M_i$ be 3-manifolds as in \Cref{thm:closure}. Suppose
  $\calF_i \subset H^1(M_i;\RR)$ are fibered faces of $M_i$ and
  $\Omega_i \subset \calF_i$ are $d$-dimensional slices for some integer
  $d \ge 1$. Let $\Sigma_i$ be the span of $\Omega_i$ in $H^1(M_i;\RR)$ and
  consider the lattice $\Lambda_i = \Sigma_i \cap H^1(M_i;\ZZ)$ in $\Sigma_i$.
  Let $g_i: \inter(\Omega_i) \to \RR_+$ be the bounding functions for the
  functions $\mu_d^{\calF_i}|_{\Omega_i}$ as in \Cref{thm:closure}.

  If there is an linear isomorphism $i: \Sigma_1 \to \Sigma_2$ such that
  $i(\Omega_1) = \Omega_2$, then
  \begin{displaymath}
    g_2(i(\phi_1)) = \theta^{\frac1d} g_1(\phi_1)
  \end{displaymath}
  holds for all $\phi_1 \in \inter(\Omega_1)$ for
  \begin{displaymath}
    \theta = \frac{\vol(\Sigma_2/\Lambda_2)}{\vol(\Sigma_2/i(\Lambda_1))},
  \end{displaymath}
  where $\vol$ is any translation-invariant volume form on $\Sigma_2$.
\end{theorem}

\subsection*{Related results}

In \cite{KinShin17}, Kin and Shin have shown that the function $\mu_1^\calC$,
defined analogously to $\mu_1$ using the curve complex instead of the arc
complex, is bounded from above on infinite subsets of slices arising from
projecting an arithmetic progression in $H^1(M;\ZZ)$ onto $\calF$. Using this,
they improved the upper bound of a result of Gadre and Tsai \cite{GadreTsai11}
stating that the minimal asymptotic translation length in the curve complex for
pseudo-Anosov maps on the closed surface $S_g$ of genus $g$ is between
$\frac{C_1}{g^2}$ and $\frac{C_2}{g^2}$ for some $C_1$ and $C_2$. Using the
bounds on $\mu_1^\calC$, the paper \cite{KinShin17} also provides upper bounds
for the minimal asymptotic translation length in the curve complex for certain
sequences of punctured surfaces (for this, see also \cite{Valdivia14}),
handlebody groups and hyperelliptic handlebody groups. It would be interesting
to investigate the implications of our more explicit description of the
function $\mu_1$ on similar questions.

In \cite{BaikShinWu18}, Baik, Shin and Wu have studied the function
$\mu_d^\calC$, defined analogously to $\mu_d$ using the curve complex instead
of the arc complex. They proved that the function $\mu^\calC_d$ is bounded from
above on compact $d$-dimensional polytopes contained in the interior of
$\calF$. In Conjecture 1 of their paper, they conjecture that their bound is
sharp in the sense that for each $d \ge 2$ there exist $M$, $\calF$, $\Omega$
such that the function $\mu^\calC_d$ is bounded away from 0 on an infinite
subset of $\Omega$. Although in the arc complex instead of the curve complex,
our \Cref{thm:closure} verifies the stronger statement that $\mu_d$ is bounded
away from zero on an infinite subset of $\Omega$ \emph{for all choices} of $M$,
$\calF$ and $\Omega$. Moreover, in addition to showing that the values are
bounded away from zero, \Cref{thm:closure} precisely specifies the values
$\mu_d$ can approach along accumulating sequences in $\Omega$.

\subsection*{Acknowledgements}

We thank the anonymous referees for their thorough review and for pointing out an error in an earlier version of the paper.

\section{Background}
\label{sec:background}

\subsection{Fibrations over the circle}
\label{sec:fibrations-over-circle}

Let $\pi\co M \to S^1$ be a fibering of $M$ over the circle with fiber
$S = \pi^{-1}(0)$. Let $\phi\in H^1(M)$ be the pullback of one of the two
generators of $H^1(S^1) \cong \ZZ$. There is an infinite cyclic cover
$S \times \RR \to M$ corresponding to the homomorphism $\phi:\pi_1(M) \to \ZZ$.
Let $h\co S\times \RR \to S\times \RR$ be the element of the deck group that
maps $S \times \{1\}$ to $S\times \{0\}$. The composition
\begin{displaymath}
  S \times \{0\} \to S\times \{1\} \stackrel{h}{\to} S\times \{0\},
\end{displaymath}
where the first map is the isotopy map $(x,0) \mapsto (x,1)$ in the product
$S \times \RR$, yields a homeomorphism of $S$. The mapping class of this
homeomorphism is the \emph{monodromy} $f$ of the fibration. The monodromy
depends on which of the two generators we pick for $S^1$.

In other words, we can present the 3-manifold $M$ as the quotient
\begin{displaymath}
  M = (S \times \RR)/\langle (x,t) \sim (\psi(x), t-1) \rangle
\end{displaymath}
for any homeomorphism $\psi:S \to S$ representing the mapping class $f$.

The map $((x,t),s) \mapsto (x,t+s)$ defines a flow $(S\times \RR) \times \RR
\to S\times \RR$ which is $h$-equivariant, therefore it descends to a map
$M \times \RR \to M$, defining the \emph{suspension flow} on $M$.

Finally, we will use the following conventions to make the discussions
more intuitive. We image the product $S \times \RR$ such that the
$\RR$-coordinate axis is vertical and where $\infty$ is up and $-\infty$ is
down. So if $t_1 < t_2$, then we will say that the slice $S \times \{t_2\}$ is
\emph{above} $S\times \{t_1\}$ and $S \times \{t_1\}$ is \emph{below}
$S\times \{t_2\}$.

\subsection{Pseudo-Anosov monodromies}
When $\pi\co M \to S^1$ is a fibration and $M$ is hyperbolic, the monodromy $f$
is pseudo-Anosov by Thurston's Hyperbolization Theorem. Let
$\lambda^{\pm} \subset S$ be the invariant singular (unmeasured) foliations of
$f$. We will refer to the foliation whose leaves are expanded by $f$ as the
\emph{horizontal foliation} and the foliation whose leaves are contracted as
the \emph{vertical foliation}.

Since $\lambda^{\pm}$ are invariant under the monodromy, their orbits under the
suspension flow are singular 2-dimensional foliations $\Lambda^{\pm}$ in $M$,
transverse to the fibers, whose singular set is the suspension of the singular
points of $\lambda^{\pm}$. Conversely, the foliations $\lambda^{\pm}$ can be
obtained from $\Lambda^{\pm}$ by taking the intersection of $\Lambda^{\pm}$
with the fibers.

\subsection{Fully-punctured fibered faces}

Let $M$ be a hyperbolic 3-manifold and let $\calF$ be a fibered face of the
Thurston norm ball of $H^1(M;\RR)$. Every integral cohomology class in the
interior of the cone $\RR_+ \calF$ corresponds to a fibration of $M$ over the
circle.

Fried \cite{Fried82} (see also McMullen \cite[Corollary 3.2]{McMullen00})
showed that the suspension foliation $\Lambda^{\pm}$ constructed from any two
fibrations in this fibered cone are the same (up to isotopy) when the singular
points of $\lambda^{\pm}$ are all at punctures of the fiber for \emph{some}
fibration. It follows that in this case $\Lambda^{\pm}$ does not have any
singular points and therefore the singular points of $\lambda^{\pm}$ are all at
punctures of the fiber for \emph{all} fibrations in this fibered cone. Such a
fibered face $\calF$ is called \emph{fully-punctured}.

\subsection{Relating different fibrations}

From now on, suppose that $\calF$ is a fully-punctured fibered face. The
\emph{maximal abelian cover} $\tM$ of $M$ is the cover corresponding to the
natural homomorphism 
\begin{displaymath}
  \pi_1(M) \to G = H_1(M;\ZZ)/\mathrm{torsion}.
\end{displaymath}
The foliations $\Lambda^{\pm}$ in $M$ lift to foliations $\tLambda^{\pm}$ in $\tM$.
The suspension flow on $M$ also lifts to a flow on $\tM$, leaving invariant the
foliations $\tLambda^{\pm}$. The leaf space of this flow is homeomorphic to a
surface $\tS$ and the foliations $\tLambda^{\pm}$ define foliations
$\tlambda^{\pm}$ on $\tS$.

Every fiber of every fibration in the cone $\RR_+ \calF$ is a quotient of the
foliated surface $(\tS,\tlambda^{\pm})$ by a covering map. One can see this as
follows. Let $\phi \in \RR_+ \calF$ be a primitive integral point with
monodromy $f$, fiber $S$ with stable and unstable foliations $\lambda^{\pm}$.
The covering $\tM \to M$ factors through the infinite cyclic covering
$S \times \RR \to M$ induced by $\phi:\pi_1(M) \to \ZZ$. The lift of $S$ to
$S\times \RR$ is an infinite collection of parallel copies of $S$. Under the
covering $\tM \to S\times \RR$, each copy lifts to a surface intersecting every
flowline of $\tM$ exactly once. This gives rise to a foliation-preserving
covering $(\tS,\tlambda^\pm) \to (S,\lambda^\pm)$ whose deck group is the
kernel of the homomorphism $G \to \ZZ$ induced by $\phi$.

\subsection{Veering triangulations}

This section recalls some facts about veering triangulations of hyperbolic
3-manifolds defined by Agol \cite{Agol11}, refined by Gu\'eritaud
\cite{Gueritaud16} and further studied by Minsky and Taylor
\cite{MinskyTaylor17}.

When the foliations $\lambda^\pm$ for some fibration are endowed with the
measure invariant under the pseudo-Anosov monodromy and this measure is lifted
to $\tlambda^\pm$, the measured foliations $\tlambda^\pm$ endow the surface
$\tS$ with a singular Euclidean metric. Denote the metric completion of $\tS$
by $\hat{S}$. Each completion point in $\hat{S}-\tS$ is an isolated
point whose small neighborhood minus the completion point cover the neighborhood of
a puncture in any fiber $S$. The metric on $\tS$ depends on the fibration
chosen in the construction, but the topology of $\hat{S}$ does not. In the
future, we will ignore the metric and consider $\hat{S}$ together with the
\emph{unmeasured} foliations $\hat\lambda^{\pm}$ obtained from $\tlambda^\pm$
by extending to the completion points.

A \emph{singularity-free rectangle} in $(\hat{S}, \hat\lambda^\pm)$ is an
immersion $[0,1]^2 \to \hat{S}$ such that the vertical and horizontal
foliations of $[0,1]^2$ map to $\hat\lambda^\pm$ and the interior of
the rectangle does not contain any completion point of $\hat{S}-\tS$. By the
\emph{interior} of the rectangle, we mean the
image of $(0,1)^2$ under the immersion. Similarly, by the
\emph{boundary} of the rectangle, we mean the image of the boundary of
$[0,1]^2$ under the immersion.

A singularity-free rectangle is \emph{maximal} if all four sides of $[0,1]^2$
contain the preimage of a completion point in their interior under the
immersion map. (Each side may contain only one completion point, since the horizontal and vertical foliations, being invariant foliations of a pseudo-Anosov map, cannot have saddle connections.) By connecting each pair of the four points with an arc
inside $[0,1]^2$ and considering the image under the immersion map, we obtain
six arcs in $\hat{S}$, forming a flattened tetrahedron in $\hat{S}$.

These arcs are defined only up to isotopy. To make the choice of the arcs
canonical, we choose a fibration in our fibered cone and---as we have seen
above---this choice endows $\hat{S}$ with a singular Euclidean metric. We choose
the arcs to be the unique geodesics in their isotopy class in this metric.

We think of the arc connecting the horizontal sides to be \emph{above} the arc
connecting the vertical sides. So the two triangles containing the arc
connecting the vertical sides are the two bottom triangles and the remaining
two triangles are the two top triangles of the tetrahedron.

Consider \emph{all} maximal singularity-free rectangles in $\hat{S}$ and all
the arcs, triangles and tetrahedra they define through this process. For each
triangle the smallest singularity-free rectangle containing it can be enlarged
in two ways to a maximal singularity-free rectangle: we can enlarge the
rectangle horizontally or vertically. In the former case, we obtain a
tetrahedron that contains our triangle as one of the two top triangles. In the
latter case, we obtain a tetrahedron that contains our triangle as one of the
two bottom triangles. Therefore the tetrahedra glue together in a layered
fashion.

The links of the triangulation around the vertices are not spheres. Instead we
glue up the \emph{ideal} tetrahedra that do not include the vertices. With more
work (see \cite{Gueritaud16}), one can check that the links of the edges are
circles, so the ideal tetrahedra glue up to a 3-manifold. Moreover, this
3-manifold is homeomorphic to $\tM \cong \tS \times \RR$ and the ideal
triangulation is called the \emph{veering triangulation} of $\tM$. The veering
triangulation is invariant under the $G$-action, and the quotient is the
veering triangulation of $M$.

We conclude by comparing the conventions regarding \emph{above} and
\emph{below} introduced in this section versus the conventions introduced
earlier. Recall from \Cref{sec:fibrations-over-circle} that for any fibration
in our fibered cone, the deck transformation $h\co S\times \RR \to S\times \RR$
satisfies $h(x,t) = (\psi(x),t-1)$ where $\psi$ is a pseudo-Anosov
homeomorphism representing the monodromy mapping class. Recall also that our
convention is that $\psi$ expands horizontally and contracts vertically.
Therefore the tetrahedra and the corresponding maximal singularity-free
rectangles become wider and shorter as we go down in the product
$\tS \times \RR$. This is consistent with the convention that the top edge of
each tetrahedron, connecting the horizontal sides of the corresponding
rectangle, has larger slope than the bottom edge, connecting the vertical
sides.

\section{Asymptotic translation length via cycles in graphs}
\label{sec:atl-via-graphs}

\subsection{Intersecting edges of the veering triangulation}
\label{sec:intersecting-edges}

Given an edge of the veering triangulation of $\tM \cong \tS \times \RR$, its
projection onto $\tS$ is an arc in $\tS$. We say that two edges
\emph{intersect} if their projections intersect in $\tS$. Otherwise we say that
the two edges are \emph{disjoint}. Recall that we have chosen these arcs to be
geodesics in a singular Euclidean metric, so the arcs are automatically in
minimal position and we do not need to be concerned about isotopies. Recall
also that the edges do not have endpoints, so if they intersect, they have to
intersect in their interiors.

For our applications, it will be important to keep track of which pairs of
edges of the veering triangulation of $\tM$ intersect and which two are
disjoint. We can organize this information as follows.

Let $E$ be the set of edges of the veering triangulation of $M$. The set $E$ is
finite, which follows from Agol's construction of the veering triangulations by
periodic train track sequences \cite{Agol11}. For each edge $e \in E$,
\emph{choose} a lift $\te$ in the veering triangulation of $\tM$. Denote the
set of these lifts by $\tE$. Each edge of the veering triangulation of $\tM$
can be uniquely written in the form $g\te$ for some $g \in G$ and
$\te \in \tE$.

When two edges $g\te$ and $g'\te'$ intersect, one of the edges is \emph{above}
the other with respect to the pseudo-Anosov flow. If $g\te$ is above $g'\te'$,
we write $g\te > g'\te'$. By our conventions, $g\te$ is above $g'\te'$ if they
intersect and the smallest singularity-free rectangle containing $g'\te'$ is
wider and shorter than the smallest rectangle containing $g\te$.

\begin{definition}[Stashing set]\label{def:stashing-set}
  For any $e, e' \in E$, introduce the notation
\begin{displaymath}
  \stash(e,e') = \{g \in G: \te > g\te' \}.
\end{displaymath}
In words, $\stash(e,e')$ is the set of deck transformations in $G$ that send (or \emph{stash}) $\te'$ below $\te$, hence we call $\stash(e,e')$ the \emph{stashing set} of $\te'$ with respect to $\te$.
\end{definition}

The knowledge of the sets $\stash(e,e')$ for all pairs $e,e'\in E$
contains all disjointness information, since $g\te > g'\te'$ if and
only if $g^{-1}g' \in \stash(e,e')$.

\subsection{Frobenius numbers}
\label{sec:frob-numbers}

We define the \emph{Frobenius number} of a function $\beta\co A \to \ZZ$ as
\begin{equation}\label{eq:frob-def}
  \Frob(\beta) = \max \big(\ZZ-\beta(A)\big)
\end{equation}
if the maximum exists.

We remark that this notion is closely related to the Frobenius coin problem
\cite{Alfonsin05} that, given relatively prime positive integers
$a_1, \ldots, a_n$, asks for the largest integer that cannot be written as a
linear combination of $a_1, \ldots, a_n$ with nonnegative integer coefficients.
Indeed, let $H$ be the free commutative monoid generated by $x_1,\ldots,x_n$
and let $\beta\co H \to \ZZ$ be a homomorphism such that
$\beta(x_1),\ldots, \beta(x_n)$ are positive. Then the Frobenius number of
$\beta$, as defined in \Cref{eq:frob-def}, is the largest integer that cannot
be written as a nonnegative integral linear combination of
$\beta(x_1), \ldots, \beta(x_n)$.

\subsection{Translation length in the arc complex via graphs}

To every primitive integral class $\phi$ in the interior of $\RR_+\calF$, we
associate a weighted directed graph $W(\phi)$ on the vertex set $E$. There is
an edge from $e$ to $e'$ if and only if there is at least one integer that is
\emph{not} contained in the subset
\begin{displaymath}
  -\phi(\stash(e,e')) \cup \phi(\stash(e',e))
\end{displaymath}
of $\ZZ$. Here $\phi$ stands for the surjective linear functional $G \to \ZZ$
associated to $\phi$. If there is an edge from $e$ to $e'$, then its weight
$w(ee')$ is defined as the largest integer not contained in the subset
$-\phi(\stash(e,e'))$ of $\ZZ$. Alternatively,
\begin{equation}
  \label{eq:weight}
  w(ee') = \Frob(\phi|_{-\stash(e,e')}).
\end{equation}
In \Cref{cor:frob-exists-for-intersection-data}, we will see that this largest
integer always exists and therefore $w(ee')$ is always well-defined.

\Cref{lemma:weight-and-disjointness} below will explain the information
contained by the weighted graph $W(\phi)$. First we need the following lemma.

\begin{lemma}\label{lemma:gap}
  Any element of $-\phi(\stash(e,e'))$ is larger than any element of
  $\phi(\stash(e',e))$. In addition, if
  $-\phi(\stash(e,e')) \cup \phi(\stash(e',e))$ is not all of $\ZZ$, then the
  difference between the smallest element of $-\phi(\stash(e,e'))$ and the
  largest element of $\phi(\stash(e',e))$ is at least 2.
\end{lemma}
\begin{proof}
  To prove the first statement, let $g_1, g_2\in G$ such that $\te > g_1\te'$ and
  $\te' > g_2\te$. Since the relation $>$ is transitive, we have $\te' > g_1g_2\te'$.
  One should think of $\phi$ as a height function: since $g_1g_2\te'$ is below
  $\te'$, we have $\phi(g_1g_2) = \phi(g_1) + \phi(g_2) < 0$. So $-\phi(g_1)$ is
  indeed larger than $\phi(g_2)$.

  Assume that the difference is between the smallest element of
  $-\phi(\stash(e,e'))$ and the largest element of $\phi(\stash(e',e))$ is 1.
  Then there are $g_1, g_2$ such that $\te > g_1\te'$ and $\te' > g_2\te$ and
  $\phi(g_1g_2) = -1$. As before, we have $\te' > g_1g_2\te'$. So
  \begin{displaymath}
    \te > g_1\te' > g_1(g_1g_2)\te' > g_1(g_1g_2)^2\te' > \cdots
  \end{displaymath}
  which means that every integer at least $-\phi(g_1)$ is contained in
  $-\phi(\stash(e,e'))$. Similarly, we obtain that every integer at most
  $\phi(g_2)$ is contained in $\phi(\stash(e',e))$. Since the gap between
  $-\phi(g_1)$ and $\phi(g_2)$ is 1, we have $-\phi(\stash(e,e')) \cup
  \phi(\stash(e',e)) = \ZZ$. This proves the second statement.
\end{proof}

In the following lemma, $S$ is the fiber of the fibration corresponding to
$\phi$, $f$ is the monodromy and $p_\phi$ is the composition
$\tM \cong \tS\times \RR \to \tS \to S$.

\begin{lemma}\label{lemma:weight-and-disjointness}
  There is an edge from $e$ to $e'$ in $W(\phi)$ if and only if there exists an
  integer $k$ such that $p_\phi(\te)$ and $f^{k}(p_\phi(\te'))$ are disjoint
  in $S$. Moreover, if there is an edge from $e$ to $e'$, then its weight
  $w(ee')$ is the largest integer $k$ such that $p_\phi(\te)$ and
  $f^{k}(p_\phi(\te'))$ are disjoint in $S$.
\end{lemma}
\begin{proof}

  One can check step by step that following are equivalent for any integer $k$:
  \begin{enumerate}
  \item The arcs $p_\phi(\te)$ and $f^{k}(p_\phi(\te'))$ are disjoint in $S$.
  \item The edge $\te$ is disjoint from \emph{all} lifts of $f^{k}(p_\phi(\te'))$ to
    $\tS$.
  \item The edge $\te$ is disjoint from $g \te'$ for all $g \in G$ with
    $\phi(g) = -k$.
  \item $\te \not> g\te'$ and
      $\te \not< g\te'$ for all $g \in G$ with $\phi(g) = -k$.
  \item $g \notin \stash(e,e')$ and $g^{-1} \notin \stash(e',e)$ for all
      $g \in G$ with
      $\phi(g) = -k$.
  \item $-k \notin \phi(\stash(e,e'))$ and $k \notin \phi(\stash(e',e))$.
  \item $k \notin -\phi(\stash(e,e')) \cup \phi(\stash(e',e))$.
  \end{enumerate}
  By definition, there is an edge from $e$ to $e'$ in $W(\phi)$ if an integer
  $k$ satisfying the last statement exists. The first statement of the lemma follows.

  Using the equivalences again, the largest integer $k$ such that $p_\phi(\te)$
  and $f^{k}(p_\phi(\te'))$ are disjoint in $S$ is the largest $k$ that is not
  contained in either $-\phi(\stash(e,e'))$ or $\phi(\stash(e',e))$. But if there
  exists such a $k$, then by \Cref{lemma:gap} it is the largest $k$ that is not
  contained in $-\phi(\stash(e,e'))$. The second statement now follows from the
  definition of $w(ee')$.
\end{proof}

\begin{corollary}\label{cor:frob-exists-for-intersection-data}
  For any pair $e,e' \in E$ and primitive integral class $\phi$ in the interior
  of $\RR_+\calF$, there exists some integer $N$ such that $n \in -\phi(\stash(e,e'))$
  for all $n>N$.
\end{corollary}
\begin{proof}
  We will show that there exists some $N$ such
  that for all $n>N$ the arcs $p_\phi(\te)$ and $f^n(p_\phi(\te'))$ are not
  disjoint. Using the equivalences in the proof of
  \Cref{lemma:weight-and-disjointness}, the statement will follow.
  
  Let $R$ and $R'$ be the rectangles with horizontal and vertical sides whose diagonals are
  $p_\phi(\te)$ and $p_\phi(\te')$, respectively.
  Let $s'$ be a horizontal side of $R'$. The side $s'$ is a starting segment of a horizontal
  separatrix emanating from a singularity. For any $n$, the segment $f^n(s')$ is the starting segment of a horizontal separatrix. There are finitely many horizontal separatrices and each such separatrix is dense in the surface. The map $f$
  stretches the surface horizontally by the stretch factor, so if $n$ is large enough, then $f^n(s')$ intersects
  the interior of $R$ and consequently both horizontal sides of $f^n(R')$ intersect the interior of $R$
  and therefore $f^n(p_\phi(\te'))$, the diagonal of $f^n(R')$, intersects $p_\phi(\te)$, the diagonal of
  $R$.
\end{proof}

We define the \emph{average weight of a cycle} $\gamma = e_1\ldots e_n e_1$ in
$W(\phi)$ as
\begin{displaymath}
  \bar{w}(\gamma) = \frac{w(e_1e_2) + \cdots + w(e_{n-1}e_n) + w(e_ne_1)}{n}.
\end{displaymath}

\begin{proposition}[Asymptotic translation length via weighted graphs]\label{prop:atl-via-weighted-graphs}
  For any primitive integral class $\phi$ in the interior of the cone $\RR_+F$,
  the asymptotic translation length in the arc complex of the pseudo-Anosov
  monodromy corresponding to $\phi$ is
  \begin{displaymath}
    \ell_\calA(\phi) = \frac{1}{\max\{\bar{w}(\gamma): \gamma \mbox{ is a cycle in } W(\phi)\}}.
  \end{displaymath}
\end{proposition}
\begin{proof}
  First note that the maximum is indeed realized, since any cycle decomposes
  to minimal cycles and the average weight of the cycle is at most the average
  weight of the minimal cycle with the largest weight. The fact that there is
  at least one cycle in $W(\phi)$ will follows from the rest of the proof.

  For any cycle $\gamma = e_1\ldots e_d e_{1}$, extend the sequence $e_i$ for
  $i\ge d+1$ such that $e_{i+d} = e_i$ for all integer $i\ge 1$. Consider the
  sequence
  \begin{displaymath}
    p_\phi(\te_1), \; f^{w(e_1e_2)}(p_\phi(\te_2)),\;
    f^{w(e_1e_2)+w(e_2e_3)}(p_\phi(\te_3)),\;  \ldots,
  \end{displaymath}
  of arcs in $S$. By
  \Cref{lemma:weight-and-disjointness}, consecutive arcs are disjoint.
  Therefore we have
  \begin{displaymath}
    d_\calA\left(p_\phi(\te_1), f^{\sum_{i=1}^{nd} w(e_{i}e_{i+1})}\left(p_\phi(\te_1)\right)\right) \le nd
  \end{displaymath}
  for any integer $n\ge 1$. This demonstrates that
  \begin{multline*}
    \ell_\calA(\phi) = \lim_{n\to \infty} \frac{d_\calA\left(p_\phi(\te_1),
        f^{\sum_{i=1}^{nd}
          w(e_{i}e_{i+1})}\left(p_\phi(\te_1)\right)\right)}{\sum_{i=1}^{nd}
      w(e_{i}e_{i+1})} \le \\ \le \lim_{n\to \infty} \frac{nd}{\sum_{i=1}^{nd} w(e_{i}e_{i+1})} = \frac{1}{\bar{w}(\gamma)}.
  \end{multline*}
  Since this inequality holds for any cycle $\gamma$, it follows that the left
  hand side in the proposition is bounded from above by the right hand side.

  A key ingredient for the inequality in the reverse direction is a result of
  Minsky and Taylor \cite[Theorem 1.4]{MinskyTaylor17} that states that there
  is a 1-Lipschitz retraction from the arc complex $\calA(S)$ to the set of
  arcs that are projections of the edges of the veering tringulation of $\tM$
  under $p_\phi$. In particular, any two arcs $f^{k}(p_\phi(\te))$ and
  $f^{k'}(p_\phi(\te'))$ in $S$ are joined by a geodesic in the arc complex
  $\calA(S)$ whose vertices are all of the form $f^{k''}(p_\phi(\te''))$.

  Fix $\te_1\in \tE$ and let $n$ be a positive integer. Denoting the distance
  $d_{\calA}(p_\phi(\te_1), f^{n}(p_\phi(\te_1)))$ by $d_n$, there is a sequence
  of arcs
  \begin{displaymath}
    f^{k_1}(p_\phi(\te_1)),\; f^{k_2}(p_\phi(\te_2)),\; \ldots,\;
    f^{k_{d_n+1}}(p_\phi(\te_{d_n+1}))
  \end{displaymath}
  such that consecutive arcs are disjoint in $S$, $\te_{d_n+1} = \te_1$, the
  $k_i$ are integers and $k_{d_n+1}-k_1 = n$. By
  \Cref{lemma:weight-and-disjointness}, there is an edge from $e_i$ to
  $e_{i+1}$ in $W(\phi)$ and we have $k_{i+1}-k_i \le w(e_{i}e_{i+1})$ for all
  $i = 1,\ldots,d_n$. Summing these inequalities, we obtain
  $n \le \sum_{i=1}^{d_n} w(e_{i}e_{i+1})$. After dividing both sides by $d_n$ and
  taking reciprocals, we have
  \begin{displaymath}
    \frac{d_n}{n} \ge \frac1{\bar{w}(\gamma_n)}
  \end{displaymath}
  where $\gamma_n$ denotes the cycle $e_1\ldots e_{d_n}e_1$. In particular,
  this shows that there is at least one cycle in $W(\phi)$.

  Note that
  \begin{displaymath}
    \ell_\calA(\phi) = \lim_{n \to \infty}
    \frac{d_{\calA}(p_\phi(\te_1),f^{n}(p_\phi(\te_1)))}{n} = \lim_{n\to \infty}\frac{d_n}{n}
    \ge \liminf_{n\to\infty}\frac{1}{\bar{w}(\gamma_n)}.
  \end{displaymath}
  The right hand side in the proposition is a lower bound for
  $\frac{1}{\bar{w}(\gamma)}$ for any cycle $\gamma$ in $W(\phi)$, therefore it is
  also a lower bound for $\ell_\calA(\phi)$. This completes the proof of the
  reverse inequality.
\end{proof}

\subsection{The graph $\Delta$}
\label{sec:delta}

In this section, we introduce a digraph $\Delta$ that serves as a model for the
veering triangulation. We will use this graph to compute the weighted graphs
$W(\phi)$ discussed in the previous section.

The vertices and edges of $\Delta$ correspond to the tetrahedra and the
triangles, respectively, of the veering triangulation of $M$. The edge
corresponding to a triangle $t$ starts at the tetrahedron that has $t$ as one
of its two bottom triangles and ends at the tetrahedron that has $t$ as one of
its two top triangles. Note that every vertex has exactly two outgoing and two
incoming edges.

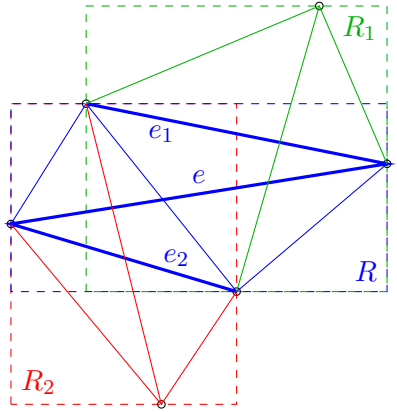
\begin{figure}
  \centering
  \begin{tikzpicture}
    \def\rad{0.05}
    \draw[red,dashed] (0,1) rectangle (3,5);
    \draw[blue,dashed] (0,2.5) rectangle (5,5);
    \draw[green!70!black,dashed] (1,2.5) rectangle (5,6.3);
    \foreach \x/\y in {1/5, 0/3.4, 2/1, 3/2.5, 5/4.2, 4.1/6.3}
    \draw (\x,\y) circle (\rad);
    \draw[blue] (5,4.2) -- (3,2.5) -- (1,5) -- (0,3.4);
    \draw[very thick,blue] (1,5) -- (5,4.2) -- node[above]
    {$e$} (0,3.4) --  (3,2.5);
    \node[blue,below] at (2,4.85) {$e_1$};
    \node[blue,above] at (2.2,2.7) {$e_2$};
    \draw[red] (3,2.5) -- (2,1) -- (1,5); \draw[red] (2,1) -- (0,3.4);
    \draw[green!70!black] (1,5) -- (4.1,6.3) -- (3,2.5); \draw[green!70!black] (4.1,6.3) --
    (5,4.2);
    \node[red,above right] at (0,1) {$R_2$};
    \node[green!70!black,below left] at (5,6.3) {$R_1$};
    \node[blue,above left] at (5,2.5) {$R$};
  \end{tikzpicture}
  \caption{The edges $e_1$ and $e_2$ representing the vertices of $\Delta$ such
    that there is an edge from those vertices to the vertex represented by $e$.}
  \label{fig:delta-edges}
\end{figure}

There is a one-to-one correspondence from the tetrahedra to the edges of the
veering triangulation that assigns to each tetrahedron its bottom edge. Using
this correspondence, we can alternatively think about the vertices of $\Delta$
as edges of the veering triangulation. For an edge $e\in E$,
\Cref{fig:delta-edges} illustrates the two other edges $e_1,e_2\in E$ such
there is an edge of $\Delta$ from $e_i$ to $e$. We can describe $e_1$ and $e_2$
as follows. Expand the smallest singularity-free rectangle containing $e$
vertically as far as possible---the four singularities on the boundary of
resulting rectangle $R$ define the tetrahedron $T$ whose bottom edge is $e$.
The edges $e_1$ and $e_2$ are the two edges of this tetrahedron that are
neither the top nor the bottom edges such that the interiors of the rectangles
$R_1$ and $R_2$ obtained by expanding the smallest singularity-free rectangle
containing $e_1$ and $e_2$ vertically cover the interior of $R$. This is
because $e_1$ and $e_2$ are the bottom edges of tetrahedra $T_1$ and $T_2$
determined by $R_1$ and $R_2$, respectively, and both $T_1$ and $T_2$ have a
bottom triangle that is top triangle of $T$.

We also label each edge of $\Delta$ by an element of $G$, called the
\emph{drift} of the edge. To do this, choose a lift $\tilde{t}$ of the triangle
$t$ corresponding to an edge of $\Delta$ in the veering triangulation of $\tM$.
If the bottom edges of the tetrahedra right above and right below $\tilde{t}$
are $g\te$ and $g'\te'$, respectively, then the drift of the edge of $\Delta$
corresponding to $t$ is $g^{-1}g'$. The drift measures how much the coefficient
of the edge $g\te$ changes as we proceed to the other edge $g'\te'$. Note
that this definition is independent of the choice of the lift $\tilde{t}$.

Recall that a digraph is strongly connected if there is path from any
vertex to any other vertex. We will need the following lemma later.
\begin{lemma}\label{lemma:Delta-strongly-connected}
  The graph $\Delta$ is strongly connected.
\end{lemma}
\begin{proof}
  Since pseudo-Anosov homeomorphisms of surfaces have dense orbits
  (\cite[Expos\'e 9]{FLP}), there is a dense flow line in $M$. A dense flow
  line visits every tetrahedron of the veering triangulation infinitely often.
  Associated to this flow line is a bi-infinite path in $\Delta$ visiting every
  vertex infinitely often. This shows that $\Delta$ is strongly connected.
\end{proof}

\subsection{The extended graph $\Delta^*$}
\label{sec:delta-star}

We also define a graph $\Delta^*$, obtained by adding some more labeled edges
to $\Delta$, one for each tetrahedron of the veering triangulation of $M$.
For each tetrahedron $T$, we create an edge from the top edge to the bottom
edge of $T$. To define the label of this edge, choose a lift of $T$
in the veering triangulation of $\tM$. If the bottom and top edges of the lift are
$g_1\te_1$ and $g_2\te_2$, then our edge points from $e_2$ to $e_1$
and has drift $g_2^{-1}g_1$. This definition is also independent of the choice
of the lift.

To distinguish between the edges of $\Delta$ and the edges of $\Delta^*$ that
are not in $\Delta$, we will call the two types of edges
\emph{triangle-edges} and \emph{tetrahedron-edges}, respectively, as a reminder
that they correspond to triangles and tetrahedra of the veering triangulation.

\subsection{Computing the stashing sets}
\label{sec:computing-intersection}

In this section we explain how the stashing sets $\stash(e,e')$ can be
computed from the digraph $\Delta^*$. We begin with a few definitions.

A \emph{path} in the digraph $\Delta^*$ is a sequence of edges
$\eps_1, \eps_2, \ldots, \eps_n$ of $\Delta^*$ ($n\ge 1$) such that the endpoint
of $\eps_i$ is the same as the starting point of $\eps_{i+1}$ for all
$i = 1, \ldots, d-1$. We caution the reader that we cannot refer to edges of
$\Delta^*$ simply by their endpoints, since there might be multiple edges between
vertices, see \Cref{fig:delta}, for example.

The \emph{drift} of a path is the product of the drifts of the edges of the
path. Formally, the drift of the path $\eps_1\eps_2\ldots\eps_n$ is
\begin{displaymath}
  \prod_{i=1}^n \drift(\eps_i) \in G
\end{displaymath}
where $\drift(\eps_i) \in G$ denotes the drift of the edge $\eps_i$.

A \emph{good path} is a path whose first edge is a tetrahedron-edge and whose
remaining edges are triangle-edges.

\begin{proposition}[Stashing sets
  via good paths]\label{prop:intersection-and-good-paths}
  We have
  \begin{displaymath}
    \stash(e,e') = \{\drift(\gamma): \mbox{$\gamma$ is a good path from
      $e$ to $e'$ in $\Delta^*$}\}.
  \end{displaymath}
\end{proposition}
\begin{proof}
  First we show that the right hand side contains the left hand side. Suppose
  $g' \in \stash(e,e')$, which means that $\te > g'\te'$. Let $p$ be the
  intersection of the images of $\te$ and $g'\te'$ in $\tS$ under the
  projection $\tM \to \tS$ by collapsing the flow lines. The preimage of $p$ is
  a flow line that intersects $\te$ and $g'\te'$. Consider the subinterval $I$ of this flow line
  between $\te$ and $g'\te'$. 
  
  If $I$ does not intersect an edge of the veering triangulation aside from its endpoints, then it passes through a sequence of tetrahedra in a way that it enters each subsequent tetrahedron through the interior of a face. Denote the bottom
  edges of these tetrahedra by $g_1\te_1, \ldots, g_k\te_k$, where $g_1\te_1$ is the bottom edge of the
  tetrahedron whose top edge is $\te$ and $g_k\te_k = g'\te'$. Then there is a tetrahedron-edge in
  $\Delta^*$ from $e$ to $e_1$ with drift $g_1$ and
  a triangle-edge from $e_i$ to $e_{i+1}$ for all $i = 1, \ldots, k-1$ with
  drift $g_i^{-1}g_{i+1}$. Hence there is indeed a good path in $\Delta^*$ from
  $e$ to $e_k = e'$ with drift $g_k = g'$.

  If $I$ does intersect an edge of the veering triangulation aside from its endpoints, then we can perturb $I$ slightly so that it passes from one tetrahedron to the next through the interior of a face. From this sequence of tetrahedra, we obtain a good path with drift $g$ just like in the previous case. The good path we get depends on how the perturbation is done, but any of them works for our purposes.
  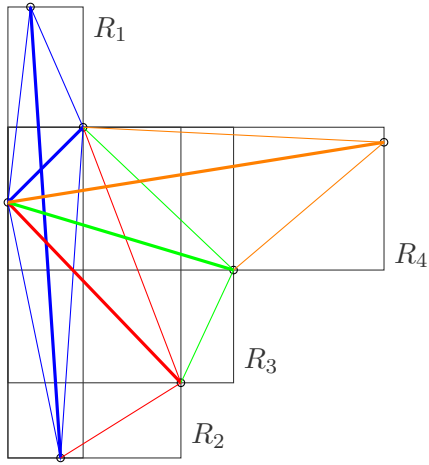
\begin{figure}
  \centering
  \begin{tikzpicture}
    \def\rad{0.05}
    \begin{scope}[black!80!white,thin]
    \draw (0,0) rectangle (1,6) node[below right] {$R_1$};
    \draw (0,0) rectangle (2.3,4.4);
    \node[above right] at (2.3,0) {$R_2$};
    \draw (0,1) rectangle (3,4.4);
    \node[above right] at (3,1) {$R_3$};
    \draw (0,2.5) rectangle (5,4.4);
    \node[above right] at (5,2.4) {$R_4$};
    \end{scope}
    \foreach \x/\y in {1/4.4, 0.3/6, 0/3.4, 0.7/0, 2.3/1, 3/2.5, 5/4.2}
    \draw (\x,\y) circle (\rad);
    \draw[blue] (1,4.4) -- (0.3,6) -- (0,3.4)
    -- (0.7,0) -- (1,4.4) -- (0,3.4); \draw[blue, very thick] (0.3,6) -- (0.7,0);
    \draw[blue, very thick] (1,4.4) -- (0,3.4);
    \draw[red] (0.7,0) -- (2.3,1) -- (1,4.4); \draw[red,very thick] (2.3,1) -- (0,3.4);
    \draw[green] (1,4.4) -- (3,2.5) -- (2.3,1); \draw[green,very thick] (3,2.5) -- (0,3.4);
    \draw[orange] (1,4.4) -- (5,4.2) -- (3,2.5); \draw[orange,very thick] (5,4.2) -- (0,3.4) ;
  \end{tikzpicture}
  \caption{The edges of the veering triangulation of $\tM$ corresponding to
    good path in $\Delta^*$. Every edge after the first one (the blue edge with
  the largest slope) is below the first edge.}
  \label{fig:good-path}
\end{figure}

For the other direction, let $\eps_1\eps_2\ldots\eps_n$ be a good path in
$\Delta^*$. Let $e_k$ be the endpoint of $\eps_k$ for all $k = 1, \ldots, n$
and let $e_0$ be the starting point of $\eps_1$. For all $k = 1, \ldots, n$,
let $R_k$ be the maximal singularity-free rectangle obtained by expanding the
smallest singularity-free rectangle containing
$\drift(\eps_1\ldots\eps_k)\te_k$ vertically as far as possible. Note that
$\te_0$ and $\drift(\eps_1)\te_1$ are the top and bottom edges of the
tetrahedron corresponding to $R_1$. In particular, $\te_0$ intersects $R_1$ at
its two horizontal sides. Each $R_{k+1}$ is shorter and wider than $R_k$, so
by induction, we see that $\te_0$ also intersects $R_n$ at its two horizontal
sides (\Cref{fig:good-path}). Therefore $\drift(\eps_1\ldots\eps_n)\te_n$,
as the arc connecting the vertical sides of $R_n$, is indeed below $\te_0$. Hence
$\drift(\eps_1\ldots\eps_n) \in \stash(e_0,e_n)$.
\end{proof}

The drift of a path in $\Delta^*$ is independent of the order of edges of a path.
So it will often be useful to think of a path in $\Delta^*$ as a nonnegative
integer-valued function on the edges of $\Delta^*$ where the value on every edge
is the number of times that edge appears in the path. This viewpoint allows us
to define the \emph{sum} of two paths by taking the sum of the corresponding
functions.

By a \emph{cycle} in the digraph $\Delta^*$, we mean a path $\eps_1\ldots\eps_n$
consisting of triangle-edges such that the starting point of $\eps_1$ coincides
with the endpoint of $\eps_n$. A \emph{minimal cycle} is a cycle that traverses
every vertex at most once. A \emph{minimal good path} in $\Delta^*$ is
a good path $\eps_1\ldots \eps_n$ such that the endpoints of $\eps_i$ are
pairwise distinct for $i=1,\ldots,n$. It is allowed, however, that the starting
point of $\eps_1$ coincides with one of the other vertices traversed. There are
finitely many minimal cycles and minimal good paths in~$\Delta^*$.

Let $\gamma = \eps_1\ldots \eps_n$ be a good path traversing the vertices
$e_0, \ldots, e_n$ and let $\gamma_1, \ldots, \gamma_k$ be cycles in $\Delta^*$.
We call the collection of paths $\gamma, \gamma_1, \ldots, \gamma_k$
\emph{connected} if
\begin{enumerate}
\item the triangle-edges appearing in these paths (that is, every edge other than
  $\eps_1$) form a connected subgraph of $\Delta^*$ when the orientations of the
  edges are ignored and
\item at least one cycle $\gamma_i$ traverses the vertex $e_1$ when $n=1$.
\end{enumerate}

We have the following decomposition lemma of good paths.

\begin{lemma}[Decompositions of good paths]\label{lemma:decomposing-good-paths}
  The sum of a connected collection of a minimal good path and a finite number of
  minimal cycles in $\Delta^*$ is a good path. Conversely, every good path in
  $\Delta^*$ can be written as such a sum.
\end{lemma}
\begin{proof}
  To prove the first statement we build up the sum step by step, adding one
  minimal cycle at a time. Denote by $\eta_0$ the minimal good path of the
  collection. There must be a minimal cycle from the collection that forms a connected union
  together with $\eta_0$. Their sum $\eta_1$ is a good path. Then, there must
  be another minimal good cycle from the collection that forms a connected union with
  $\eta_1$. Their sum is a good path $\eta_2$. Repeating this process until all
  cycles are added, we obtain the first statement.

  For the second statement, let $\gamma = \eps_1\ldots,\eps_n$ be a good path in $\Delta^*$. If it is
  minimal, we are done. If it is not minimal, then there are $1 \le i < j \le
  n$ such that the endpoints of $\eps_i$ and $\eps_j$ agree. Moreover, we can
  choose $i$ and $j$ so that $j-i$ is as small as possible. Then the subpath
  $\eps_{i+1}\ldots \eps_j$ is a minimal cycle and $\gamma$ can be written as a
  sum of this minimal cycle and a good path shorter than $\gamma$. We can
  repeat this process of removing minimal cycles until the remaining good path
  is minimal. It is straightforward to verify that the collection of summands
  is connected. Hence we obtain the second statement.
\end{proof}

An immediate corollary of \Cref{prop:intersection-and-good-paths} and
\Cref{lemma:decomposing-good-paths} is the following.

\begin{corollary}\label{cor:intersection-as-sum-of-drifts}
  Denote by $\calP_{e,e'} \subset G$ the set of drifts of minimal good paths from $e$ to $e'$
  and by $\calB \subset G$ be the set of drifts of the minimal cycles of
  $\Delta^*$. The set $\stash(e,e')$ is the set of products
  $p b_1^{\alpha_1} \ldots b_k^{\alpha_k}$ where $p \in \calP_{e,e'}$, $k\ge 0$
  is an integer, $b_1, \ldots, b_k \in \calB$, $\alpha_1, \ldots, \alpha_k$ are
  positive integers and $p, b_1, \ldots, b_k$ are drifts of a minimal good path
  and minimal cycles that form a connected collection.
\end{corollary}

We will use \Cref{cor:intersection-as-sum-of-drifts} to compute exact values of
the asymptotic translation length in the arc complex in \Cref{sec:example}. For
the proof of \Cref{thm:closure}, the following approximation of the stashing sets will be
more convenient.

\begin{corollary}\label{cor:intersection-easy-estimate}
  Denote by $\calP_{e,e'} \subset G$ the set of drifts of minimal good paths
  from $e$ to $e'$
  and by $\calB \subset G$ be the set of drifts of the minimal cycles of
  $\Delta^*$. Furthermore, let $\calP'_{e,e'} = \calP_{e,e'}\prod_{b\in \calB}b$.
  Then
  \begin{displaymath}
    \calP'_{e,e'} \langle \calB \rangle_{\ZZ_{\ge 0}} \subset \stash(e,e') \subset \calP_{e,e'} \langle \calB \rangle_{\ZZ_{\ge 0}},
  \end{displaymath}
  where $\langle \calB \rangle_{\ZZ_{\ge 0}}$ denotes the monoid generated by $\calB$.
\end{corollary}

By the product of two sets $X$ and $Y$, we mean
\begin{displaymath}
  XY = \{xy : x\in X,\, y\in Y\}.
\end{displaymath}

\begin{proof}
  The second containment is a trivial consequence of
  \Cref{cor:intersection-as-sum-of-drifts}. The first containment follows from
  \Cref{cor:intersection-as-sum-of-drifts} and the fact the union of any
  minimal good path from $e$ to $e'$ with all the minimal cycles is always a
  connected collection. This is because by
  \Cref{lemma:Delta-strongly-connected}, the graph of triangle-edges is
  strongly connected, so the union of all cycles or, equivalently, the union of
  all minimal cycles is a strongly connected graph containing all vertices.
\end{proof}

\section{Estimating the stashing sets}
\label{sec:cone-estimates}

\subsection{Monoids and cones}
\label{sec:monoids-and-cones}

We begin this section by proving some general lemmas. We will use these lemmas
to estimate the stashing sets at the end of the section.

For any $\calB \subset \RR^n$ and $\calE\subset \RR$, we introduce the notation
\begin{displaymath}
  \langle \calB \rangle_{\calE} = \{\sum_{i=1}^k\eta_ib_i: b_i \in \calB, \eta_i
  \in \calE\}
\end{displaymath}
for the set generated by $\calB$ with coefficients in $\calE$. (The empty sum is
allowed in the definition and it is defined to be zero.) For example,
$\langle \calB \rangle_{\RR_{\ge 0}}$ is the cone generated by $\calB$ and
$\langle \calB \rangle_{\ZZ_{\ge 0}}$ is the monoid generated by $\calB$.

\begin{lemma}\label{lemma:elementary-conelike}
  Let $\calB \subset \ZZ^n$ be a finite set and let
  $C = \langle \calB \rangle_{\RR_{\ge 0}}$ be the cone generated by $\calB$ in
  $\RR^n$. Then there exists some $x^* \in \ZZ^n$ such that
  \begin{displaymath}
    \langle \calB \rangle_{\ZZ_{\ge 0}} \cap (x^*+C)= \langle \calB \rangle_{\ZZ} \cap (x^*+C).
  \end{displaymath}
\end{lemma}
In words, the lemma says that the sets $\langle \calB \rangle_{\ZZ_{\ge 0}}$
and $\langle \calB \rangle_{\ZZ}$ are equal inside the translated cone
$x^* + C$. From this viewpoint, it is clear that the lemma also holds for any
element of the cone $x^*+C$ instead of $x^*$.
\begin{proof}
  The left hand side is clearly contained in the right hand side for any
  $x^*\in \ZZ^n$. We will find some $x^*\in \ZZ^n$ such that the reverse
  containment also holds. Let $\calB = \{b_1,\ldots, b_m\}$ and consider the
  compact subset $K = \{\sum_{i=1}^m\kappa_ib_i: 0\le \kappa_i \le 1\}$ of $C$.
  Each element of $K \cap \langle \calB \rangle_{\ZZ}$ can be represented in the form
  $\sum_{i=1}^m\eta_{i}b_i$ where $\eta_{i} \in \ZZ$. By choosing such an expression for each element, we may choose a positive integer $\eta^*$
  so that all $\eta_i$ that appear in these finitely many representations satisfy
  $-\eta^* \le \eta_i$.

  We claim that the reverse containment in the lemma holds for
  $x^* = \eta^*\sum_{i=1}^mb_i$. To see this, let $x \in \langle \calB \rangle_{\ZZ} \cap (x^*+C)$. Since
  $x \in x^*+C$, we have $x = \sum_{i=1}^m\alpha_{i}b_i$ for
  $\alpha_i \in \RR$ and $\alpha_i\ge \eta^*$. We can rewrite this
  representation of $x$ as
  \begin{displaymath}
    x = \sum_{i=1}^m[\alpha_i]b_i + \sum_{i=1}^m\{\alpha_i\}b_i
  \end{displaymath}
  where $[\alpha_i]$ and $\{\alpha_i\}$ denote the integer and fractional parts
  of $\alpha_i$. The first of the two terms on the right is in
  $\langle \calB \rangle_{\ZZ}$ and so is $x$, therefore the second term on the
  right is also in $\langle \calB \rangle_{\ZZ}$. It is also in $K$, therefore
  we can replace it by $\sum_{i=1}^m\eta_ib_i$ where $\eta_i\in \ZZ$ and
  $\eta_i \ge -\eta^*$. Since $[\alpha_i]\ge \eta^*$, we obtain a representation of $x$
  as a sum of the $b_i$ with nonnegative integer coefficients. Therefore
  $x \in \langle \calB \rangle_{\ZZ_{\ge 0}}$ and the right hand side in the
  lemma is indeed contained in the left hand side.
\end{proof}

\begin{lemma}\label{lemma:complicated-conelike-fixed}
  Let $n\ge 2$ and let $\calB \subset \ZZ^n$ be a finite set. Let
  $\calD \subset \Hom(\RR^n,\RR)$ such that the cone
  $\langle \calD \rangle_{\RR_{\ge 0}}$ has nonempty interior in the
  $n$-dimensional vector space $\Hom(\RR^n,\RR)$. Assume that
  $\Frob(\phi|_{\langle \calB \rangle_{\ZZ_{\ge 0}}}) < \infty$ for every
  primitive integral point $\phi$ in $\langle \calD \rangle_{\RR_{\ge 0}}$.
  Then the following statements hold:
  \begin{enumerate}[(i)]
      \item\label{item:not-proper-subgroup} $\langle \calB \rangle_{\ZZ} = \ZZ^n$.
      \item\label{item:full-dimensional-cone} The cone $\langle \calB \rangle_{\RR_{\ge 0}}$ has nonempty interior.
      \item\label{item:cone-contains-all-integral-points} There exists $x \in \ZZ^n$ such that
    $\ZZ^n \cap (x + \langle \calB \rangle_{\RR_{\ge 0}}) \subset \langle
    \calB \rangle_{\ZZ_{\ge 0}}$.
  \end{enumerate}
\end{lemma}
\begin{proof}

  We begin by proving (\ref{item:not-proper-subgroup}).
  If $\langle \calB \rangle_{\ZZ}$ was a proper subgroup of $\ZZ^n$, then there would exist a surjective homomorphism $\phi\co \ZZ^n \to \ZZ$ such that $ \phi(\langle \calB \rangle_{\ZZ})$ is not all of $\ZZ$. Since $\phi$ is surjective, it is a primitive integral point in $\Hom(\RR^n, \RR)$.
  
  Let $\phi_0$ be any primitive integral
    point in the interior of $\langle \calD \rangle_{\RR_{\ge 0}}$ that is not
    a scalar multiple of $\phi$. For a large enough positive integer $N$, $N\phi_0 + \phi$ is in
    the interior of $\langle \calD \rangle_{\RR_{\ge 0}}$. We may choose a
    basis for the 2-dimensional lattice obtained as the intersection of
    $\Hom(\ZZ^n,\ZZ)$ with the 2-dimensional subspace spanned by $\phi_0$ and
    $\phi$ in $\Hom(\RR^n,\RR)$ such that $\phi_0$ has coordinates $(0,1)$. Let
    $(q,r)$ be the coordinates of $\phi$. Note that $q\ne 0$ and $q$ and $r$ are relatively prime. From this, we see that for any
    integer $N$, the point 
    $\psi_N = Nq\phi_0 + \phi$ is primitive since it has coordinates $(q, r+Nq)$.  Since $\phi(\langle \calB \rangle_{\ZZ})$ is not all of $\ZZ$, there is some integer $d\ge 2$ that divides every element of this image. Choosing $N=ad$ for some large positive integer $a$, we can see that all elements of $\psi_{ad}(\langle \calB \rangle_{\ZZ})$ are divisible by $d$.
    But this contradicts the fact
    that $\Frob(\psi_{ad}|_{\langle \calB \rangle_{\ZZ_{\ge 0}}}) < \infty$.

The statement (\ref{item:full-dimensional-cone}) is a straightforward corollary of (\ref{item:not-proper-subgroup}).

    By \Cref{lemma:elementary-conelike}, there is some $x \in \ZZ^n$
    such that
    \begin{displaymath}
      \langle \calB \rangle_{\ZZ_{\ge 0}} \cap (x+\langle \calB \rangle_{\RR_{\ge 0}})= \langle \calB \rangle_{\ZZ} \cap (x+\langle \calB \rangle_{\RR_{\ge 0}}).
    \end{displaymath}
    Using that $\langle \calB \rangle_{\ZZ} = \ZZ^n$ from (\ref{item:not-proper-subgroup}) and that the left hand side is contained in $\langle \calB \rangle_{\ZZ_{\ge 0}}$, we obtain (\ref{item:cone-contains-all-integral-points}).
\end{proof}

\begin{lemma}\label{lemma:good-loops}
  For all $e \in E$ we have $\stash(e,e) \subset \langle \calB \rangle_{\ZZ_{\ge 0}}$.
\end{lemma}
\begin{proof}
  Observe that
  for any tetrahedron-edge of $\Delta^*$ from $e$ to $e'$, there is a path of
  triangle edges from $e$ to $e'$ with the same drift. To see this, choose a
  flow line close to $e$ that does not intersect any edges of the veering
  triangulation. This flow line has a subarc that starts at the tetrahedron
  whose bottom edge is $e$ and ends at the tetrahedron whose bottom edge is
  $e'$ and whose top edge is $e$. Between the two tetrahedra, the flow line
  intersects a sequence of tetrahedra. This sequence defines a path from $e$ to
  $e'$ with the required properties. As a consequence, for any good path in $\Delta^*$ from starting and ending at the
  same vertex, there is a cycle in $\Delta$ with the same drift. So the statement follows by
  \Cref{prop:intersection-and-good-paths}.
\end{proof}

We are now equipped with the tools to prove the following statement, describing
the structure of the stashing sets.

\begin{proposition}\label{prop:automaton-and-intersecting-edges}
  Let $\calB \subset G$ be the set of drifts of the minimal cycles of
  $\Delta^*$. Then the cone $\langle \calB \rangle_{\RR_{\ge 0}}$ has nonempty
  interior and there exists $g \in G$ such that
  \begin{displaymath}
    G \cap g \langle \calB \rangle_{\RR_{\ge 0}} \subset \langle
    \calB \rangle_{\ZZ_{\ge 0}}.
  \end{displaymath}
  Moreover, for every $e,e' \in E$, there exist $g_1,g_2 \in G$ such that
  \begin{displaymath}
    G \cap g_1\langle \calB \rangle_{\RR_{\ge 0}} \subset \stash(e,e') \subset
    G \cap g_2\langle \calB \rangle_{\RR_{\ge 0}}.
  \end{displaymath}
\end{proposition}
\begin{proof}
  By \Cref{cor:frob-exists-for-intersection-data},
  \begin{equation}\label{eq:finite-frob}
    \Frob(\phi|_{\stash(e,e)}) < \infty
  \end{equation}
  for every primitive integral point $\phi$ in $-\RR_+\calF$. If we replace $\stash(e,e)$
  with the set $\langle \calB \rangle_{\ZZ_{\ge 0}}$ (which is larger by \Cref{lemma:good-loops}) in \Cref{eq:finite-frob}, the statement remains true.
  Therefore we may use \Cref{lemma:complicated-conelike-fixed} with $G\cong \ZZ^n$ and
  any generator set $\calD$ for the cone $-\RR_+\calF$ to obtain that $\langle \calB \rangle_{\RR_{\ge 0}}$ has
  nonempty interior and that there exists $g \in G$ such that
  \begin{equation}\label{eq:B-cone-contains-all-integral-points}
    G \cap g \langle \calB \rangle_{\RR_{\ge 0}} \subset \langle
    \calB \rangle_{\ZZ_{\ge 0}}.
  \end{equation}
    
  Moreover, we obtain that there are
  $g_1, g_2 \in G$ such that
  \begin{displaymath}
    G \cap g_1\langle \calB \rangle_{\RR_{\ge 0}} \subset \calP'_{e,e'} \langle \calB \rangle_{\ZZ_{\ge 0}} \subset \stash(e,e')
    \subset \calP_{e,e'} \langle \calB \rangle_{\ZZ_{\ge 0}} \subset G \cap
    g_2\langle \calB \rangle_{\RR_{\ge 0}},
  \end{displaymath}
  where the first containment follows from (\ref{eq:B-cone-contains-all-integral-points}), the second and third containments were shown in \Cref{cor:intersection-easy-estimate}, and the last containment follows from the simple observation that $g_2$ can be chosen such that $g_2^{-1} \calP_{e,e'} \subset \langle \calB \rangle_{\RR_{\ge 0}}$ since $\langle \calB \rangle_{\RR_{\ge 0}}$ has non-empty interior.
 \end{proof}

\subsection{Duality of cones}
\label{sec:duality}

The goal of this section is to prove \Cref{prop:cone-description} below that
states that the cone over the fibered face consists of precisely those
cohomology classes that take non-positive values on the cone $\langle \calB \rangle_{\RR_{\ge 0}}$.

\begin{proposition}\label{prop:cone-description}
  Let $\calB \subset G$ be the set of drifts of the minimal cycles of
  $\Delta^*$. Then the interior of the cone
  $\RR_+\calF \subset H^1(M;\RR)$ can be described as
  \begin{displaymath}
    \{\phi \in H^1(M;\RR): \phi(x) < 0\mbox{ for all }x \in \langle \calB \rangle_{\RR_{\ge 0}}\}.
  \end{displaymath}
\end{proposition}

We remark that an analogous statement was proven by Fried in \cite[Theorem
D]{Fried82a}. In that paper, Fried has defined the set of \emph{homology
  directions} of a flow on an $n$-dimensional closed manifold and showed that the
integral cohomology classes that correspond to a fibration of the manifold over
the circle are exactly the ones that take positive values on the set of
homology directions. Unfortunately, Fried's proof assumes that the manifold is
closed, so the theorem cannot be directly applied for our case. In the end of
the introduction, Fried mentions that under appropriate hypotheses the results
carry over also to compact manifolds by doubling the manifold along the
boundary, but details are not given. In order to make the proof of
\Cref{prop:cone-description} as transparent as possible, instead of extending
Fried's theorem to the non-closed case and relating the set of homology
directions to our cone $\langle \calB \rangle_{\RR_{\ge 0}}$, we 
give a direct proof for \Cref{prop:cone-description} following Fried's
strategy but in the combinatorial spirit of this paper.

Before giving the proof of \Cref{prop:cone-description}, we prove a few brief lemmas.

Given \emph{any} nonzero $\phi \in H^1(M;\ZZ)$, not necessarily in the fibered
cone, consider the infinite cyclic covering $M_\phi \to M$ corresponding of
$\phi$. This covering induces an infinite cyclic covering
$\Delta_\phi \to \Delta$ of the graph $\Delta$ modeling the veering
triangulation of $M$. We define the drift of each edge of $\Delta_\phi$ as the
drift of the its projection in $\Delta$.

Define an integer-valued function on the vertices of $\Delta_\phi$ as follows.
By associating to each tetrahedron of the veering triangulation its bottom
edge, each tetrahedron in the veering triangulation of $\tM$ can be referred to
as $g\te$ for some $g \in G$ and $e \in E$. Since $M_\phi$ is a quotient of
$\tM$ where two edges $g_1\te$ and $g_2\te$ have the same image if and only if
$\phi(g_1) = \phi(g_2)$, the integer $\phi(g)$ is a well-defined invariant of
the image of any edge $g\te$ in $M_\phi$. This way we obtain an integer
associated to each vertex $v$ of $\Delta_{\phi}$ which we will denote by
$\phi(v)$.

\begin{lemma}\label{lemma:phi-difference}
  For any nonzero $\phi \in H^1(M;\ZZ)$, there exists some $Q > 0$ such
  that if $\gamma$ is a path in $\Delta_\phi$ starting at $v$ and ending at
  $v'$, then
  \begin{displaymath}
    \phi(v') - \phi(v) = q + \sum_{i=1}^k \phi(b_i)
  \end{displaymath}
  for some $q\in \ZZ$ with $|q|\le Q$ and some $k\ge 0$ and $b_i \in \calB$.
\end{lemma}
\begin{proof}
  It is straightforward to verify from the definition of $\phi(v)$ and
  $\phi(v')$ that
  \begin{displaymath}
    \phi(v') - \phi(v) = \phi(\drift(\gamma)) = \phi(\drift(\pi(\gamma)))
  \end{displaymath}
  where $\pi(\gamma)$ is the projection of $\gamma$ in the graph $\Delta$ with
  finitely many vertices. As in \Cref{lemma:decomposing-good-paths}, we can
  decompose $\pi(\gamma)$ as a sum of minimal cycles and a path that does not
  contain a cycle. By setting
  \begin{displaymath}
    Q = \max\{|\phi(\drift(\delta))|: \delta \mbox{ is a
      path in
      $\Delta$ containing no cycles}\},
  \end{displaymath}
  we obtain the statement of the lemma.
\end{proof}

\begin{lemma}\label{lemma:inf-to-neg-inf}
  Let $\phi \in H^1(M;\ZZ)$ such that $\phi(b) < 0$ for all $b \in \calB$. If
  $\cdots \to v_{-1} \to v_0 \to v_1 \to \cdots$ is a bi-infinite path in
  $\Delta_{\phi}$, then $\lim_{n \to \infty} \phi(v_n) = -\infty$ and
  $\lim_{n \to -\infty} \phi(v_n) = \infty$.
\end{lemma}
\begin{proof}
  Using \Cref{lemma:phi-difference} and its notation, we have
  \begin{displaymath}
    \lim_{n\to \infty} \phi(v_n) - \phi(v_0) = \lim_{n\to \infty} q_n +
    \sum_{i=1}^{k_n} \phi(b_{i,n})
  \end{displaymath}
  where $q_n \in \ZZ$ with $|q_n| \le Q$ and $b_{i,n}\in \calB$. Since the
  $\phi(b_{i,n})$ are negative integers and $\lim_{n\to \infty} k_n = \infty$,
  we have $\lim_{n \to \infty} \phi(v_n) = -\infty$. The proof of the limit as
  $n\to -\infty$ is analogous.
\end{proof}

\begin{lemma}\label{lemma:every-large-every-small}
  Let $\phi \in H^1(M;\ZZ)$ such that $\phi(b) < 0$ for all $b \in \calB$. Let $v_0$ be a vertex of
  $\Delta_\phi$ and let $V_+$ be the set of vertices (including $v_0$)
  that are endpoints of a path starting at $v_0$. Then there exist $N_1,N_2 \in
  \ZZ$ such that
  \begin{displaymath}
    \{v\in \Delta_\phi: \phi(v) \le N_1\} \subset V_+ \subset \{v\in \Delta_\phi: \phi(v) \le N_2\}.
  \end{displaymath}
\end{lemma}
\begin{proof}
  If $v \in V_+$, then by \Cref{lemma:phi-difference}, we have
  \begin{displaymath}
    \phi(v') = \phi(v_0) + q + \sum_{i=1}^{k} \phi(b_{i})
  \end{displaymath}
  where $q \in \ZZ$ with $|q| \le Q$ and $b_{i}\in \calB$. Hence the second
  containment in the lemma holds with $N_2 = \phi(v_0) + Q$. For the first
  containment, observe that there is some $N_1 < 0$ such that every integer
  less than $N_1$ can be written in the form $\sum_{i=1}^{k} \phi(b_{i})$. This
  follows from the fact that the cone $\langle \calB \rangle_{\RR_{\ge 0}}$ has
  nonempty interior and that the 
  monoid $\langle \calB \rangle_{\ZZ\ge 0}$ contains every integral point in
  some translate of $\langle \calB \rangle_{\RR_{\ge 0}}$
  (\Cref{prop:automaton-and-intersecting-edges}).
\end{proof}

\begin{proof}[Proof of \Cref{prop:cone-description}]
  First we will show that if $\phi$ is a primitive integral point in the
  interior of $\RR_+\calF$, then $\phi(g) < 0$ for all
  $g \in \langle \calB \rangle_{\RR_{\ge 0}}$. It suffices to show this for all
  $g \in \calB$. Let $\gamma$ be a cycle in $\Delta$ with drift $g$. Corresponding to the cycle
  $\gamma$ is a sequence of tetrahedra $T_0, \ldots, T_m$ in the veering triangulation of $\tM$ such that
  for all $i=1, \ldots, m$, the tetrahedra $T_{i-1}$ and $T_{i}$ share a face
  and $T_{i}$ is below $T_{i-1}$. Moreover, $T_m = gT_0$. Therefore
  multiplication by $g$ translates $T_0$ to a tetrahedron
  below it.

  By convention (see \Cref{sec:fibrations-over-circle}) the cohomology class
  $\phi$ evaluates to positive integers on loops of $M$ whose lift ``goes up''
  (the endpoint of the lift is higher than the starting point) in the infinite
  cyclic cover $S\times \RR \to M$ corresponding to $\phi$. As we see from the
  tetrahedron sequence, loops representing $g$ lift to paths that ``go
  down'' in $S \times \RR$. Therefore $\phi(g) < 0$ indeed.

  Consider the (open) cone
  \begin{displaymath}
    D = \{\phi \in H^1(M;\RR): \phi(g) < 0 \mbox{ for all }g\in \langle \calB \rangle_{\RR_{\ge 0}}\}\subset H^1(M;\RR).
  \end{displaymath}
  What we have just proved implies that the interior of $\RR_+\calF$ is
  contained in $D$. To prove the proposition, we need to prove that,
  conversely, $D$ is contained in the interior of $\RR_+\calF$. If this was not
  true, then $D$ would contain a primitive integral class on the boundary of
  $\RR_+\calF$. (The boundary faces of the cone $\RR_+\calF$ are defined by
  rational equations, so primitive integral points are projectively dense on
  the boundary of $\RR_+\calF$.) Since primitive integral classes on the
  boundary of $\RR_+\calF$ are known not to correspond to fibrations, it
  suffices to show that if $\phi \in D$ is a primitive integral class, then
  $\phi$ is dual to a fibration.

  A \emph{cut} of $\Delta_\phi$ is a way of dividing the vertices of
  $\Delta_\phi$ into two disjoint nonempty sets $V_{-\infty}$ and $V_{\infty}$
  that are closed under ``going forward'' and ``going backward'', respectively.
  More precisely, if there is an edge from $v_1$ to $v_2$ in $\Delta_\phi$ then
  $v_1 \in V_{-\infty}$ implies $v_2 \in V_{-\infty}$ and $v_2 \in V_{\infty}$
  implies $v_1 \in V_{\infty}$.

  To see that cuts exists, let $v$ be a vertex of $\Delta_{\phi}$ and let
  $V_{-\infty}$ be the set of vertices (including $v$) that are endpoints of a
  path starting at $v$ and let $V_{\infty}$ be the set of the remaining
  vertices. It is clear that $V_{-\infty}$ and $V_{\infty}$ are closed under
  going forward and going backward, respectively. It follows from
  \Cref{lemma:inf-to-neg-inf,lemma:every-large-every-small} that both
  $V_{-\infty}$ and $V_{\infty}$ are nonempty.

  Next, we associate an embedded surface in $M_\phi$ to each cut. Given a cut
  $V_{-\infty} \cup V_{\infty}$, let $\Sigma$ be the union of triangles of the
  veering triangulation corresponding to the edges starting at a point of
  $V_\infty$ and ending at a point of $V_{-\infty}$. To show that $\Sigma$ is
  a surface, we need to prove that there are two triangles meeting at every
  edge. (We are gluing together ideal triangles---their vertices are not part of the 3-manifold---, therefore we do not need to check that the links of the vertices are circles.) Let $e$ be an edge of the veering triangulation of $M_\phi$ and let $T$
  and $T'$ be the tetrahedra whose bottom and top edges are $e$, respectively.
  The tetrahedra adjacent to $e$ define an immersed subgraph $\Gamma_e$ of
  $\Delta_\phi$ with the structure shown on \Cref{fig:graph_around_edge}.

  \begin{figure}[htb]
    \labellist
    \small\hair 2pt
    \pinlabel {$T$} [ ] at 15 31
    \pinlabel {$T'$} [ ] at 291 29
    \endlabellist
    \centering
    \includegraphics[scale=0.7]{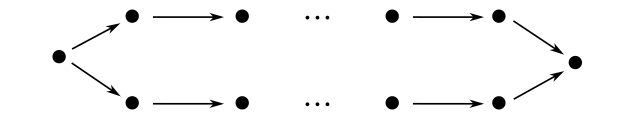}
    \caption{The immersed subgraph $\Gamma_e$ in $\Delta$ whose vertices
      correspond to the tetrahedra adjacent to $e$.}
    \label{fig:graph_around_edge}
  \end{figure}

  Observe that either 
  \begin{enumerate}
      \item all vertices of $\Gamma_e$ are in $V_{-\infty}$, or
      \item all vertices of $\Gamma_e$ are in $V_\infty$, or
      \item $T \in V_\infty$, $T' \in V_{-\infty}$ and
  exactly two edges of $\Gamma_e$ start in $V_\infty$ and end in $V_{-\infty}$, with one
  edge on each of the two paths from $T$ to $T'$ in $\Gamma_e$.
  \end{enumerate}
  Hence there are
  indeed either zero or two triangles meeting at every edge and $\Sigma$ is an
  embedded surface.

  Next, observe that each flow line in $M_\phi$ intersects $\Sigma$ at exactly
  one point. This is because the tetrahedra intersected by the flow line give
  rise to a bi-infinite path $\cdots \to v_{-1} \to v_0 \to v_1 \to \cdots$ in
  $\Delta_\phi$. By \Cref{lemma:inf-to-neg-inf} and
  \Cref{lemma:every-large-every-small}, there exists some $i_0 \in \ZZ$ such
  that $v_i \in V_{-\infty}$ if $i \ge i_0$ and $v_i \in V_{\infty}$ otherwise.
  Therefore the flow line intersects exactly one triangle of $\Sigma$: the one
  corresponding to the edge from $v_{i_0-1}$ to $v_{i_0}$. When the flow line
  intersects some edges of the veering triangulation, then the corresponding
  bi-infinite path is not unique, but it is straightforward to verify that such
  flow lines also intersect $\Sigma$ in one point.

  As a corollary, we obtain a homeomorphism $\Sigma \times \RR \to M_\phi$
  defined by the formula $(x,t) \to g_t(x)$ where $g_t$ denotes the flow on
  $M_\phi$.

  Let $h\co M_\phi \to M_\phi$ be the generator of the deck group of the
  covering $M_\phi \to M$ such that $\phi(h(v)) = \phi(v) - 1$ for every vertex
  $v$ of $\Delta_\phi$. Our final step is to replace $\Sigma$ with a homotopic
  surface $\Sigma'$ in $M_\phi$ such that $h(\Sigma')$ is disjoint from and
  homotopic to $\Sigma'$. This will show that the covering $M_\phi \to M$ comes
  from a fibration.

  Every surface $\Sigma'$ in $M_\phi$ intersecting every flow line once can be
  represented by a continuous function $u\co \Sigma \to \RR$ such that
  $\Sigma' = \{g_{u(x)}(x): x \in \Sigma\}$. For example, the function
  corresponding to $\Sigma$ is the constant zero function.

  Let $n$ be a positive integer and consider the surfaces $\Sigma$,
  $h(\Sigma)$, \ldots, $h^{n}(\Sigma)$ that all correspond to cuts of
  $\Delta_\phi$, therefore intersect every flow line once. Let
  $u_0, \ldots, u_{n} \co \Sigma \to \RR$ be the corresponding functions. Let
  $u' = \frac{u_0 + \cdots +u_{n-1}}{n}$ and let $\Sigma'$ be the corresponding
  surface. The function corresponding to $h(\Sigma')$ is
  $\frac{u_1 + \cdots + u_{n}}{n}$, which is strictly larger than $u'$ at every
  point of $\Sigma$ if $n$ is large enough, since $u_n > u_0$ if $n$ is large
  enough. Therefore $\Sigma'$ is an embedded surface intersecting every
  flowline exactly once such that $h(\Sigma')$ is homotopic to and disjoint
  from $\Sigma'$. Hence $M_\phi \to M$ comes from a fibration, and that is what
  we wanted to show.
\end{proof}

\section{Lemmas on cones, lattices and volumes}
\label{sec:frob-lemmas}

\newcommand\be{\mathbf{e}}
\newcommand\one{\mathbf{1}}
\newcommand\occ{\operatorname{occ}}
\newcommand\occs{\operatorname{occs}}

This section contains various lemmas on cones, lattices and volumes in
Euclidean spaces that will be used to prove the main theorems. All results in
this section are self-contained and independent of 3-manifold theory.

\subsection{Occupancy coefficients}
\label{sec:occ-coeff}

Let $V$ be an $n$-dimensional real vector space, let $K \subset V$ be a compact
set with nonempty interior and let $\Lambda \subset V$ be a lattice.
The \emph{occupancy coefficient} $\occ(\Lambda, K)$ of $K$ with respect to
the lattice $\Lambda$ is the ratio
\begin{equation}\label{eq:occ-def}
  \occ(\Lambda, K) = \frac{\vol(K')}{\vol(V/\Lambda)}
\end{equation}
where $K' = a + bK$ ($a, b\in \RR$) is a set of maximal volume that is obtained from $K$
by dilatation and translation and does not contain any point of $\Lambda$ in its
interior.

It seems difficult to compute occupancy coefficients in general, but some
basics facts can easily be deduced.

\begin{lemma}\label{lemma:occ-dim-one}
  The occupancy coefficient of a connected set in a 1-dimensional vector space
  equals 1.
\end{lemma}
\begin{proof}
  Identifying $V$ by $\RR$, we have $\Lambda = a\ZZ$ for some $a > 0$. Our
  connected set $K$ is an interval. The longest interval $K'$ that does not
  contain a point of $\Lambda$ in its interior has length $a$. By
  \Cref{eq:occ-def}, we have $\occ(\Lambda, K) = \frac{a}{a} = 1$.
\end{proof}

\begin{lemma}\label{lemma:occ-lower-bound}
  The occupancy coefficient is always at least 1.
\end{lemma}
\begin{proof}
  If a set $K \subset V$ has volume less than $\vol(V/\Lambda)$, then its image
  in $n$-torus $V/\Lambda$ is not everything, therefore there is a translate of
  $K$ that is disjoint from $\Lambda$. So $\vol(K') \ge \vol(V/\Lambda)$ for
  the set $K'$ with maximal volume and therefore the occupancy coefficient is
  at least 1.
\end{proof}

\subsection{Occupancy coefficients as the lattice is varied}
\label{sec:min-occ}

For any compact set $K \subset V$ with nonempty interior, introduce the
notation
\begin{displaymath}
  \min\occ(K) = \inf_{\Lambda\subset V} \occ(\Lambda, K)
\end{displaymath}
where $\Lambda$ ranges over the lattices in $V$. By
\Cref{lemma:occ-lower-bound}, $\min\occ(K) \ge 1$ holds for all $K$.

\begin{lemma}\label{lemma:occs}
  Let $K \subset V$ be a compact connected set with nonempty interior in an
  $n$-dimensional vector space $V$. Consider the set
  \begin{equation}\label{eq:occs-def}
    \occs(K) = \{\occ(\Lambda, K): \Lambda \subset V \mbox{ is a lattice}\}.
  \end{equation}
  If $n=1$, then $\occs(K) = \{1\}$. If $n \ge 2$, then $\occs(K)$ is a
  half-infinite interval whose left endpoint is $\min\occ(K)$.
\end{lemma}
\begin{proof}
  The $n=1$ case follows from \Cref{lemma:occ-dim-one}.

  The space of lattices is connected and the occupancy coefficient is a
  continuous function, so $\occs(K)$ is an interval. It is clear that the left
  endpoint of this interval is $\min\occ(K)$. It remains to show that $\Lambda$
  can be chosen so that $\occ(\Lambda, K)$ is arbitrarily large when $n\ge 2$.

  Let $\be_1, \ldots, \be_n$ be a basis for $V$ and consider the sequence of
  lattices generated by $c\be_1, c\be_2, \ldots, c\be_{n-1}, \be_{n}$ as
  $c \to 0$. The covolumes of these lattices go to zero. However, a dilated and translated
  copy of $K$ that lies between the hyperplanes $\Sigma$ and $\be_n + \Sigma$
  where $\Sigma$ is the hyperplane spanned by $\be_1, \ldots, \be_{n-1}$ is
  disjoint from all these lattices, hence $\vol(K')$ in \Cref{eq:occ-def} is
  bounded from below as $c\to \infty$. So indeed, the occupancy coefficient can be arbitrarily
  large.
\end{proof}

\subsection{The main technical lemma on Frobenius numbers}
\label{sec:main-technical-lemma}

The following technical lemma is at the heart of the proof of \Cref{thm:closure}.

\begin{lemma}\label{lemma:main-frob-bound}
  Let $\Lambda$ be a lattice in an $n$-dimensional real vector space $V$. Let
  $C =\langle \calB \rangle_{\RR_{\ge 0}}$ be a cone with nonempty
  interior, generated by a finite set $\calB \subset V$. Let
  $e_0$ be a point in the interior of $C$ and let $x_1, x_2 \in V$ be
  arbitrary. Then there is a constant $K = K(\Lambda, C, e_0, x_1, x_2) > 0$ such
  that the following holds.

  Let $\Lambda' \subset \Lambda$ be such that
  \begin{equation}\label{eq:containment}
    \Lambda \cap (x_1 + C) \subset \Lambda' \subset \Lambda \cap (x_2 + C).
  \end{equation}
  Let $\beta : V \to \RR$ be a linear function with $\beta(e_0) = 1$ that
  takes rational values on $\Lambda$ and positive values on $C-\{0\}$. Let
  $\bar\beta$ be the unique positive scalar multiple of the $\beta$ such that
  $\bar\beta(\Lambda) = \ZZ$. Let $P$ be a polytope that is the intersection
  of the hyperplane $\beta^{-1}(0)$ and $y-C$ for some $y\in V$ with
  $\beta(y) > 0$. Then

  \begin{displaymath}
    \left| \Frob(\bar\beta|_{\Lambda'}) -
      \sqrt[n-1]{\frac{\occ(\Lambda \cap \beta^{-1}(0),P)\cdot \vol(\RR^n/\Lambda)}{n\vol(C \cap \beta^{-1}([0,1]))}}
      \bar\beta(e_0)^{1+\frac{1}{n-1}} \right| \le K \bar\beta(e_0).
  \end{displaymath}

\end{lemma}
\begin{proof}
  Let $P_0$ be a polytope in the hyperplane $\beta^{-1}(0)$, obtained from
  $P$ by a dilatation and translation whose interior does not contain any point of
  $\Lambda$ and whose $n-1$-dimensional volume is the maximal with respect to
  this property. Let $y_0 \in V$ such that
  $P_0 = \beta^{-1}(0) \cap (y_0 - C)$. Let
  $\Lambda_{0} = \beta^{-1}(0) \cap \Lambda$ be the lattice in the hyperplane
  $\beta^{-1}(0)$.

  \setcounter{step}{0}
  \begin{step}[Upper bound on the Frobenius number] The following inequalities hold:
    \begin{equation}\label{eq:frob-upper-bound}
      \begin{aligned}
        \Frob(\bar\beta|_{\Lambda'}) &\le \Frob(
        \bar\beta|_{\Lambda \cap (x_1+C)})  \\
        &= \Frob(\bar\beta|_{\Lambda \cap (x_1+\Lambda_0+C)})
         \\
        &\le \Frob(\bar\beta|_{\Lambda \cap \{z\in \RR^n: \beta(z) \ge
          \beta(x_1+y_0)\}}) \\ &< \bar\beta(x_1+y_0).
      \end{aligned}
    \end{equation}
    The first inequality follows from the containment
    $\Lambda \cap (x_1+C) \subset \Lambda'$. The equality holds, because
    $\bar\beta(\Lambda_0) = 0$. For the second inequality, note that if $z$
    satisfies $\beta(z) > \beta(y_0)$, then $z \in C + \Lambda_0$. This is
    because the polytope $\beta^{-1}(0) \cap (z - C)$ is a scaled-up copy of $P_0$, so
    it contains some $x^* \in \Lambda_0$ in its interior and therefore
    $z \in x^* + C \subset \Lambda_0 + C$. So if $\beta(z) > \beta(y_0 + x_1)$,
    then $z \in x_1 + \Lambda_0 + C$. Finally, the last inequality simply
    follows from the definition of the Frobenius number.

  \begin{figure}[ht]
    \centering \begin{tikzpicture}
  \clip (-3,-1.7) rectangle (9,1.5);
  \filldraw[thick,blue!20!white] (0,0.72) -- (3.6,0) -- (0,0);
  \filldraw[thick,red!20!white] (1,0) -- (6,-0.96) -- (6,0);
  \filldraw[blue!40!white] (0,0) rectangle (0.6,0.6);
  \filldraw[red!40!white] (6,0) rectangle +(-0.8,-0.8);
  \foreach \i in {-10,...,15}{
    \foreach \j in {-10,...,15}
    \filldraw (\i*0.2 + \j*1.033, \i*0.2-\j*0.166) circle (0.01);
  }
  \filldraw (0,0) node[below] {$0$} circle (0.03) -- (0,4);
  \draw (0,0) -- (10,0);
  \draw[fill=black] (0.6,0.6) circle (0.03);
  \node[fill=white,inner sep=1pt,rounded corners,above right] at (0.6,0.6) {$e_0$};
  \draw[very thick,green!70!black] (-7,1.4) --  (10,-2);
  \node[fill=white,inner sep=1pt,below=0.1cm] at (-1.5,0.3) {$\beta^{-1}(0)$};
  \node[red!50!black] at (3.9,-0.3) {$U$};
  \node[blue!20!black] at (0.3,0.3) {$Q$};
  \draw (6,-1.2) -- +(0,4);
  \filldraw (6,-1.2) circle (0.03) -- +(4,0);
  \node[fill=white,inner sep=1pt,above right=0.05cm] at (6,0) {$y_0$};
  \filldraw (5.57,-0.45) circle (0.03) node[left] {$p$};
\end{tikzpicture}
    \caption{}
    \label{fig:triangles}
  \end{figure}
\end{step}

\begin{step}[Lower bound on the Frobenius number]
  Let $Q$ be the polytope that is the intersection of $C$ and $e_0-C$ (see
  \Cref{fig:triangles}). Let $m = m(C,e_0,\Lambda) >0$ be a number such that
  any translate of $mQ$ contains some point of the lattice $\Lambda$ in its
  interior. We claim that
  \begin{equation}\label{eq:frob-lower-bound}
    \Frob(\bar\beta|_{\Lambda'}) \ge \Frob(\bar\beta|_{\Lambda \cap (x_2+C)}) > \bar\beta(y_0+x_2) - m\bar\beta(e_0).
  \end{equation}

  The first inequality is a consequence of the containment
  $\Lambda' \subset\Lambda \cap (x_2+C)$. To see why the second inequality
  holds, consider the polytope
  $U = \beta^{-1}([\beta(y_0)-m,\beta(y_0)]) \cap (y_0-C)$ that contains the
  polytope $y_0-mQ$.

  The translate $x_2 + U$ of $U$ contains a point $p \in \Lambda$ in its
  interior. The side of $x_2 + U$ opposite to $x_2 + y_0$ is contained in the
  level set $\beta^{-1}(\beta(x_2+y_0)-m)$. Therefore
  $\beta(p) > \beta(y_0+x_2) - m$ and
  \begin{displaymath}
    \bar\beta(p) > \bar\beta(y_0+x_2) - m\bar\beta(e_0).
  \end{displaymath}
  The second inequality in (\ref{eq:frob-lower-bound}) now follows from putting
  this together with the inequality
  $\Frob(\bar\beta|_{\Lambda \cap (C+x_2)}) \ge \bar\beta(p)$ which holds
  because the elements of $\Lambda$ on which $\bar\beta$ takes the value
  $\bar\beta(p)$ are exactly the points of $p + \Lambda_0$, none of which are
  contained in $C+x_2$, since the set $U + \Lambda_0$ is disjoint from $C+x_2$.
  (\Cref{fig:triangles} shows the case $x_2 = 0$.)
\end{step}

Let $Y_{y_0}$ be the pyramid $(y_0-C) \cap \beta^{-1}([0,\infty))$.

\begin{step}[Expressing the volume of the pyramid $Y_{y_0}$, first way]
  We claim that
  \begin{equation}
    \label{eq:pyramid-comparison}
    \vol(Y_{y_0}) = \frac{1}{n} \bar\beta(y_0) \occ(\Lambda_0, P) \vol(\RR^n/\Lambda).
  \end{equation}
  One can see this by comparing the pyramid $Y_{y_0}$ with a pyramid $Y$ whose
  base is a parallelepiped spanned by a basis of the lattice $\Lambda_0$ in the hyperplane
  $\beta^{-1}(0)$ and whose tip is some $v \in \Lambda$ with
  $\bar\beta(v) = 1$. Note that a parallelepiped $Z$ that contains the base of $Y$ as a face and $e_0$ as a vertex is a fundamental domain for $\Lambda$, since it is spanned by $n$ linearly independent elements of $\Lambda$ and it only contains elements of $\Lambda$ at its vertices, since $\bar\beta$ is a primitive integral class with $\bar\beta(e_0) = 1$. So $\vol(Z) = \vol(\RR^n/\Lambda)$.
  
  Recall that the volume formulas for a pyramid and a
  parallelepided are $\frac{1}{n}bh$ and $bh$ respectively, where $b$ is the
  $n-1$-dimensional area of the base and $h$ is the height. Since $Y$ and $Z$ have the same base and same height, but $Y$ is a pyramid and $Z$ is a parallelepiped, we have
  \begin{equation}\label{eq:pyramid-volume}
    \vol(Y) = \frac{1}{n}\vol(Z) =  \frac{1}{n}\vol(\RR^n/\Lambda).
  \end{equation}

  Both $Y_{y_0}$ and $Y$ are pyramids with a base on the hyperplane
  $\beta^{-1}(0)$. To compare the volumes, we need to compare their heights and
  the areas of their bases. The base of $Y_{y_0}$ is $P_0$ and the base of $Y$
  is a fundamental domain for $\Lambda_0$. Therefore the ratio of the areas of
  the bases is the occupancy coefficient
  $\frac{\vol(P_0)}{\vol(\beta^{-1}(0)/\Lambda_0)} = \occ(\Lambda_0, P_0) =
  \occ(\Lambda_0, P)$. The tips of $Y_{y_0}$ and $Y$ are at the level sets
  $\bar\beta^{-1}(\bar\beta(y_0))$ and $\bar\beta^{-1}(1)$ respectively,
  therefore the height of $Y_{y_0}$ is $\bar\beta(y_0)$ times the height of
  $Y$. The formula \Cref{eq:pyramid-comparison} follows from
  \Cref{eq:pyramid-volume} and the comparisons between the bases and the
  heights.
\end{step}

\begin{step}[Expressing the volume of the pyramid $Y_{y_0}$, second way]
  Let $Y_{e_0}$ be the pyramid $(e_0-C) \cap \beta^{-1}([0,\infty))$. Using the
  similarity between the polytopes $Y_{y_0}$ and $Y_{e_0}$, we have
  \begin{equation}\label{eq:similarity}
    \vol(Y_{y_0}) = \frac{\bar\beta(y_0)^n}{\bar\beta(e_0)^n}\vol(Y_{e_0}).
  \end{equation}
\end{step}

\begin{step}[Expressing $\bar\beta(y_0)$]
  Using the equality between the right hands sides of
  \Cref{eq:pyramid-comparison} and \Cref{eq:similarity} and solving for
  $\bar\beta(y_0)$, we obtain
  \begin{equation}\label{eq:barbetay_0}
    \bar\beta(y_0) = \sqrt[n-1]{\frac{\occ(\Lambda_0,P)\vol(\RR^n/\Lambda)}{n\vol(Y_{e_0})}}
    \bar\beta(e_0)^{1+\frac{1}{n-1}}.
  \end{equation}
\end{step}

\begin{step}[Conclusion]
  The pyramid $Y_{e_0}$ is isometric to the pyramid $C \cap \beta^{-1}([0,1])$
  The statement of the lemma now follows from \Cref{eq:barbetay_0}, the upper and
  lower bounds on the Frobenius number (\Cref{eq:frob-upper-bound} and
  \Cref{eq:frob-lower-bound}) and the fact that all the error terms
  ($\bar\beta(x_1)$, $\bar\beta(x_2)$ and $m\bar\beta(e_0)$) are constant
  multiples of $\bar\beta(e_0)$ where the constant depends only on
  $x_1, x_2, C, \Lambda$ and $e_0$.
\end{step}
\end{proof}

\begin{remark}\label{rem:beta-one-goes-to-infty}
  For a given $C$, $e_0$ and $\Lambda$ as in \Cref{lemma:main-frob-bound} and
  any $K>0$, there are only finitely many linear functions
  $\bar\beta: V \to \RR$ taking positive values on
  $C-\{0\}$, integer values on $\Lambda$ such that $\bar\beta(e_0) < K$.
  To see this, let $v_1, \ldots, v_n \in C$ be a basis for
  $\Lambda$ such that $e_0$ is in the interior of the cone generated by the
  $v_i$. That is, $e_0 = c_1v_1 + \ldots c_nv_n$ with $c_i > 0$ for
  $i=1,\ldots, n$. Then
  $\bar\beta(e_0) = c_1\bar\beta(v_1) + \cdots + c_n \bar\beta(v_n)$ where
  the $\bar\beta(v_i)$ are positive integers. From this, we see that there are only
  finitely many choices for the $\bar\beta(v_i)$ that make this sum less than $K$.
  Since the $\bar\beta(v_i)$ determine $\bar\beta$, we obtain that there are
  indeed finitely many possibilities for $\bar\beta$.
\end{remark}

From \Cref{lemma:main-frob-bound} and \Cref{rem:beta-one-goes-to-infty}, we
obtain the following.
\begin{corollary}\label{cor:frob-limit}
  Let $V, C, e_0, \Lambda, \Lambda'$ be as in \Cref{lemma:main-frob-bound}.
  Consider a sequence of pairwise distinct linear functions
  $\{\beta_k\}_{k\in\NN}$ and associated polytopes $\{P_k\}_{k\in\NN}$ as in
  \Cref{lemma:main-frob-bound}. Then
  \begin{displaymath}
    \lim_{k \to \infty}\frac{\Frob(\bar\beta_k|_{\Lambda'})}{\bar\beta_k(e_0)^{1+\frac{1}{n-1}}} -
    \sqrt[n-1]{\frac{\occ(\Lambda \cap \beta_k^{-1}(0),P_k)\cdot \vol(V/\Lambda)}{n\vol(C \cap \beta_k^{-1}([0,1]))}} = 0
  \end{displaymath}
\end{corollary}

\subsection{Cones with a tetrahedron base}
\label{sec:dim-two}

In the following lemmas, $\pi_i \co \RR^n \to \RR$ denotes the
projection to the $i$th coordinate.

\newcommand\balpha{{\boldsymbol\alpha}}

\begin{lemma}\label{lemma:triangle-area}
  Let $n\ge 1$ be an integer and let
  $\balpha = (\alpha_1,\ldots, \alpha_n)$ be such that
  $\sum_{i=1}^n \alpha_i = 1$ and $\alpha_i > 0$ for each $i$. Let
  $\beta_\balpha = \sum_{i=1}^n \alpha_i \pi_i$ and denote by $T_\balpha$ the
  tetrahedron $\RR_{\ge 0}^n \cap \beta_\balpha^{-1}([0,1])$. Then
  \begin{displaymath}
    \vol(T_\balpha) = \frac{1}{n!}\prod_{i=1}^n \frac{1}{\alpha_i}.
  \end{displaymath}
\end{lemma}
\begin{proof}
  One vertex of $T_\balpha$ is the origin, the other $n$ vertices are the
  intersections of the hyperplane $\beta_\balpha^{-1}(1)$ with the coordinate
  axes. For example, the intersection with the first axis is the point
  $(x_1,0,\ldots,0)$ that satisfies
    \begin{displaymath}
    1 = \beta_\balpha(x_1,0,\ldots,0) = \alpha_1x_1,
  \end{displaymath}
  which yields $x_1 = \frac{1}{\alpha_1}$. Similarly, we obtain that the only
  nonzero coordinates of the other intersection points are
  $\frac{1}{\alpha_i}$. The parallelepided spanned by the intersection points
  has volume $\prod_{i=1}^n \frac{1}{\alpha_i}$ and the tetrahedron spanned by
  them has volume $\frac{1}{n!}$ times that.
\end{proof}

\subsection{Projective convergence of lattices}
\label{sec:projective-convergence}

The occupancy coefficient term in \Cref{cor:frob-limit} is not very
well-behaved, since the lattices
$\Lambda \cap \beta_k^{-1}(0)$ may vary wildly even when the linear functions $\beta_k$ converge. It will be useful to single out
subsequences where the lattices are stable in a sense. This section introduces
some lemmas and terminology for this.

\begin{lemma}\label{lemma:concrete-hyperplane-lattices}
  If $a_1, \ldots, a_n \in \ZZ$ and $v_1, \ldots, v_n \in \RR^{n+1}$ are the
  columns of the $(n+1) \times n$ matrix
  \begin{equation}\label{eq:lattice-columns}
    \begin{pmatrix}
      a_1 & 1 & 0 & 0 &  \cdots & 0 \\
      0 & a_2 & 1 & 0 & \cdots & 0 \\
      0 & 0 & a_3 & 1 & \ddots & 0\\
      \vdots & \vdots & \ddots & \ddots & \ddots & \\
      0 & 0 & 0 & \ddots & a_{n-1} & 1\\
      0 & 0 & 0 & \cdots & 0 & a_n \\
      1 & 0 & 0 & \cdots & 0 & 0 \\
    \end{pmatrix},
  \end{equation}
  then $\langle v_1, \ldots, v_n \rangle_\ZZ = \langle v_1, \ldots, v_n
  \rangle_\RR \cap \ZZ^{n+1}$.
\end{lemma}
\begin{proof}
  The abelian group $\langle v_1, \ldots, v_n \rangle_\ZZ$ is a finite index
  subgroup of $\langle v_1, \ldots, v_n \rangle_\RR \cap \ZZ^{n+1}$. To show
  that they are equal, it suffices to show that there is some $v_{n+1}\in
  \ZZ^{n+1}$ such that $v_1, \ldots, v_{n+1}$ form a basis for $\ZZ^{n+1}$. The
  vector $v_{n+1} = (0, \ldots, 0, 1, 0)^T$ has this property, since by adding
  this vector as the last column of the matrix above, we obtain a matrix whose
  determinant is $\pm 1$.
\end{proof}

If $\Lambda$ and $\{\Lambda_k\}_{k\in \NN}$ are discrete subgroups of rank
$r$ in a vector space $V \cong \RR^n$, then we say that
$\Lambda_k \to \Lambda$ \emph{projectively} if there is a basis $v^1, \ldots, v^r$ of
$\Lambda$, a positive constant $c_k$ and a basis $v_k^1, \ldots, v_k^r$ of
$\Lambda_k$ for every $k \in \NN$ such that $\lim_{k\to \infty} c_k v_k^i =
v^i$ for every $i = 1, \ldots, r$.

\begin{lemma}\label{lemma:projective-lattice-convergence}
  Let $\Lambda$ be a lattice in an $n$-dimensional vector space $V$. Let
  $\Sigma$ be a hyperplane in $V$ and let $\Lambda_0$ be any lattice in $\Sigma$. Then there exists a
  sequence $\{\Sigma_k\}_{k\in \NN}$ of hyperplanes in $V$ such that $\Sigma_k
  \cap \Lambda$ is a lattice in $\Sigma_k$ for all $k$ and $\Sigma_k \cap
  \Lambda \to \Lambda_0$ projectively.
\end{lemma}
\begin{proof}
  First we prove the statement in a special case and then we use this special case to prove the general case.
  
  In the special case, we assume that $\Sigma \cap \Lambda$ is a lattice in $\Sigma$. In this case, we can
  choose an isomorphism $\iota: V \to \RR^n$ that identifies $\Lambda$ with
  $\ZZ^n$ and the hyperplane $\Sigma$ with the orthogonal complement of
  $(0,\ldots,0,1)^T$. In this special case, we further assume that $\Lambda_0$ is \emph{rational}, that
  is, there is a positive constant $c$ such that $c\Lambda_0 \subset \Lambda$. Then
  we may choose the isomorphism $\iota$ such that the columns of the
  $n \times (n-1)$ matrix
  \begin{displaymath}
    \begin{pmatrix}
      b_1 & 0 &  \cdots & 0 \\
      0 & b_2 &  \ddots & 0 \\
      \vdots & \ddots & \ddots& \vdots \\
      0 & 0 & \cdots & b_{n-1} \\
      0 & 0 & \cdots & 0 \\
    \end{pmatrix}
  \end{displaymath}
  form a basis for $\iota(c\Lambda_0)$ for some $b_1, \ldots, b_{n-1} \in \ZZ$. Applying
  \Cref{lemma:concrete-hyperplane-lattices} with $a_i^{(k)} = kb_i$ for
  $i=1,\ldots,n-1$ and all $k \in \NN$ yields a sequence of hyperplanes
  (spanned by the columns of the matrix \Cref{eq:lattice-columns}) whose
  pullbacks by $\iota$ satisfy the required properties.

  To prove the general case, we take a sequence of hyperplanes $\Sigma_m$ converging to $\Sigma_0$ such that $\Sigma_m \cap \Lambda$ is a lattice in $\Sigma_m$ for all $m$ and for each $m$ we take a rational lattice
  $\Lambda_0^{(m)}$ in $\Sigma_m$ so that the lattices $\Lambda_0^{(m)}$ converge to
  $\Lambda_0$ projectively. As we have already shown, we can construct
  sequences of lattices of the required properties converging projectively to
  each $\Lambda_0^{(m)}$. From this, the statement of the lemma follows also
  for $\Lambda_0$.
\end{proof}

\subsection{Projection of cones}

The following lemma will be used for projections of the cone dual to the cone
$\RR_+\calF$. Such projections naturally correspond to slices of the fibered face
$\calF$.

\begin{lemma}\label{lemma:image-of-integral-cone}
  Let $C = \langle \calB \rangle_{\RR_{\ge 0}} \subset \RR^n$ be the cone
  generated by some subset $\calB \subset \RR^n$. Suppose $C$ has nonempty
  interior and let $p\co \RR^n \to \RR^k$ be a linear map such that $p(\ZZ^n)$
  is a lattice in $\RR^k$. Then there exists $x \in \RR^k$ such that
  $p(\ZZ^n) \cap (x + p(C)) \subset p(\ZZ^n \cap C)$.
\end{lemma}
\begin{proof}
  Since $p(\ZZ^n)$ is a lattice in $\RR^k$, the subset
  $A = p^{-1}(0) \cap \ZZ^n$ is a lattice in the $n-k$-dimensional kernel
  $p^{-1}(0)$. Let $K>0$ be large enough so that $B(y,K)$, the ball of radius
  $K$ centered at $y$, contains some $a \in A$ for all $y \in p^{-1}(0)$.
  Since the cone $C$ has nonempty interior, there is some $c_0 \in C$ such that
  $B(c_0,K) \subset C$. But then $B(c,K) \subset C$ for all $c \in c_0+C$.

  We claim that any point $b \in p(\ZZ^n)$ that lies in
  $p(c_0+C) = p(c_0) + p(C)$ is the image of some point of $\ZZ^n$ that lies in
  $C$. This will prove the lemma with $x = p(c_0)$. To see this, let
  $c \in c_0+C$ such that $p(c) = b$. We know that $B(c,K) \subset C$, and the
  ball $B(c,K)$ is centered around a point of the translate $p^{-1}(0) + c$ of
  the subspace $p^{-1}(0)$. Moreover, this translated subspace $p^{-1}(0)+c$, being
  a preimage of $b$, contains some point in $\ZZ^n$ hence contains a translate
  of $A$. As a consequence, $(p^{-1}(0)+c) \cap B(c,K) \cap \ZZ^n$ is nonempty
  and any element of it is a point of $\ZZ^n$ that lies in $C$ and maps to $b$.
\end{proof}

\section{Proof of the main theorem}
\label{sec:proof}

We are now ready to prove our main theorem.

In the statement of the theorem, we use the following notions of duality. Let
$V$ be a vector space, let $C = \langle \calB \rangle_{\RR_{\ge 0}}$ be cone
with nonempty interior, generated by some subset $\calB \subset V$, let
$\Lambda \subset V$ be a lattice. Then the dual of the triple $(V, C, \Lambda)$
is the triple $(V^*, C^*, \Lambda^*)$ where $V^*$ is the dual vector space of
$V$, the cone $C^* \subset V^*$ is the set linear functions $V \to \RR$ that
take nonnegative values on $C$ and $\Lambda^* \subset V^*$ is the lattice
consisting of linear functions that take integer values on $\Lambda$. Note that
if $\phi: V \to \RR$ is a linear function that takes positive values on $C$,
then $\phi \in \inter(C^*)$.

\Cref{thm:closure} is a direct corollary of the following theorem which, in
addition to \Cref{thm:closure}, describes the bounding function $g$.

\begin{theorem}\label{thm:closure-concrete}
  Let $M$ be a connected 3-manifold that admits a complete finite-volume
  hyperbolic metric. Let $\calF$ be a fully-punctured fibered face of the unit
  ball of the Thurston norm on $H^1(M;\RR)$. Let
  $1 \le d \le \dim(H^1(M;\RR))-1$, let $\Omega$ be a rational $d$-dimensional
  slice of $\calF$ cut out by the $d+1$-dimensional subspace $\Sigma$, let $C$
  be the cone $\langle \Omega \rangle_{\RR_{\ge0}}$ in $\Sigma$ and consider
  the lattice $\Lambda = \Sigma \cap H^1(M;\ZZ)$ in $\Sigma$. Consider the dual
  triple $(\Sigma^*, C^*, \Lambda^*)$ of the triple $(\Sigma, C, \Lambda)$.

  Let $\graph(\mu_d|_\Omega) \subset \Omega \times \RR$ be the graph of the
  normalized asymptotic translation length function $\mu_d$, restricted to
  $\Omega$. Let $g \co \inter(\Omega) \to \RR_+$ be the function defined by the
  formula
  \begin{equation}
    \label{eq:bounding-function}
    g(\phi) = \sqrt[d]{\frac{(d+1)\vol_{\Lambda^*}(C^* \cap \beta_{\phi}^{-1}([0,1]))}{\min\occ(C^* \cap \beta_{\phi}^{-1}(1))}}
  \end{equation}
  where $\vol_{\Lambda^*}$ is the translation-invariant volume-form on
  $\Sigma^*$ with respect to which $\Lambda^*$ has covolume 1 and $\beta_{\phi}$ denotes
  the linear function $\Sigma^* \to \RR$ corresponding to the element $\phi \in
  \Sigma$ in the dual space of $\Sigma^*$.

  Then the set of accumulation points of the graph $\graph(\mu_d|_\Omega)$
  is
  \begin{displaymath}
    \{(\omega,g(\omega)) : \omega \in \inter(\Omega)\}
  \end{displaymath}
  if $d=1$ and
  \begin{displaymath}
    \{(\omega,r) : \omega \in \inter(\Omega), \,0 \le r \le g(\omega)\} \cup
    (\bdy\Omega \times [0,\infty))
  \end{displaymath}
  if $d\ge 2$.
  Moreover, the function $g$ is continuous and $g(\phi) \to \infty$ as
  $\phi \to \bdy\Omega$.
\end{theorem}

\begin{proof}
  Once again, we break the proof into several steps.

  \setcounter{step}{0}
  \begin{step}[Asymptotic behavior of the weighted graphs $W(\phi)$]
    Every element of $H_1(M;\RR)$ defines a linear function $H^1(M;\RR) \to
    \RR$ and this linear function restricts to a linear function $\Sigma \to
    \RR$. This way we obtain a natural linear map $p: H_1(M;\RR) \to \Sigma^*$. It is easy
    to see that $p(G) = \Lambda^*$.

    By \Cref{prop:cone-description}, the elements of the cone $C$ take
    nonpositive values on the cone $\langle \calB\rangle_{\RR_{\ge 0}}$,
    therefore $p(\langle \calB\rangle_{\RR_{\ge 0}}) \subset -C^*$. Conversely,
    every element of $-C^*$ is a linear function $\Sigma \to \RR$ that takes
    nonpositive values on $C$ and every such linear function is a restriction
    of a linear function $H^1(M;\RR) \to \RR$ that takes nonpositive values on
    the cone $\RR^+\calF$ over the fibered face. By the other direction of
    \Cref{prop:cone-description}, this last linear function corresponds to an
    element of $\langle \calB\rangle_{\RR_{\ge 0}}$. Hence we have
    $p(\langle \calB\rangle_{\RR_{\ge 0}}) = -C^*$.

    Fix some $e, e' \in E$. By \Cref{prop:automaton-and-intersecting-edges},
    there are $g_1, g_2 \in G$ such that
    \begin{displaymath}
      G \cap g_1 \langle \calB\rangle_{\RR_{\ge 0}} \subset \stash(e,e') \subset G \cap g_2 \langle \calB\rangle_{\RR_{\ge 0}}.
    \end{displaymath}

    We claim that it follows that there are some $v_1, v_2 \in \Sigma^*$ such that
    \begin{displaymath}
      \Lambda^* \cap (v_1 + C^*) \subset -p(\stash(e,e')) \subset \Lambda^* \cap (v_2 + C^*),
    \end{displaymath}
    using additive notation in the vector space $\Sigma^*$. The existence of
    $v_1$ follows from \Cref{lemma:image-of-integral-cone}. The existence of
    $v_2$ is simply a consequence of the identity
    $p(A \cap B) \subset p(A) \cap p(B)$.

    Let $e_0: \Sigma \to \RR$ be the linear function that takes the value 1 on
    $\Omega$. Note that $e_0$ is in the interior of $C^*$.

    We now wish to apply \Cref{cor:frob-limit} with
    $V = \Sigma^*, C = C^*, e_0, \Lambda =\Lambda^*, \Lambda'=-p(\stash(e,e'))$ to conclude that
    \begin{equation}
      \label{eq:first-form}
      \lim_{\phi \to \phi_0}\frac{\Frob(\bar\beta_{\phi}|_{\Lambda'})}{\bar\beta_{\phi}(e_0)^{1+\frac{1}{d}}} -
      \sqrt[d]{\frac{\occ(\Lambda^* \cap \beta_{\phi}^{-1}(0),P_\phi)\cdot \vol(\Sigma^*/\Lambda^*)}{(d+1)\vol(C^* \cap \beta_{\phi}^{-1}([0,1]))}} = 0,
    \end{equation}
    where $P_\phi = \beta_{\phi}^{-1}(0) \cap (y-C^*)$ for some
    $y \in \inter(C^*)$ and $\vol(\cdot)$ is any translation-invariant
    volume-form on $\Sigma^*$. The cohomology classes $\phi$ are
    rational points of the interior of $\Omega$ and $\phi_0$ is an arbitrary
    element of $\Omega$.

    It is straightforward to check that the hypotheses of \Cref{cor:frob-limit}
    are satisfied. For example, the classes $\beta_\phi$ take positive values
    on $C^*$, since $\phi$ is assumed to be in the interior of $\Omega$ and
    hence in the interior of $C$. The equality
    $\beta_\phi(e_0) = 1$ holds because $e_0$ takes the value 1 on $\Omega$ and
    $\phi$ is in $\Omega$. The rest of the hypotheses are also satisfied,
    therefore \Cref{cor:frob-limit} applies and \Cref{eq:first-form} holds.

    Denote by $w_\phi(ee')$ the weight of the edge $ee'$ in the weighted graph
    $W(\phi)$. Using the definition \Cref{eq:weight}, we have
    \begin{displaymath}
      w_\phi(ee') = \Frob(\bar\phi|_{-\stash(e,e')}) =  \Frob(\bar\beta_\phi|_{\Lambda'}) .
    \end{displaymath}
    To see that the second equality holds, first note that
    $\beta_\phi(p(x)) = \phi(x)$ for all $x\in H_1(M;\RR)$. So
    $\phi(G) \subset \QQ$ is the same discrete subgroup as
    $\beta_\phi(\Lambda^*) = \beta_\phi(p(G)) \subset \QQ$, so
    $\bar\phi = c \phi$ and $\bar\beta_{\phi} = c\beta_\phi$ hold with the same
    constant $c>0$. Therefore $\bar\phi(x) = \bar\beta_\phi(p(x))$ for all
    $x\in H_1(M;\RR)$.

    Finally, we have $\beta_\phi(e_0) =1 = \|\phi\|$ whenever
    $\phi \in \Omega$. After multiplying by $c$, we obtain
    $\bar\beta_\phi(e_0) = \|\bar\phi\|$.
    Applying these
    substitutions to \Cref{eq:first-form} yield
    \begin{equation}\label{eq:weight-limit}
      \lim_{\phi \to \phi_0}\frac{w_\phi(ee')}{\left\|\bar\phi\right\|^{1+\frac{1}{d}}} -
      \sqrt[d]{\frac{\occ(\Lambda^* \cap \beta_{\phi}^{-1}(0),P_\phi)\cdot
          \vol(\Sigma^*/\Lambda^*)}{(d+1)\vol(C^* \cap \beta_{\phi}^{-1}([0,1]))}} = 0.
    \end{equation}
  \end{step}

  \begin{step}[Accumulation points in $\inter(\Omega) \times \RR$]
    Assume that $\phi_0$ is in the interior of $\Omega$. Let us apply
    \Cref{lemma:projective-lattice-convergence} with $V=\Sigma^*$, $\Lambda = \Lambda^*$, $\Sigma = \beta_{\phi_0}^{-1}(0)$ and for a lattice $\Lambda_0$ in $\Sigma$ such that the $\occ(\Lambda_0, P_{\phi_0}) = \alpha$ for some $\alpha \in \occs(P_{\phi_0})$. Such a lattice $\Lambda_0$ exists by the definition of the set
    $\occs(P_{\phi_0})$ in \Cref{eq:occs-def}. \Cref{lemma:projective-lattice-convergence} guarantees that there exists a sequence
    $\phi_k \to \phi_0$ such that
    \begin{displaymath}
      \lim_{k\to \infty} \occ(\Lambda^* \cap \beta_{\phi_k}^{-1}(0), P_{\phi_k}) = \occ(\Lambda_0, P_{\phi_0}) = \alpha
    \end{displaymath}
    and therefore
    \begin{equation}\label{eq:edge-weight-limit}
      \lim_{k\to \infty}\frac{w_{\phi_k}(ee')}{\left\|\bar\phi_k\right\|^{1+\frac{1}{d}}} = \sqrt[d]{\frac{\alpha \vol(\Sigma^*/\Lambda^*)}{(d+1)\vol(C^* \cap
          \beta_{\phi_0}^{-1}([0,1]))}}
    \end{equation}
    for each edge $ee'$. Since $\phi_0$ is in the interior of $\Omega$, the set
    $C^* \cap \beta_{\phi_0}^{-1}([0,1])$ is a pyramid of finite volume.

    Now recall from \Cref{eq:normalized-atl} that
    $\mu_d(\phi) = \lVert \bar\phi \rVert^{1+\frac{1}{d}}
    \ell_\calA(\bar\phi)$. Moreover, by \Cref{prop:atl-via-weighted-graphs},
    $\ell_\calA(\bar\phi)$ equals the reciprocal of the maximal average cycle
    weight in the graph $W(\phi)$. According to \Cref{eq:edge-weight-limit},
    the weight of each edge has the same asymptotics for the sequence
    $\phi_k$, therefore the average weight of every cycle also has the same
    asymptotics. So we can replace $w_{\phi_k}(ee')$ by
    $\frac{1}{\ell_\calA(\bar\phi_k)}$ in \Cref{eq:edge-weight-limit} and the
    limit still holds. Taking the reciprocal of both sides, we obtain
    \begin{align*}
      \lim_{k\to \infty}\mu_d(\phi_k) &= \sqrt[d]{\frac{(d+1)\vol(C^* \cap
                                        \beta_{\phi_0}^{-1}([0,1]))}{\alpha \vol(\Sigma^*/\Lambda^*)}} \\&= \sqrt[d]{\frac{(d+1)\vol_{\Lambda^*}(C^* \cap
      \beta_{\phi_0}^{-1}([0,1]))}{\alpha}}.
    \end{align*}
    Such a sequence $\phi_k \to \phi$ exists for every
    $\alpha \in \occs(P_{\phi_0})$. Using \Cref{lemma:occs} and the fact that
    the polytope $P_{\phi_0}$ in the hyperplane $\beta_{\phi_0}^{-1}(0)$ and
    the polytope $C^* \cap \beta_{\phi_0}^{-1}(1)$ in the hyperplane
    $\beta_{\phi_0}^{-1}(1)$ are homothetic, we obtain that the accumulation
    points of $\graph(\mu_d|_\Omega)$ in $\inter(\Omega) \times \RR$ are as
    specified in the theorem.
  \end{step}

  \begin{step}[Accumulation points in $\bdy\Omega \times \RR$] If
    $\phi \to \phi_0$, then the pyramid $C^* \cap \beta_\phi^{-1}([0,1])$
    converges to the set $C^* \cap \beta_{\phi_0}^{-1}([0,1])$. If
    $\phi_0 \in \bdy \Omega$ and hence $\phi_0 \in \bdy C$, then this limit set
    is unbounded and has infinite volume. By \Cref{lemma:occ-dim-one}, the
    occupancy coefficient is 1 if $d=1$, therefore the expression under the
    root in \Cref{eq:weight-limit} goes to 0. Hence
    $\lim_{\phi \to
        \phi_0}\frac{w_\phi(ee')}{\left\|\bar\phi\right\|^{1+\frac{1}{d}}} = 0$
    and $\lim_{\phi\to\phi_0} \mu_d(\phi) = \infty$ when
    $\phi_0\in \bdy \Omega$ and $d=1$. Therefore $\graph(\mu_d|_{\Omega})$ does
    not have any accumulation point in $\bdy\Omega \times \RR$ when $d=1$.

    Similarly, the limit of the pyramids $C^* \cap \beta_\phi^{-1}([0,1])$ in
    the definition of $g$ is an unbounded set when $\phi \to \bdy \Omega$. This
    shows that $g(\phi) \to \infty$ whenever $\phi \to \bdy \Omega$. It is
    clear that $g$ is continuous.

    It is now automatic that the set of accumulation points of
    $\graph(\mu_d|_{\Omega})$ in $\bdy\Omega\times \RR$ is
    $\bdy\Omega \times [0,\infty)$ when $d \ge 2$, since $g$ goes to infinity
    at $\bdy\Omega$ and the set of accumulation points is closed.
  \end{step}
\end{proof}

Next, we prove \Cref{thm:simplex}.

\begin{proof}[Proof of \Cref{thm:simplex}]
  Choose $\omega_1$, \ldots, $\omega_{d+1}$, the vertices of the simplex
  $\Omega$, as the basis for $\Sigma$. This choice of basis naturally defines
  coordinates on $\Sigma$ and the dual space $\Sigma^*$. With these
  coordinates, we have $C \cong \RR_{\ge0}^{d+1} \cong C^*$.

  By \Cref{thm:closure-concrete}, the function $g$ takes the form
  \begin{displaymath}
    g(\phi) = \sqrt[d]{\frac{(d+1)\vol_{\Lambda^*}(\RR_{\ge0}^{d+1} \cap \beta_\phi^{-1}([0,1]))}{\min\occ(\RR_{\ge0}^{d+1} \cap \beta_\phi^{-1}(1))}}.
  \end{displaymath}
  Using the definition of $g^*$ in the theorem and denoting $\beta_{\phi}$
  for $\phi = \sum_{i=1}^{d+1} \alpha_i \omega_i$ by $\beta_{\balpha}$ where $\balpha = (\alpha_1, \ldots, \alpha_{d+1})$, we can
  rewrite this as
  \begin{equation}\label{eq:g-star}
    g^*(\balpha) = \sqrt[d]{\frac{(d+1)\vol_{\Lambda^*}(\RR_{\ge0}^{d+1} \cap \beta_\balpha^{-1}([0,1]))}{\min\occ(\RR_{\ge0}^{d+1} \cap \beta_\balpha^{-1}(1))}}.
  \end{equation}

  Note that $\RR_{\ge0}^{d+1} \cap \beta_\balpha^{-1}(1)$ is a $d$-dimensional simplex for all
  $\balpha$. All $d$-dimensional simplices have the same minimal occupancy
  coefficient, since they differ only by a linear transformation. Therefore the
  denominator is a constant $O_d$, depending only on $d$. (It is
  straightforward to check that this definition of $O_d$ is equivalent to the
  definition provided in the introduction after \Cref{thm:simplex}.)

The volume in the numerator can be written as
\begin{equation}\label{eq:fraction}
  \vol_{\Lambda^*}(\RR_{\ge0}^{d+1} \cap \beta_\balpha^{-1}([0,1])) = \frac{\vol(\RR_{\ge0}^{d+1} \cap \beta_\balpha^{-1}([0,1]))}{\vol(\RR^{d+1}/\Lambda^*)}
\end{equation}
where $\vol(\cdot)$ denotes the standard volume form on $\RR^{d+1}$. Note that
$\beta_{\balpha} = \sum_{i=1}^{d+1} \alpha_{i} \pi_i$ where $\pi_i: \RR^{d+1}
\to \RR$ is the projection to the $i$th coordinate. So we can apply
\Cref{lemma:triangle-area} to obtain that
\begin{equation}\label{eq:tetrahedron-volume}
  \vol(\RR_{\ge0}^{d+1} \cap \beta_{\balpha}^{-1}([0,1])) = \frac{1}{(d+1)!}\prod_{i=1}^{d+1} \frac{1}{\alpha_i}.
\end{equation}

Finally, the covolume of $\Lambda^*$ equals the reciprocal of the covolume of
$\Lambda$. By our choice of basis, the volume form on $\Sigma \cong \RR^{d+1}$
is the one with respect to which the lattice $\Gamma = \langle \omega_1, \ldots,
\omega_{d+1} \rangle_{\ZZ}$ has covolume 1. So
\begin{equation}\label{eq:flip-volume}
  \vol(\RR^{d+1}/\Lambda^*) = \frac{1}{\vol(\RR^{d+1}/\Lambda)} =
  \frac{1}{\vol_{\Gamma}(\Sigma/\Lambda)} =
\vol_\Lambda(\Sigma/\Gamma).
\end{equation}

  Putting together \Cref{eq:g-star}, \Cref{eq:fraction},
  \Cref{eq:tetrahedron-volume} and
  \Cref{eq:flip-volume}, we obtain that
  \begin{displaymath}
        g^*(\balpha) = \sqrt[d]{\frac{\frac{1}{d!}\prod_{i=1}^{d+1}
            \frac{1}{\alpha_i}}{O_d \vol_\Lambda(\Sigma/\Gamma)}}.
  \end{displaymath}
  This is what we wanted to prove. The fact that $O_d = 1$ follows from \Cref{lemma:occ-dim-one}.
\end{proof}

Using \Cref{thm:closure-concrete}, we can also prove \Cref{thm:shape}.

\begin{proof}[Proof of \Cref{thm:shape}]
  The isomorphism $i: \Sigma_1 \to \Sigma_2$ induces a dual isomorphism
  $i^*: \Sigma_2^* \to \Sigma_1^*$ of the dual spaces. By indexing the objects
  in the statement of \Cref{thm:closure-concrete} by 1 and 2, corresponding to
  the manifolds $M_1$ and $M_2$, respectively, the isomorphism $i^*$ identifies
  $C_2^*$ with $C_1^*$. From
  \begin{displaymath}
    \theta = \frac{\vol(\Sigma_2/\Lambda_2)}{\vol(\Sigma_2/i(\Lambda_1))},
  \end{displaymath}
  we obtain that
  \begin{displaymath}
    \theta = \frac{\vol(\Sigma_1/\Lambda^*_1)}{\vol(\Sigma_1/i^*(\Lambda_2))}
  \end{displaymath}
  and therefore $\vol_{i^*(\Lambda_2)} = \theta \vol_{\Lambda^*_1}$.
  The functions $\beta_{\phi}$ are also
  identified in the sense that for $\phi_1 \in \Omega_1$, we have
  $\beta_{i(\phi_1)}(x) = \beta_{\phi_1}(i^*(x))$ for all $x \in \Sigma_2$.
  So
  \begin{align*}
    g_2(i(\phi_1))
    &= \sqrt[d]{\frac{(d+1)\vol_{\Lambda^*_2}(C^*_2
      \cap \beta_{i(\phi_1)}^{-1}([0,1]))}{\min\occ(C^*_2 \cap
      \beta_{i(\phi_1)}^{-1}(1))}} \\
    &= \sqrt[d]{\frac{(d+1)\vol_{i^*(\Lambda^*_2)}(C^*_1 \cap
      \beta_{\phi_1}^{-1}([0,1]))}{\min\occ(C^*_1 \cap
      \beta_{\phi_1}^{-1}(1))}} \\
    &= \sqrt[d]{\frac{(d+1) \theta \vol_{\Lambda^*_1}(C^*_1 \cap
      \beta_{\phi_1}^{-1}([0,1]))}{\min\occ(C^*_1 \cap
      \beta_{\phi_1}^{-1}(1))}} = \theta^{\frac1d} g_1(\phi_1)
  \end{align*}
  by \Cref{thm:closure-concrete}.
\end{proof}

Finally, we show that the bounding function $g$ in \Cref{thm:simplex} is convex.

\begin{lemma}\label{lemma:convex-function}
  For any $\delta > 0$ and integer $n\ge 1$, the function 
  \begin{displaymath}
     f(\alpha_1, \ldots, \alpha_n) =\prod_{i=1}^{n}\alpha_i^{-\delta}
  \end{displaymath}
  is convex on its natural domain $\{(\alpha_1, \ldots, \alpha_{n}): \alpha_i>0 \mbox{ for } i=1, \ldots n\}$.
\end{lemma}
\begin{proof}
  The Hessian of $\log{f}$ is a diagonal matrix with diagonal entries $\frac{\delta}{\alpha_i^2}$. This matrix is positive definite, therefore $f$ is logarithmically convex. Every logarithmically convex function is also convex, since a composition of the convex function with the increasing convex function $e^x$ is also convex. Hence $f$ is indeed convex.
\end{proof}

\section{An example}
\label{sec:example}

In this section, we consider the simplest pseudo-Anosov braid on three strands,
describe the veering triangulation of its mapping torus and compute the
asymptotic translation length in the arc complex for infinitely many fibrations
of this 3-manifold. The purpose of this computation is two-fold: to illustrate the
methods of \Cref{sec:atl-via-graphs} on a concrete example and to show that it
seems very difficult to find an explicit formula for the normalized asymptotic
translation length functions $\mu_d$ defined in \Cref{eq:normalized-atl}.

A good reference for the pseudo-Anosov theory appearing in this section
(invariant train tracks, measured foliations and translation surfaces) is
Chapters 14 and 15 in \cite{FarbMargalit12}.

\begin{figure}[ht]
  \centering
  \begin{tikzpicture}[scale=1]
    \def\rad{0.02}
    \draw (0,0) circle (2);
    \node[below] at (0,2) {$P$};
    \foreach \i in {1,2,3}
    \filldraw (\i-2,0) circle (\rad) node[left] {$P_\i$};
    \foreach \j in {1,2}{
      \pgfmathsetmacro\i{\j-1.5}
      \draw[->] (\i-0.4,0.2) .. controls (\i-0.2,0.6) and (\i+0.2,0.6) ..
      node[above] {$\sigma_{\j}$}
      (\i+0.4,0.2);
      \draw[->] (\i+0.4,-0.2) .. controls (\i+0.2,-0.6) and (\i-0.2,-0.6) ..
      (\i-0.4,-0.2);
      }
  \end{tikzpicture}
  \caption{The half-twists $\sigma_1$ and $\sigma_2$.}
  \label{fig:sigmas}
\end{figure}
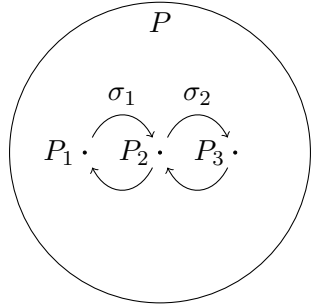

Let $S$ be the sphere punctured at four points $P_1$, $P_2$, $P_3$ and
$P$. Let $\sigma_1$ and $\sigma_2$ be the half-twists illustrated
on \Cref{fig:sigmas}.

\subsection*{Invariant train tracks}
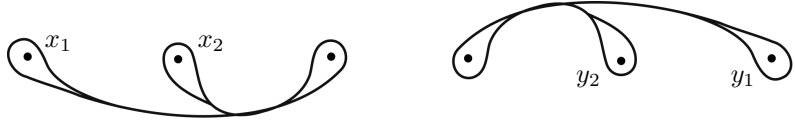
\begin{figure}[htb]
\labellist
\small\hair 2pt
 \pinlabel {$x_1$} [ ] at 32 47
 \pinlabel {$x_2$} [ ] at 128 47
 \pinlabel {$y_2$} [ ] at 365 25
 \pinlabel {$y_1$} [ ] at 462 25
\endlabellist
\centering
\includegraphics[scale=0.6]{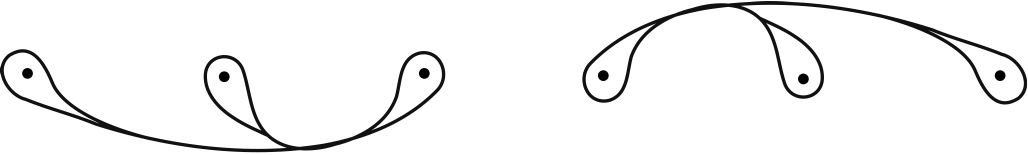}
\caption{The invariant train tracks $\tau$ and $\tau^{-1}$.}
\label{fig:invariant_train_tracks}
\end{figure}

The train track $\tau$ on the left of
\Cref{fig:invariant_train_tracks} is invariant under
$f = \sigma_1 \sigma_2^{-1}$ (read from left to right) and the train
track $\tau^{-1}$ on the right is invariant under $f^{-1}$. On each
train track, measures are parametrized by measures on two of the
branches. The action of $f$ and $f^{-1}$ on the measures are
$(x_1,x_2) \mapsto (x_2+2x_1,x_1+x_2)$ and
$(y_1,y_2) \mapsto (y_2+2y_1,y_1+y_2)$. In other words, both maps are
described by the matrix $
\begin{pmatrix}
  2 & 1 \\
  1 & 1 \\
\end{pmatrix}
$ whose eigenvalues are $\varphi^2$ and $\varphi^{-2}$, where $\varphi$ is the
golden ratio, the largest root of $x^2-x-1$. The eigenvector corresponding to
$\varphi^2$ is $(\varphi,1)$. Therefore the unstable foliation $\calF^u$ is
represented by the measure $(x_1,x_2) = (\varphi,1)$ on $\tau$ and the stable
foliation $\calF^s$ is represented by $(y_1,y_2) = (1,\varphi^{-1})$ on
$\tau^{-1}$. (The invariant measured foliations are well-defined only up to
scaling. We choose the scaling in a way that will be convenient later on.)

\subsection*{The half-translation surface}
\begin{figure}[ht]
  \centering
  \begin{tikzpicture}[scale = 3]
    \def\rad{0.02}
    \foreach \i in {-1,0,1}
    \filldraw (\i,0) circle (\rad);

    \begin{scope}
    \draw (0,0) node[below] {$P_2$} -- node[below right]
    {$(\varphi^2,1)$} (0,1) node[above] {$P$};
    \draw (-1,0) node[left] {$P_1$} -- node[left] {$(\varphi,\varphi)$} (0,1);
    \draw (1,0) node[right] {$P_3$} -- node[right] {$(1,\varphi^2)$} (0,1);
    \draw (-1,0) -- node[above] {$(1,\varphi^2)$} (0,0) -- node[above]
    {$(\varphi,\varphi)$} (1,0) .. controls (0.3,-1) and (-0.3,-1) ..
    node[below] {$(\varphi^2,1)$} (-1,0);
    \end{scope}
  \end{tikzpicture}
  \caption{The measures of the edges with respect to $\calF^u$ (second
    coordinate) and $\calF^s$ (first coordinate).}
  \label{fig:triangulation}
\end{figure}
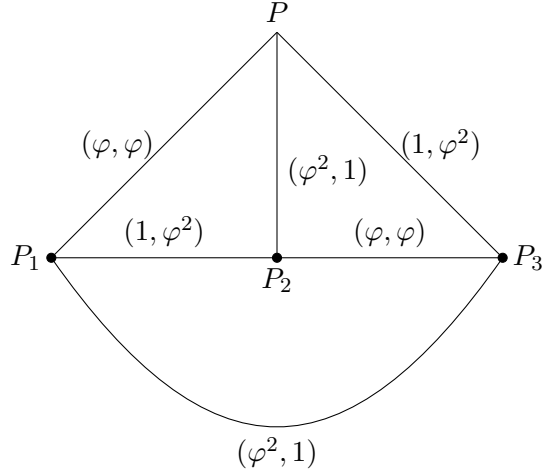

Our next goal is to draw a picture of the half-translation surface whose horizontal
foliation is $\calF^u$ and whose vertical foliation is $\calF^s$. For this,
consider the ideal triangulation of $S$ consisting of four triangles, shown on
\Cref{fig:triangulation}. The measures of the edges with respect to $\calF^u$
and $\calF^s$ can be obtained from the measured train tracks. These measures on
the edges provide the widths and heights of the edges in the half-translation
surface. Using these coordinates for the edges, we obtain the upper left
picture on \Cref{fig:translation-surface} showing the half-translation surface
defined by $\calF^u$ and $\calF^s$.

\newcommand\TriangleOne[9]{
    \def\x{1.618}
    \draw[thick,blue] (-\x,\x) -- node[below left] {$#1$} (0,0)  --
    node[below left] {$#2$} (\x,-\x);

    \draw[thick,violet] (-\x,\x) -- node[above] {$#3$}  ++(\x*\x,1)
    -- node[above] {$#4$} ++(\x*\x,1);

    \draw[thick,violet] (\x,-\x) -- node[right] {$#5$}
    ++(1,\x*\x) -- node[right] {$#6$}
    ++(1,\x*\x);

    \draw[thick,red] (0,0) -- node[below] {$#7$} (\x*\x,1);
    \draw[thick,red] (0,0) -- node[right] {$#8$} (1,\x*\x);
    \draw[thick,green] (\x*\x,1) -- node[above right] {$#9$} (1,\x*\x);

    \node[left] at (-\x,\x) {$P$};
    \node[below] at (0,0) {$P_1$};
    \node[above] at (\x*\x-\x,1+\x) {$P_2$};
    \node[right] at (1+\x,\x*\x-\x) {$P_3$};
  }


\newcommand\TriangleTwo[9]{
    \draw[thick,blue] (-\x*\x*\x,1/\x) -- node[below] {$#1$} (0,0) --
    node[below] {$#2$} (\x*\x*\x,-1/\x) ;
    \draw[thick,violet] (-\x*\x*\x,1/\x) -- node[above] {$#3$} ++(\x*\x,1)
     -- node[above] {$#4$} ++(\x*\x,1);
    \draw[thick,blue] (\x*\x*\x,-1/\x) -- node[right] {$#5$}
    ++(-\x,--\x) -- node[right] {$#6$} ++(-\x,--\x);
    \draw[thick,red] (0,0) -- node[right] {$#7$} (\x*\x,1);
    \draw[thick,green] (0,0) -- node[left] {$#8$} (-\x,--\x);
    \draw[thick,green] (\x*\x,1) --  node[above] {$#9$}(-\x,--\x);

    \node[left] at (-\x*\x*\x,1/\x) {$P$};
    \node[below] at (0,0) {$P_3$};
    \node[above] at (\x*\x-\x*\x*\x,\x) {$P_2$};
    \node[right] at (\x*\x*\x-\x,\x-1/\x) {$P_1$};
}

\newcommand\TriangleThree[9]{
    \draw[thick,violet] (-\x*\x*\x*\x,-1/\x/\x) -- node[below] {$#1$}
    (0,0)  -- node[below] {$#2$}
    (--\x*\x*\x*\x,--1/\x/\x);
    \draw[thick,violet] (-\x*\x*\x*\x,-1/\x/\x) --  node[above] {$#3$} ++(\x*\x,1) -- node[above] {$#4$} ++(\x*\x,1);
    \draw[thick,blue] (--\x*\x*\x*\x,--1/\x/\x) -- node[above] {$#5$}
    ++(-\x*\x*\x,--1/\x) -- node[above] {$#6$} ++(-\x*\x*\x,--1/\x);
    \draw[thick,red] (0,0) -- node[right] {$#7$}(\x*\x,1);
    \draw[thick,green] (0,0) -- node[left] {$#8$} (-\x*\x*\x,--1/\x);
    \draw[thick,red] (\x*\x,1) -- node[above] {$#9$} (-\x*\x*\x,--1/\x);

    \node[left] at (-\x*\x*\x*\x,-1/\x/\x) {$P$};
    \node[below] at (0,0) {$P_1$};
    \node[above] at (-\x*\x*\x*\x+\x*\x,1-1/\x/\x) {$P_2$};
    \node[above] at (\x*\x*\x*\x-\x*\x*\x,1/\x/\x+1/\x) {$P_3$};
}

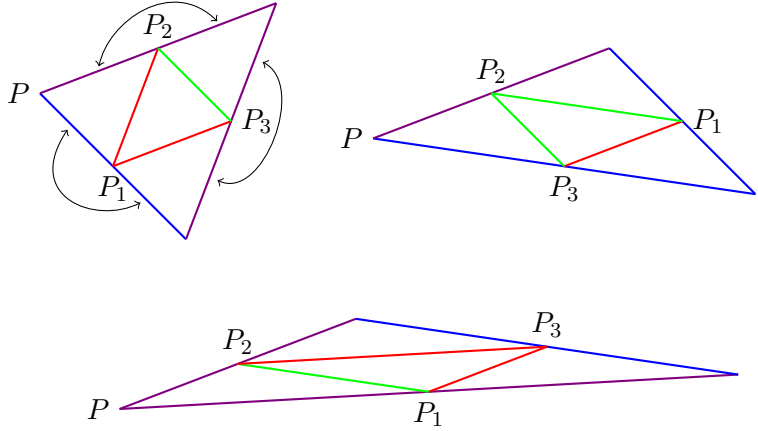
\begin{figure}[ht]
  \centering
  \def\x{1.618}
  \begin{tikzpicture}[scale=0.6]
    \TriangleOne{}{}{}{}{}{}{}{}{}

    \draw[<->] (-\x+0.5*\x*\x,\x+0.6) .. controls +(0.2,1) and +(-1,1)
    .. (-\x+1.5*\x*\x,\x+1.6);

    \draw[<->] (-0.5*\x-0.2,0.5*\x) .. controls +(-1,-1.5) and +(-1,-0.5)
    .. (0.5*\x-0.2,-0.5*\x);

    \draw[<->] (\x*\x - 0.5 +0.2,1 - 0.5*\x*\x) .. controls +(1,-0.5)
    and +(1,-0.5) .. (\x*\x + 0.5 +0.2,1 + 0.5*\x*\x);

    \begin{scope}[xshift=10cm]
      \TriangleTwo{}{}{}{}{}{}{}{}{}
    \end{scope}

    \begin{scope}[yshift=-5cm,xshift=7cm]
      \TriangleThree{}{}{}{}{}{}{}{}{}
    \end{scope}

  \end{tikzpicture}
  \caption{The half-translation surface defined by $\calF^u$ and $\calF^s$.
  Pairs of boundary edges are identified by 180 degree
  rotations.
  }
  \label{fig:translation-surface}
\end{figure}

\subsection*{The veering triangulation}
We can now use Gu\'eritaud's construction to find the veering triangulation of
the mapping torus $M$. Flipping the edges $P_1P_2$ and $PP_3$ yields the upper
right triangulation on \Cref{fig:translation-surface}. Then, flipping $P_2P_3$
and $PP_1$ yields the third triangulation on \Cref{fig:translation-surface}.
This triangulation is the image of the initial triangulation under $f$
(stretched horizontally by $\varphi^2$ and compressed vertically by
$\varphi^{-2}$). So the veering triangulation $\tau$ is obtained by gluing two
tetrahedra below the initial triangulation, then two tetrahedra under that, and
finally mapping the top (initial) triangulation to the bottom (final)
triangulation by $f$.

Therefore $\tau$ consists of four tetrahedra, four edges and eight
faces. The four edges are colored by blue, red, purple and green. Note
that the top and bottom edges of the tetrahedra and either blue and
red or green and purple. For each tetrahedra, the other four edges are
colored by four different colors.

\subsection*{The infinite cyclic cover of $S$}

The homology of $S$ is generated by the loops $c_1,c_2,c_3$ around the
punctures $P_1,P_2,P_3$, respectively. We have $f(c_1) = c_3$, $f(c_2) = c_1$
and $f(c_3) = c_2$. Therefore the $f$-invariant cohomology is
$H^1(S;\ZZ)^f = \langle \alpha \rangle$ where
$\alpha(c_1) = \alpha(c_2) = \alpha(c_3) = 1$. Let $t$ be a generator for
$H = \Hom(H^1(S;\ZZ)^f,\ZZ) \cong \ZZ$. By evaluating elements of
$H^1(S;\ZZ)^f$ on loops, we obtain a surjective homomorphism $\pi_1(S) \to H$.
Corresponding to this homomorphism is an infinite cyclic covering $\tS \to S$. For
more details about the theory, see Section 3 of \cite{McMullen00}.

\newcommand\TriangleOneB[9]{
    \def\x{1.618}
    \draw[thick,blue] (-\x,\x) -- node[below left] {$#1$} (0,0)  --
    node[below left] {$#2$} (\x,-\x);

    \draw[thick,violet] (-\x,\x) -- node[above] {$#3$}  ++(\x*\x,1)
    -- node[above] {$#4$} ++(\x*\x,1);

    \draw[thick,violet] (\x,-\x) -- node[right] {$#5$}
    ++(1,\x*\x) -- node[right] {$#6$}
    ++(1,\x*\x);

    \draw[thick,red] (0,0) -- node[below] {$#7$} (\x*\x,1);
    \draw[thick,red] (0,0) -- node[right] {$#8$} (1,\x*\x);
    \draw[thick,green] (\x*\x,1) -- node[above right] {$#9$} (1,\x*\x);
  }


\newcommand\TriangleTwoB[9]{
    \draw[thick,blue] (-\x*\x*\x,1/\x) -- node[below] {$#1$} (0,0) --
    node[below] {$#2$} (\x*\x*\x,-1/\x) ;
    \draw[thick,violet] (-\x*\x*\x,1/\x) -- node[above] {$#3$} ++(\x*\x,1)
     -- node[above] {$#4$} ++(\x*\x,1);
    \draw[thick,blue] (\x*\x*\x,-1/\x) -- node[right] {$#5$}
    ++(-\x,--\x) -- node[right] {$#6$} ++(-\x,--\x);
    \draw[thick,red] (0,0) -- node[right] {$#7$} (\x*\x,1);
    \draw[thick,green] (0,0) -- node[left] {$#8$} (-\x,--\x);
    \draw[thick,green] (\x*\x,1) --  node[above] {$#9$}(-\x,--\x);
}

\newcommand\TriangleThreeB[9]{
    \draw[thick,violet] (-\x*\x*\x*\x,-1/\x/\x) -- node[below] {$#1$}
    (0,0)  -- node[below] {$#2$}
    (--\x*\x*\x*\x,--1/\x/\x);
    \draw[thick,violet] (-\x*\x*\x*\x,-1/\x/\x) --  node[above] {$#3$} ++(\x*\x,1) -- node[above] {$#4$} ++(\x*\x,1);
    \draw[thick,blue] (--\x*\x*\x*\x,--1/\x/\x) -- node[above] {$#5$}
    ++(-\x*\x*\x,--1/\x) -- node[above] {$#6$} ++(-\x*\x*\x,--1/\x);
    \draw[thick,red] (0,0) -- node[right] {$#7$}(\x*\x,1);
    \draw[thick,green] (0,0) -- node[left] {$#8$} (-\x*\x*\x,--1/\x);
    \draw[thick,red] (\x*\x,1) -- node[above] {$#9$} (-\x*\x*\x,--1/\x);
}

\begin{figure}
  \centering
  \begin{tikzpicture}[scale=0.5]
    \def\x{1.618}
    \TriangleOneB{t^{-1}}{1}{tu}{u}{t^{-1}}{1}{t^{-1}u}{1}{1}
    \begin{scope}[yshift=-6cm]
      \TriangleOneB{1}{t}{t^2u}{tu}{1}{t}{u}{t}{t}
    \end{scope}

    \begin{scope}[yshift=-12cm]
      \TriangleOneB{t}{t^2}{t^3u}{t^2u}{t}{t^2}{tu}{t^2}{t^2}
    \end{scope}

    \begin{scope}[yshift=4cm,xshift=8cm]
      \TriangleTwoB{t^{-2}u}{t^{-1}u}{u}{t^{-1}u}{t^{-1}}{1}{t^{-1}u}{t^{-1}}{u}
    \end{scope}
    \begin{scope}[yshift=-2cm,xshift=8cm]
      \TriangleTwoB{t^{-1}u}{u}{tu}{u}{1}{t}{u}{1}{tu}
    \end{scope}

    \begin{scope}[yshift=-8cm,xshift=8cm]
      \TriangleTwoB{u}{tu}{t^2u}{tu}{t}{t^2}{tu}{t}{t^2u}
    \end{scope}

    \begin{scope}[yshift=-14cm,xshift=8cm]
      \TriangleTwoB{tu}{t^2u}{t^3u}{t^2u}{t^2}{t^3}{t^2u}{t^2}{t^3u}
    \end{scope}

    \begin{scope}[yshift=2cm,xshift=18cm]
      \TriangleThreeB{u^2}{tu^2}{u}{t^{-1}u}{t^{-1}u}{u}{u}{u}{t^{-1}u^2}
    \end{scope}

    \begin{scope}[yshift=-4cm,xshift=18cm]
      \TriangleThreeB{tu^2}{t^2u^2}{tu}{u}{u}{tu}{tu}{tu}{u^2}
    \end{scope}

    \begin{scope}[yshift=-10cm,xshift=18cm]
      \TriangleThreeB{t^2u^2}{t^3u^2}{t^2u}{tu}{tu}{t^2u}{t^2u}{t^2u}{tu^2}
    \end{scope}

  \end{tikzpicture}
  \caption{Part of the 2-skeleton of the veering triangulation on the
    maximal Abelian cover $\tM$ of $M$. The picture continues in all
    four directions indefinitely. The deck transformation $t$ acts by
    translating down by one triangle. The deck transformation $u$ acts
    by translating to the right by two columns, rotating each triangle
    by 180 degrees, stretching horizontally by $\varphi^2$ and
    vertically by $1/\varphi^2$.}
  \label{fig:abelian-cover}
\end{figure}

To construct $\tS$ explicitly, cut the upper left surface on
\Cref{fig:translation-surface} along the edges $PP_1$, $PP_2$ and $PP_3$, take
infinitely many copies of this cut-up surface, and reglue the edges according
to the labelling on the left column of \Cref{fig:abelian-cover} to obtain a
surface $\tS$. (Ignore the meaning of the labels for now, we will elaborate on
that later.) The action of the deck transformation $t$ is translating each
triangle to the triangle below it. To check that this is the right infinite
cyclic covering, all we need to check is that the loops $c_1$, $c_2$ and $c_3$,
oriented clockwise, all lift to paths in $\tS$ connecting some point $x$
to $tx$, therefore $c_1$, $c_2$ and $c_3$ all map to $t$ under the
homomorphism $\pi_1(S) \to H$.

\subsection*{The maximal Abelian cover of $M$}

To construct the maximal Abelian cover $\tM$ of $M$ and its veering
triangulation, we start with the triangulated surface in the left column of
\Cref{fig:abelian-cover} and we build down by gluing tetrahedra below it in the
same way as we did for the construction of the veering triangulation of $M$.

First, we glue tetrahedra to all quadrilaterals whose diagonals are lifts of
the edges $P_1P_2$ and $PP_3$. The bottom of the resulting cell complex is
triangulated as shown in the second column (after some rearranging of the
triangles to make the second column look similar to the first column). Then we
glue another round of tetrahedra to the bottom again obtain a cell complex
whose bottom is triangulated as shown in the third column.

Take infinitely many copies of this cell complex, indexed by $\ZZ$. Choose a
lift $\tilde{f}$ of $f$ identifying the triangulated surface in the first
column with the triangulated surface in the third column, and use this to glue
together the top of copy $i$ with the bottom of the copy $i+1$. There is an
isomorphism $u$ of the resulting cell complex that maps copy $i+1$ to copy $i$.
(One should think of copy $i+1$ to be above copy $i$ in the flow. Therefore $u$
shifts downward.) Our 3-manifold $M$ is the quotient of this cell complex by the
group generated by $t$ and $u$.


\subsection*{Labelling the edges by $G$}

The labelling of the edges on \Cref{fig:abelian-cover} can be found as
follows. First, for each color, label exactly one edge on the left
column by 1. These edges are the chosen lifts of the four edges of the
veering triangulation of $M$. For simplicity, we have chosen all
four lifts in the upper left triangle.

Using the $t$-action, the translates of the four edges in the middle
left and bottom left triangle should get the label $t$ and $t^2$,
respectively. Using the identification of the boundary edges, we can
label all edges in the left column except one red and two purple edges
in each triangle.

Now using the $u$-action, we can label all edges in the right column
except one red and two purple edges in each triangle.

The second column is obtained from the first column by flipping two
edges for each triangle, and the third column is obtained from the
second column analogously. This yields identifications between certain
edges in the first and second columns and also in the second and third
columns. In fact, the one purple and two red edges in each triangle in
the first column are present in the third column where they are
already labelled. Copying this labelling to the first column yields a
complete labelling of edges there.

Now using the $u$-action, we obtain a complete labelling of the third
column as well. Finally, using the identifications between the first
and second and the second and third columns, respectively, it is
possible to fully label the second column as well.

\subsection{The graph $\Delta^*$}
\label{sec:automaton-example}

\begin{figure}[htb]
\labellist
\small\hair 2pt
 \pinlabel {$t^{-1}u$} [ ] at 145 160
 \pinlabel {$t$} [ ] at 145 132
 \pinlabel {$1$} [ ] at 49 70
 \pinlabel {$t^{-1}u$} [ ] at 88 70
 \pinlabel {$u$} [ ] at 220 83
 \pinlabel {$t$} [ ] at 249 83
 \pinlabel {$t^2u$} [ ] at 150 -5
 \pinlabel {$t^{-2}$} [ ] at 147 25
 \pinlabel {$t^{-1}u$} [ ] at 2 121
 \pinlabel {$t^{-1}u^2$} [ ] at 33 105
 \pinlabel {$tu$} [ ] at 274 111
 \pinlabel {$tu^2$} [ ] at 294 126
\endlabellist
\centering
\includegraphics[scale=1.0]{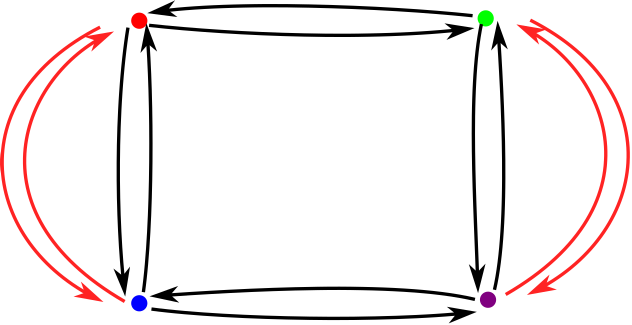}
\caption{The graph $\Delta^*$ corresponding to the fibered face containing the
  pseudo-Anosov braid $f=\sigma_1 \sigma_2^{-1}$. The black edges are the
  triangle-edges and the red edges are the tetrahedron-edges. Therefore the
  graph $\Delta$ is the subgraph consisting of black edges.}
\label{fig:delta}
\end{figure}
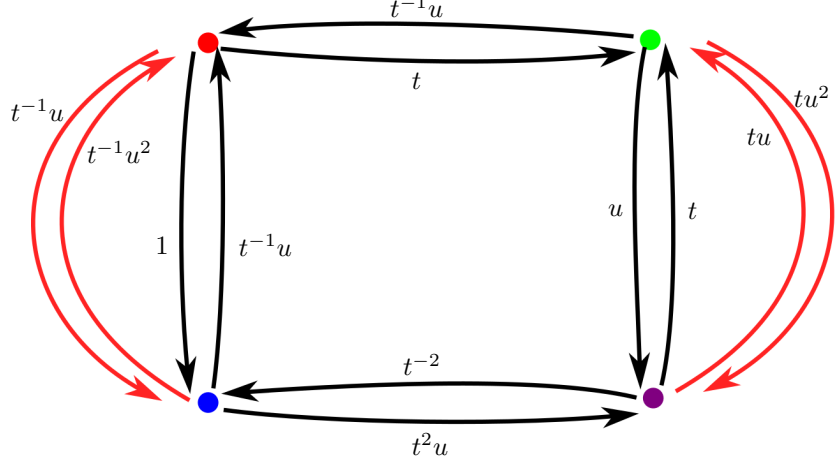

Using \Cref{fig:abelian-cover}, it is straightforward to construct the graph
$\Delta^*$ defined in \Cref{sec:delta-star}. See \Cref{fig:delta}.

\begin{table}
  \centering
  \begin{tabular}{l|cccccc}
    Cycle & RBR & RGR & BPB & GPG & RBPGR & RGPBR \\
    \hline
    Drift & $t^{-1}u$ & $u$ & $u$ & $tu$ & $t^2u^2$ & $t^{-2}u^2$ \\
  \end{tabular}
  \vspace{0.3cm}
  \caption{The drifts of minimal cycles.}
  \label{tab:minimal-cycles}
\end{table}

\subsection{Minimal cycles and minimal good paths}
\label{sec:minimal-path-example}

The minimal cycles are listed in \Cref{tab:minimal-cycles} with their
drift. Denoting the set of drifts of minimal cycles by $\calB$ as in
\Cref{sec:computing-intersection}, we have
\begin{displaymath}
  \calB = \{t^{-1}u, u, tu, t^{-2}u^2,t^2u^2\}
\end{displaymath}
and
\begin{displaymath}
  \langle \calB \rangle_{\RR_{\ge 0}} = \{t^a u^b : |a| \le b \}.
\end{displaymath}

The minimal good paths are listed in \Cref{tab:minimal-good-paths}.

\begin{table}
  \centering
  \begin{tabular}{l|ccccccc}
    Path & $RB$ & $RBR$ & $RBP$ & $RBRG$ & $RBPG$ & $RBRGP$ & $RBPGR$ \\
    \hline
    Drift & $t^{-1}u$ & $t^{-2}u^2$ & $t u^{2}$ & $t^{-1}u^2$ & $t^2u^2$ &
                                                                           $t^{-1}u^3$
                                                & $tu^3$     \vspace{0.2cm}
\\
    Path & $BR$ & $BRB$ & $BRG$ & $BRBP$ & $BRGP$ & $BRBPG$ & $BRGPB$ \\
    \hline
    Drift & $t^{-1}u^2$ & $t^{-1}u^2$ & $u^{2}$ & $t u^{3}$ & $u^{3}$ & $t^{2}
                                                                        u^{3}$
                                                            & $t^{-2}u^3$ \vspace{0.2cm}
    \\
    Path & $GP$ & $GPG$ & $GPB$ & $GPGR$ & $GPBR$ & $GPGRB$ & $GPBRG$ \\
    \hline
    Drift & $tu^2$ & $t^{2}u^2$ & $t^{-1}u^{2}$ & $tu^3$ & $t^{-2}u^3$ &
                                                                           $tu^3$
                                                & $t^{-1}u^3$     \vspace{0.2cm}
\\
    Path & $PG$ & $PGP$ & $PGR$ & $PGPB$ & $PGRB$ & $PGPBR$ & $PGRBP$ \\
    \hline
    Drift & $tu$ & $tu^2$ & $u^{2}$ & $t^{-1}u^{2}$ & $u^{2}$ & $t^{-2} u^{3}$ & $t^{2}u^3$ \\
  \end{tabular}
\vspace{0.2cm}
\caption{The drifts of minimal good paths.}
  \label{tab:minimal-good-paths}
\end{table}

\subsection{Determining the stashing sets $\stash(e,e')$}
\label{sec:determining-calI}

For a minimal good path $\gamma$, denote by $\stash(\gamma)$ the set of drifts of
good paths that decompose as the union of $\gamma$ with minimal cycles
(cf.~\Cref{lemma:decomposing-good-paths}).

\begin{proposition}\label{prop:drift-cones-of-paths}
  If $\gamma$ is a minimal good path in $\Delta^*$, then we have
  $\stash(\gamma) = \drift(\gamma) \cdot D$, where
  \begin{displaymath}
    D =
    \begin{cases}
      \langle \calB \rangle_{\ZZ_{\ge0}} - \{tu\} & \mbox{if }\gamma = RB, BR,
      RBR, BRB \\
      \langle \calB \rangle_{\ZZ_{\ge0}} - \{t^{-1}u\} & \mbox{if $\gamma = GP,
        PG, GPG, PGP$}\\
      \langle \calB \rangle_{\ZZ_{\ge0}} & \mbox{otherwise} \\
    \end{cases}
  \end{displaymath}
\end{proposition}
\begin{proof}
  If $\gamma$ consists of 3 or 4 edges, then it forms a connected collection
  with all minimal cycles. So in these cases, we have
  $D = \langle \calB \rangle_{\ZZ_{\ge0}}$.

  If $\gamma =RBP$, then the only minimal cycle $\gamma$ does not form a
  connected collection with is $RGR$. But the drift of the cycle $BPB$ is the
  same as the drift of the cycle $RGR$, so we still have
  $D = \langle \calB \rangle_{\ZZ_{\ge0}}$. We obtain
  $D = \langle \calB \rangle_{\ZZ_{\ge0}}$ similarly for $\gamma=BRG$, $GPB$
  and $PGB$.

  For the remaining four possibilities for $\gamma$ starting with $R$ or $B$
  ($RB$, $BR$, $RBR$, $BRB$), the cycles $RBR$, $RBPGR$ and $RGPBR$ form a
  connected collection with $\gamma$, so $t^{-1}u, t^2u^2, t^{-2}u^2 \in D$.
  Also, at least one of $RGR$ and $BPB$ forms a connected collection with
  $\gamma$, so $u \in D$. However, the cycle $GPG$ does not form a connected
  collection with $\gamma$, so $tu \notin D$.

  It remains to show that all elements of
  $\langle \calB \rangle_{\ZZ_{\ge0}} - \{tu\}$ that are not nonnegative
  integral linear combinations of $u, t^{-1}u$ and $t^2u^2$ are in $D$. This
  follows from the fact that once we add to $\gamma$ a minimal cycle with drift
  $u$ or $t^2u^2$, we can now add the cycle $GPG$ to obtain a connected
  collection. So the translated cones $u\langle \calB \rangle_{\ZZ_{\ge0}}$ and
  $t^2u^2\langle \calB \rangle_{\ZZ_{\ge0}}$ are contained in $D$. This
  completes the proof in the case that $\gamma$ starts with $R$ or $B$.

  The proof is analogous in case $\gamma$ starts with $G$ or $P$.
\end{proof}

\newcommand\dotbase{
  \def\rad{0.1}
  \draw (-5,0) -- (5,0);
  \draw (0,0) -- (0,5);

  \foreach \i in {0,1,2,3,4,5}{
    \foreach \j in {-\i,...,\i}
    \draw[thin,gray] (\j,\i) circle (\rad);
  }
}

\newcommand\rpoint[2]{
  \draw (#1-0.2,#2-0.2) -- (#1+0.2,#2+0.2);
}

\newcommand\lpoint[2]{
  \draw (#1-0.2,#2+0.2) -- (#1+0.2,#2-0.2);
}

\newcommand\lrpoint[2]{
  \rpoint{#1}{#2}
  \lpoint{#1}{#2}
}

\begin{figure}\centering

\newcommand\picone{  \dotbase
  \foreach \i in {-5,...,3}
  \filldraw (\i,5) circle (\rad);
  \foreach \i in {-4,...,2}
  \filldraw (\i,4) circle (\rad);
  \foreach \i in {-3,-2,1}
  \filldraw (\i,3) circle (\rad);
  \foreach \i in {-2}
  \filldraw (\i,2) circle (\rad);
  \lpoint{-2}{2}
  \lrpoint{1}{3}
}

\begin{tikzpicture}[scale=0.4]
  \picone
  \node[above] at (4,0) {\small $\stash(R, R)$};
\end{tikzpicture}\hspace{5mm}
\begin{tikzpicture}[scale=0.4,xscale=-1]
  \picone
  \node[above] at (4,0) {\small $\stash(G, G)$};
\end{tikzpicture}\hspace{5mm}
\newcommand\pictwo{
  \dotbase
  \foreach \i in {-5,...,3}
  \filldraw (\i,5) circle (\rad);
  \foreach \i in {-4,...,2}
  \filldraw (\i,4) circle (\rad);
  \foreach \i in {-3,...,1}
  \filldraw (\i,3) circle (\rad);
  \foreach \i in {-2,-1}
  \filldraw (\i,2) circle (\rad);
  \filldraw (-1,1) circle (\rad);
  \lpoint{-1}{1}
}
\begin{tikzpicture}[scale=0.4]
  \pictwo
  \node[above] at (4,0) {\small $\stash(R, B)$};
\end{tikzpicture}\hspace{5mm}
\begin{tikzpicture}[scale=0.4,xscale=-1]
  \dotbase
  \foreach \i in {-4,...,2}
  \filldraw (\i,5) circle (\rad);
  \foreach \i in {-3,...,1}
  \filldraw (\i,4) circle (\rad);
  \foreach \i in {-2,...,-1}
  \filldraw (\i,3) circle (\rad);
  \filldraw (-1,2) circle (\rad);
  \lpoint{-1}{2}
  \node[above] at (4,0) {\small $\stash(G,P)$};
\end{tikzpicture}\hspace{5mm}
\newcommand{\picthree}{
  \dotbase
  \foreach \i in {-4,...,5}
  \filldraw (\i,5) circle (\rad);
  \foreach \i in {-3,...,4}
  \filldraw (\i,4) circle (\rad);
  \foreach \i in {-2,...,3}
  \filldraw (\i,3) circle (\rad);
  \foreach \i in {-1,2}
  \filldraw (\i,2) circle (\rad);
  \lrpoint{2}{2}
  \lrpoint{-1}{2}
}
\begin{tikzpicture}[scale=0.4]
  \picthree
  \node[above] at (4,0) {\small $\stash(R, G)$};
\end{tikzpicture}\hspace{5mm}
\begin{tikzpicture}[scale=0.4,xscale=-1]
  \dotbase
  \foreach \i in {-3,...,4}
  \filldraw (\i,5) circle (\rad);
  \foreach \i in {-2,...,3}
  \filldraw (\i,4) circle (\rad);
  \foreach \i in {-1,2}
  \filldraw (\i,3) circle (\rad);
  \lrpoint{2}{3}
  \lrpoint{-1}{3}
  \node[above] at (4,0) {\small $\stash(G, R)$};
\end{tikzpicture}\hspace{5mm}
\newcommand{\picfour}{  \dotbase
  \foreach \i in {-3,...,4}
  \filldraw (\i,5) circle (\rad);
  \foreach \i in {-2,...,3}
  \filldraw (\i,4) circle (\rad);
  \foreach \i in {-1,...,2}
  \filldraw (\i,3) circle (\rad);
  \foreach \i in {1}
  \filldraw (\i,2) circle (\rad);
  \lrpoint{1}{2}
  \lrpoint{-1}{3}
}
\begin{tikzpicture}[scale=0.4]
  \picfour
  \node[above] at (4,0) {\small $\stash(R, P)$};
\end{tikzpicture}\hspace{5mm}
\begin{tikzpicture}[scale=0.4,xscale=-1]
  \picfour
  \node[above] at (4,0) {\small $\stash(G, B)$};
\end{tikzpicture}\hspace{5mm}
\newcommand{\picfive}{
  \dotbase
  \foreach \i in {-4,...,2}
  \filldraw (\i,5) circle (\rad);
  \foreach \i in {-3,...,1}
  \filldraw (\i,4) circle (\rad);
  \foreach \i in {-2,-1}
  \filldraw (\i,3) circle (\rad);
  \foreach \i in {-1}
  \filldraw (\i,2) circle (\rad);
  \lpoint{-1}{2}
}
\begin{tikzpicture}[scale=0.4]
  \picfive
  \node[above] at (4,0) {\small $\stash(B, R)$};
\end{tikzpicture}\hspace{5mm}
\begin{tikzpicture}[scale=0.4,xscale=-1]
  \pictwo
  \node[above] at (4,0) {\small $\stash(P, G)$};
\end{tikzpicture}\hspace{5mm}
\newcommand{\picsix}{
  \dotbase
  \foreach \i in {-4,...,2}
  \filldraw (\i,5) circle (\rad);
  \foreach \i in {-3,...,1}
  \filldraw (\i,4) circle (\rad);
  \foreach \i in {-2,-1}
  \filldraw (\i,3) circle (\rad);
  \foreach \i in {-1}
  \filldraw (\i,2) circle (\rad);
  \lpoint{-1}{2}
  \lrpoint{-2}{3}
}
\begin{tikzpicture}[scale=0.4]
  \picsix
  \node[above] at (4,0) {\small $\stash(B, B)$};
\end{tikzpicture}\hspace{5mm}
\begin{tikzpicture}[scale=0.4,xscale=-1]
  \picsix
  \node[above] at (4,0) {\small $\stash(P, P)$};
\end{tikzpicture}\hspace{5mm}
\newcommand{\picseven}{
  \dotbase
  \foreach \i in {-2,...,3}
  \filldraw (\i,5) circle (\rad);
  \foreach \i in {-1,...,2}
  \filldraw (\i,4) circle (\rad);
  \foreach \i in {0,1}
  \filldraw (\i,3) circle (\rad);
  \lrpoint{0}{3}
  \lrpoint{1}{3}
}
\begin{tikzpicture}[scale=0.4]
  \picseven
  \node[above] at (4,0) {\small $\stash(B, P)$};
\end{tikzpicture}\hspace{5mm}
\begin{tikzpicture}[scale=0.4,xscale=-1]
  \dotbase
  \foreach \i in {-3,...,4}
  \filldraw (\i,5) circle (\rad);
  \foreach \i in {-2,...,3}
  \filldraw (\i,4) circle (\rad);
  \foreach \i in {-1,...,2}
  \filldraw (\i,3) circle (\rad);
  \foreach \i in {0,1}
  \filldraw (\i,2) circle (\rad);
  \lrpoint{0}{2}
  \lrpoint{1}{2}
  \node[above] at (4,0) {\small $\stash(P, B)$};
\end{tikzpicture}\hspace{5mm}
\newcommand{\piceight}{
  \dotbase
  \foreach \i in {-3,...,4}
  \filldraw (\i,5) circle (\rad);
  \foreach \i in {-2,...,3}
  \filldraw (\i,4) circle (\rad);
  \foreach \i in {-1,...,2}
  \filldraw (\i,3) circle (\rad);
  \foreach \i in {0}
  \filldraw (\i,2) circle (\rad);
  \lrpoint{0}{2}
  \lrpoint{2}{3}
}\hspace{5mm}
\begin{tikzpicture}[scale=0.4]
  \piceight
  \node[above] at (4,0) {\small $\stash(B, G)$};
\end{tikzpicture}\hspace{5mm}
\begin{tikzpicture}[scale=0.4,xscale=-1]
  \piceight
  \node[above] at (4,0) {\small $\stash(P, R)$};
\end{tikzpicture}\hspace{0mm}

\caption{The sets $\stash(e,e')$ for all pairs $e,e' \in \{R, B, G, P\}$. }
\label{fig:cone-pictures}
\end{figure}

\Cref{prop:drift-cones-of-paths} allows us to determine the stashing set $\stash(e,e')$ for every
pair $e, e'\in \{B, R, G, P\}$ using the formula
\begin{displaymath}
  \stash(e, e') = \bigcup\{\stash(\gamma): \gamma \mbox{ is a minimal good path
    from  $e$ to  $e'$}\}.
\end{displaymath}
For example,
\begin{align*}
  \stash(R, R) &= \stash(RBR) \cup \stash(RBPGR) \\ &= t^{-2}u^2(\langle \calB
  \rangle_{\ZZ_{\ge0}} - \{tu\}) \cup tu^3 \langle \calB \rangle_{\ZZ_{\ge0}}
\end{align*}
since $\drift(RBR) = t^{-2}u^2$ and $\drift(RBPGR) = tu^3$. We can visualize
that computation as follows.

Draw the cone $\langle t^{-1}u, u, tu \rangle_{\ZZ_{\ge0}}$ as on
\Cref{fig:cone-pictures}. The horizontal axis is the $t$-axis and the vertical
axis is the $u$-axis. Mark the point $t^{-2}u^2$ with a ``left tick'' and the
point $tu^3$ with a cross, indicating the coloring of the sets
$t^{-2}u^2(\langle \calB \rangle_{\ZZ_{\ge0}} - \{tu\})$ and
$tu^3\langle \calB \rangle_{\ZZ_{\ge0}}$. The union of these two sets forms
$\stash(R, R)$. \Cref{fig:cone-pictures} illustrates the computation for all
$\stash(e, e')$.

\subsection{The fibered face}
\label{sec:asdf}

Denote by $\phi_{a,b}$ the element of $H^1(M)$ such that
$\phi_{a,b}(t) = a$ and $\phi_{a,b}(u) = b$.

\begin{lemma}
  The monodromy $f = \sigma_1 \sigma_2^{-1}$ corresponds to $\phi_{0,-1}$.
\end{lemma}
\begin{proof}
  The homology class $t$ can be represented by loops in $M$ that are in the
  fiber dual to $f$ and loops representing $u$ wind around the fibration dual to
  $f$ once in the opposite direction of the flow.
\end{proof}

Let $\calF$ be the fibered face such that the cone $\RR_+\calF$ contains $\phi_{0,-1}$.

\begin{lemma}
  $\RR_+\calF = \{\phi_{a,b} : |a| \le -b, b < 0\}$.
\end{lemma}
\begin{proof}
  By \Cref{prop:cone-description}, the cone $\RR_+\calF$ contains precisely
  those cohomology classes that take nonpositive values on
  $\langle \calB \rangle_{\RR_{\ge0}} = \langle t^{-1}u, tu \rangle_{\RR_{\ge0}}$.
\end{proof}

\begin{lemma}\label{lemma:norm-of-f}
  $\lVert \phi_{0,-1} \rVert = 2$.
\end{lemma}
\begin{proof}
  The fiber dual to $f$ is a four times punctured sphere and has Euler
  characteristic $-2$.
\end{proof}

\begin{lemma}\label{lemma:norm}
  We have $\lVert \phi_{a,b} \rVert = -2b$ whenever $\phi_{a,b} \in \RR_+\calF$.
\end{lemma}
\begin{proof}
  This follows from \Cref{lemma:norm-of-f} and the fact that the Thurston norm
  on the cone $\RR_+\calF$ has the symmetry
  $\lVert \phi_{a,b} \rVert = \lVert \phi_{-a,b} \rVert$. This symmetry can be
  seen from the following symmetry $\delta$ of the veering triangulation of
  $M$. The symmetry $\delta$ maps the upper left triangulation on
  \Cref{fig:translation-surface} to the upper right triangulation by vertical
  reflection and a horizontal stretch. This map extends to a symmetry of the
  whole veering triangulation and exchanges the red edge with the green edge
  and the blue edge with the purple edge. The symmetry $\delta: M \to M$ maps
  the fiber $S$ to a homotopic fiber by an orientation-reversing map, but
  preserves the orientation of the flow. Therefore the action on homology is
  $\delta_*(u) = u$ and $\delta_*(t) = -t$. So in the cone $\RR_+\calF$, two
  integral points $\phi_{-a,b}$ and $\phi_{a,b}$ are dual to homeomorphic
  fibers and hence their Thurston norms are equal.
\end{proof}

\begin{corollary}
  $\calF = \{\phi_{a,-\frac12} : |a| \le \frac12\}$
\end{corollary}

\subsection{Accumulation points of the graph of $\mu_1$}
\label{sec:mu-two-example}

We will apply \Cref{thm:simplex} to the fibered face $\calF$ to find the
accumulation points of the graph of the function $\mu_1$. According to the
theorem, the set of accumulation points is the graph of a continuous function
$g: \inter(\calF) \to \RR_+$.
\begin{proposition}\label{prop:g-example}
\begin{displaymath}
  g\left(\phi_{a,-\frac12}\right) = \frac{2}{(\frac12 -a)(\frac12 + a)}.
\end{displaymath}
\end{proposition}
\begin{proof}
  By \Cref{thm:simplex}, we have
  \begin{displaymath}
    g(\alpha\omega_1 + (1-\alpha)\omega_2) =
    g^*(\alpha, 1-\alpha) = \frac{1}{
      \vol_\Lambda(\Sigma/\langle \omega_1,
      \omega_{2}\rangle_\ZZ) \cdot \alpha(1-\alpha)}
  \end{displaymath}
  where $\omega_1 = \phi_{-\frac12,-\frac12}$ and
  $\omega_2 = \phi_{\frac12,-\frac12}$, $\Sigma = H^1(M;\RR)$ and
  $\Lambda = H^1(M;\ZZ)$. Since
  \begin{displaymath}
    \det\begin{pmatrix}
      -\frac12 & \frac12 \\
      -\frac12 & -\frac12 \\
    \end{pmatrix} = \frac12,
  \end{displaymath}
  the covolume in the denominator is $\frac12$. So
  \begin{displaymath}
    g\left(\phi_{-\frac{\alpha}2 + \frac{1-\alpha}{2}, -\frac12}\right) = \frac{2}{\alpha(1-\alpha)}
  \end{displaymath}
  and by substituting $a = \frac12 - \alpha$, we obtain the desired formula.
\end{proof}

\subsection{Exact values of the function $\mu_1$}
\label{sec:exact-values}

\begin{proposition}\label{prop:translation-lengths}
  For the primitive integral cohomology classes $\phi$ listed in the table
  below, the asymptotic translation length $\ell_\calA(\phi)$ of the monodromy
  corresponding to $\phi$ is as shown in the table.

  \begin{center}
    \begin{tabular}{lll}
      $\phi$ & $\ell_\calA(\phi)$ & maximal cycles \\
      \hline
      $\phi_{0,-1}$ & $\frac{2}3$ & BPB, RGR\\
      $\phi_{1,-2}$ & $\frac{1}6$ & BB \\
      $\phi_{1,-3}$ & $\frac19$ & BB, RR\\
      $\phi_{1,-4}$ & $\frac{1}{13}$ & RR \\
      $\phi_{1,-k}$ ($k\ge 5$ odd) & $\frac{2}{(k+1)^2}$ & BB \\
      $\phi_{1,-k}$ ($k\ge 6$ even) & $\frac2{k^2+2k-1}$ & BPB \\
    \end{tabular}
  \end{center}
    For each fibration, the table also shows the cycles of $\Delta^*$ with
    maximal average weight.
  \end{proposition}
  We remark that these cycles correspond to bi-infinite geodesics in the arc
  complex of the fiber that are invariant under some power of the monodromy.

  \begin{proof}
    Using the stashing sets $\stash(e,e')$ shown on \Cref{fig:cone-pictures}, we can
    determine the weighted graphs $W(\phi)$ using the definition
  \Cref{eq:weight}.

  \begin{proofcase}[$\phi_{0,-1}$, $\phi_{1,-2}$, $\phi_{1,-3}$ and
    $\phi_{1,-4}$]
    In these cases, we find the graphs $W(\phi)$ shown on
    \Cref{fig:cycle_graph_examples} by a case-by-case inspection of each set on
    \Cref{fig:cone-pictures}. The maximum averages are $\frac{3}2$, $6$, 9 and
    13, therefore the asymptotic translation lengths are $\frac{2}3$,
    $\frac16$, $\frac19$ and $\frac1{13}$, respectively, by
    \Cref{prop:atl-via-weighted-graphs}.
  \end{proofcase}

\begin{figure}[htb]
      \labellist
      \small\hair 2pt
      \pinlabel {$1$} [ ] at 32 1
      \pinlabel {$1$} [ ] at 34 197
      \pinlabel {$1$} [ ] at 200 199
      \pinlabel {$1$} [ ] at 202 9
      \pinlabel {$0$} [ ] at 15 106
      \pinlabel {$1$} [ ] at 48 106
      \pinlabel {$2$} [ ] at 118 187
      \pinlabel {$1$} [ ] at 118 154
      \pinlabel {$1$} [ ] at 185 106
      \pinlabel {$0$} [ ] at 217 106
      \pinlabel {$2$} [ ] at 118 18
      \pinlabel {$1$} [ ] at 117 50
      \pinlabel {$1$} [ ] at 142 140
      \pinlabel {$1$} [ ] at 160 117
      \pinlabel {$1$} [ ] at 76 117
      \pinlabel {$1$} [ ] at 95 140
      \endlabellist
      \centering
      \includegraphics[scale=0.6]{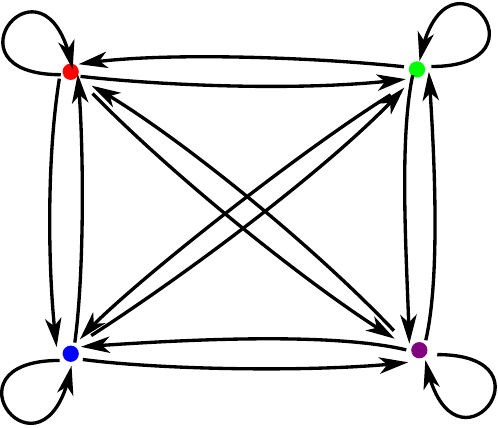}
      \quad
  \labellist
  \small\hair 2pt
  \pinlabel {$6$} [ ] at 32 1
  \pinlabel {$4$} [ ] at 34 197
  \pinlabel {$1$} [ ] at 200 199
  \pinlabel {$2$} [ ] at 202 9
  \pinlabel {$4$} [ ] at 15 106
  \pinlabel {$6$} [ ] at 48 106
  \pinlabel {$4$} [ ] at 118 187
  \pinlabel {$1$} [ ] at 118 154
  \pinlabel {$2$} [ ] at 185 106
  \pinlabel {$0$} [ ] at 217 106
  \pinlabel {$4$} [ ] at 118 18
  \pinlabel {$3$} [ ] at 117 50
  \pinlabel {$4$} [ ] at 142 140
  \pinlabel {$3$} [ ] at 160 117
  \pinlabel {$2$} [ ] at 76 117
  \pinlabel {$3$} [ ] at 95 140
    \endlabellist
    \centering
    \includegraphics[scale=0.6]{cycle_graph}\\
    \vspace{0.5cm}
  \labellist
  \small\hair 2pt
 \pinlabel {$9$} [ ] at 32 1
 \pinlabel {$9$} [ ] at 34 197
 \pinlabel {$5$} [ ] at 200 199
 \pinlabel {$6$} [ ] at 202 9
 \pinlabel {$6$} [ ] at 15 106
 \pinlabel {$9$} [ ] at 48 106
 \pinlabel {$9$} [ ] at 118 187
 \pinlabel {$5$} [ ] at 118 154
 \pinlabel {$6$} [ ] at 185 106
 \pinlabel {$3$} [ ] at 217 106
 \pinlabel {$7$} [ ] at 118 18
 \pinlabel {$5$} [ ] at 117 50
 \pinlabel {$6$} [ ] at 142 140
 \pinlabel {$5$} [ ] at 160 117
 \pinlabel {$6$} [ ] at 76 117
 \pinlabel {$7$} [ ] at 95 140
    \endlabellist
    \centering
    \includegraphics[scale=0.6]{cycle_graph}
    \quad
    \labellist
    \small\hair 2pt
    \pinlabel {$12$} [ ] at 38 5
    \pinlabel {$13$} [ ] at 36 197
    \pinlabel {$11$} [ ] at 198 199
    \pinlabel {$12$} [ ] at 198 9
    \pinlabel {$8$} [ ] at 15 106
    \pinlabel {$12$} [ ] at 48 106
    \pinlabel {$13$} [ ] at 118 187
    \pinlabel {$8$} [ ] at 118 154
    \pinlabel {$12$} [ ] at 185 106
    \pinlabel {$8$} [ ] at 217 106
    \pinlabel {$13$} [ ] at 118 18
    \pinlabel {$10$} [ ] at 117 50
    \pinlabel {$10$} [ ] at 142 140
    \pinlabel {$9$} [ ] at 160 117
    \pinlabel {$9$} [ ] at 76 117
    \pinlabel {$10$} [ ] at 95 140
    \endlabellist
    \centering
    \includegraphics[scale=0.6]{cycle_graph}
    \caption{The graphs $W(\phi_{0,-1})$ (top left), $W(\phi_{1,-2})$
      (top right),
      $W(\phi_{1,-3})$ (bottom left) and $W(\phi_{1,-4})$ (bottom
      right). }
    \label{fig:cycle_graph_examples}
  \end{figure}

  \begin{proofcase}[$k\ge 5$ is odd]
    Note that $t^{\frac{k-1}2}u^{\frac{k+3}2}$ and
    $t^{-\frac{k+1}2}u^{\frac{k+1}2}$ are not in $\stash(B,B)$, and both
    evaluate to $-\frac{(k+1)^2}2$ by $\phi_{1,-k}$. Hence
    $\frac{(k+1)^2}2 \notin -\phi_{1,-k}(\stash(B,B))$. One easily verifies that
    $\frac{(k+1)^2}2$ is in fact the largest integer that is not contained in
    $-\phi_{1,-k}(\stash(B,B))$. Since $\stash(B,R) = \stash(B,B)$, it is also the
    largest integer not contained in $-\phi_{1,-k}(\stash(B,R))$.

    The set
    \begin{displaymath}
    Z_k = \left\{t^au^b: -\frac{k-1}2 \le a \le \frac{k-1}2, b \ge \frac{k+3}2\right\}
  \end{displaymath}
is
    contained in $\stash(R,R)$, $\stash(R,B)$, $\stash(R,G)$, $\stash(R,P)$,
    $\stash(B,G)$, $\stash(G,G)$, $\stash(G,R)$, $\stash(G,B)$, $\stash(P,G)$,
    $\stash(P,B)$, $\stash(P,R)$, so the weights of the corresponding edges are
    less than
    \begin{displaymath}
      k\frac{k+3}{2} - \frac{k-1}{2} = \frac{k^2+2k+1}2 = \frac{(k+1)^2}2.
    \end{displaymath}
    The set $tZ_k$ is contained in $\stash(G,P) = \stash(P, P)$, so the weights of
    the corresponding edges are strictly less than $\frac{(k+1)^2}2-1$. For the
    only remaining pair, $(B,P)$, one can check that $w(B,P) = \frac{(k+1)^2}2-1$.

    Therefore $w(B, B) = w(B, R) = \frac{(k+1)^2}2$ and the weights of the
    other edges are strictly smaller. Therefore the largest average cycle
    weight is $\frac{(k+1)^2}2$, realized only by the loop on the blue vertex.
  \end{proofcase}

  \begin{proofcase}[$k\ge 6$ is even]
    The set
    \begin{displaymath}
      W_k = \left\{t^au^b: -\frac{k}2 \le a \le \frac{k-2}2, b \ge \frac{k+2}2\right\} -
      \left\{t^{\frac{k-2}2}u^{\frac{k+2}2}\right\}
  \end{displaymath}
is contained in $\stash(R,R)$,
    $\stash(R,B)$, $\stash(R, G)$, $\stash(B,R)$, $\stash(B, B)$, $\stash(G, R)$,
    $\stash(G, B)$, therefore the weights of the corresponding edges are at most
    \begin{displaymath}
      X_k = k\frac{k+2}2 - \frac{k-2}2 =\frac{k^2+k+2}2.
    \end{displaymath}

    The set $tW_k$ is contained in $\stash(R,P)$, $\stash(B, G)$, $\stash(G,G)$,
    $\stash(P, G)$, $\stash(P,B)$, $\stash(P, R)$, therefore the weights of the
    corresponding edges are at most $X_k-1$.

    The set $t^2W_k$ is contained in $\stash(G,P)$, $\stash(P, P)$, therefore the
    weights of the corresponding edges are at most $X_k-2$.

    There is one remaining edge: $BP$. Using that neither
    $t^{-\frac{k-2}2}u^{\frac{k+2}2}$, nor $t^{\frac{k+2}2}u^{\frac{k+4}2}$ are
    in $\stash(B,P)$ and both expressions evaluate to $\frac{k^2+3k-2}2$, one
    can verify that $w(B, P) = \frac{k^2+3k-2}2$. One can also check that
    $w(P,B)$ is in fact exactly $X_k-1$, therefore the cycle $BPB$ has average
    weight $\frac{k^2+2k-1}2$. Since this weight is larger than $X_k$, no other
    cycle can have the same of larger average weight.
  \end{proofcase}
\end{proof}

Finally, we give the proof of \Cref{thm:example}.

\begin{proof}[Proof of \Cref{thm:example}]
  Parametrize the fibered face $\calF$ with the interval $[-1,1]$, using the
  map $\phi_{a,b} \mapsto \frac{a}{b}$.

  Using \Cref{lemma:norm}, we have $\|\phi_{0,-1}\|^2 = 4$ and
  $\|\phi_{\pm 1,-k}\|^2 = 4k^2$ for every integer $k \ge 2$. Together with
  \Cref{prop:translation-lengths}, we obtain the values of $\mu_1(t)$ in
  \Cref{thm:example} for $t\le 0$. By the symmetry discussed in the proof of
  \Cref{lemma:norm}, we have $\mu_1(t) = \mu_1(-t)$ for all $t\in (-1,1)$,
  which yields the claimed values of $\mu_1(t)$ when $t > 0$.

  Using the substitution $t = \frac{a}{-\frac12}$, hence $a=-\frac12 t$, in
  \Cref{prop:g-example}, we obtain that the set of accumulation points of the
  graph of $\mu_1(t)$ is the graph of
  \begin{displaymath}
    \frac{2}{(\frac12+\frac12 t)(\frac12-\frac12 t)} = \frac{8}{1-t^2}
  \end{displaymath}
  as claimed. Finally, it is straightforward to check that $\mu_1(t) <
  \frac{8}{1-t^2}$ for all values of $t$ for which we have determined $\mu_1(t)$.
\end{proof}

\bibliographystyle{alpha}
\bibliography{mybibfile}

\end{document}